\documentclass{article}

\usepackage[sc]{mathpazo}

\usepackage{amsmath, amssymb, amsthm, mathrsfs, enumitem, tikz, stmaryrd, datetime, pifont, comment, xcolor}
\usepackage[colorlinks=true, allcolors=violet]{hyperref}
\usepackage[all]{xy}
\usepackage[section]{placeins}

\renewcommand{\today}{\the\day\ \monthname[\the\month]\ \the\year}

\renewcommand{\AA}{\mathbb{A}}
\newcommand{\cA}{\mathcal{A}}

\newcommand{\sA}{\mathscr{A}}
\newcommand{\CC}{\mathbb{C}}

\newcommand{\DD}{\mathbb{D}}
\newcommand{\fD}{\mathfrak{D}}
\newcommand{\fE}{\mathfrak{E}}

\newcommand{\FF}{\mathbb{F}}

\newcommand{\sF}{\mathscr{F}}

\newcommand{\GG}{\mathbb{G}}

\newcommand{\fG}{\mathfrak{G}}
\newcommand{\sG}{\mathscr{G}}

\newcommand{\cH}{\check{H}}

\newcommand{\sI}{\mathscr{I}}

\newcommand{\sJ}{\mathscr{J}}

\newcommand{\sL}{\mathscr{L}}

\newcommand{\OO}{\mathcal{O}}

\newcommand{\QQ}{\mathbb{Q}}
\newcommand{\fQ}{\mathfrak{Q}}
\newcommand{\RR}{\mathbb{R}}
\newcommand{\bbS}{\mathbb{S}}
\newcommand{\sS}{\mathscr{S}}

\newcommand{\VV}{\mathbb{V}}

\newcommand{\ZZ}{\mathbb{Z}}

\newcommand{\bfb}{\mathbf{b}}
\newcommand{\fc}{\mathfrak{c}}

\newcommand{\fk}{\mathfrak{k}}

\newcommand{\bft}{\mathbf{t}}

\newcommand{\Hom}{\operatorname{Hom}}
\newcommand{\sHom}{\mathscr{H}\kern -1.5pt om}

\newcommand{\im}{\operatorname{im}}

\newcommand{\GL}{\operatorname{GL}}

\newcommand{\GSp}{\operatorname{GSp}}
\newcommand{\tr}{\operatorname{tr}}

\newcommand{\id}{\operatorname{id}}
\newcommand{\End}{\operatorname{End}}
\newcommand{\Aut}{\operatorname{Aut}}

\newcommand{\Res}{\operatorname{Res}}

\newcommand{\coker}{\operatorname{coker}}

\newcommand{\supp}{\operatorname{supp}}

\newcommand{\Gal}{\operatorname{Gal}}

\renewcommand{\mod}{\operatorname{mod}}

\newcommand{\Frob}{\operatorname{Frob}}

\newcommand{\Ig}{{\operatorname{Ig}}}
\newcommand{\vol}{\operatorname{vol}}
\newcommand{\Rep}{\operatorname{Rep}}
\newcommand{\Sh}{{\operatorname{Sh}}}
\newcommand{\Isom}{\operatorname{Isom}}

\newcommand{\Int}{\operatorname{Int}}
\newcommand{\Fix}{\operatorname{Fix}}
\newcommand{\Groth}{\operatorname{Groth}}

\newcommand{\Isoc}{\operatorname{Isoc}}

\newcommand{\bpm}{\begin{pmatrix}}
\newcommand{\epm}{\end{pmatrix}}
\newcommand{\bbm}{\begin{bmatrix}}
\newcommand{\ebm}{\end{bmatrix}}

\newcommand{\ab}{{\rm ab}}
\newcommand{\ad}{{\rm ad}}

\newcommand{\ur}{{\rm ur}}

\newcommand{\et}{{\text{\'et}}}
\newcommand{\crys}{\text{crys}}

\newcommand{\one}{\mathbb{1}}

\newcommand{\dlim}{\varinjlim}
\newcommand{\ilim}{\varprojlim}

\DeclareFontFamily{U}{wncy}{}
\DeclareFontShape{U}{wncy}{m}{n}{<->wncyr10}{}
\DeclareSymbolFont{mcy}{U}{wncy}{m}{n}
\DeclareMathSymbol{\Sha}{\mathord}{mcy}{"58}

\usepackage{titlesec}
\usepackage{xcolor}
\usepackage{chngcntr}

\titleformat{\subsubsection}
   [runin]{\normalfont\normalsize\bfseries}{\thesubsubsection}{4pt}{}
\renewcommand{\paragraph}{\subsubsection}

\makeatletter
\let\c@equation\c@subsubsection

\makeatother

\newtheorem{theorem}[subsubsection]{Theorem}
\newtheorem{proposition}[subsubsection]{Proposition}
\newtheorem{lemma}[subsubsection]{Lemma}

\newtheorem*{theorem*}{Theorem}
\newtheorem*{proposition*}{Proposition}
\newtheorem*{lemma*}{Lemma}
\newtheorem*{corollary*}{Corollary}
\newtheorem*{conjecture*}{Conjecture}
\newtheorem*{claim*}{Claim}

\theoremstyle{definition}
\newtheorem{definition}[subsubsection]{Definition}

\newtheorem*{definition*}{Definition}
\newtheorem*{example*}{Example}
\newtheorem*{exercise*}{Exercise}
\newtheorem*{fact*}{Fact}
\newtheorem*{remark*}{Remark}
\newtheorem*{problem*}{Problem}
\newtheorem*{question*}{Question}

\newcommand{\LRP}{\mathcal{LRP}} 
\newcommand{\KP}{\mathcal{KP}} 
\newcommand{\sa}{{\rm sa}} 
\newcommand{\ba}{{\bfb\text{-adm}}} 
\newcommand{\der}{{\rm der}} 
\renewcommand{\gg}{{\rm gg}} 
\newcommand{\bas}{{\rm bas}} 
\newcommand{\osS}{\overline{\mathscr{S}}} 
\renewcommand{\cH}{\mathcal{H}} 
\newcommand{\AM}{\mathcal{AM}} 
\newcommand{\Comp}{\operatorname{Comp}} 
\newcommand{\Rell}{\RR-{\rm ell}} 
\newcommand{\Fl}{\sF\!\ell} 
\newcommand{\oFFp}{\overline{\FF}_p}
\newcommand{\fr}{{fr}}

\title{Counting Points on Igusa Varieties of Hodge Type}
\author{Sander Mack-Crane}
\date{\today}

\begin{document}
\maketitle

\begin{abstract}
Igusa varieties are algebraic varieties that arise in the study of special fibers of Shimura varieties, and have demonstrated many applications in the Langlands program via a Langlands--Kottwitz style point-counting formula due to Shin \cite{MR2484281} in the case of PEL type.

In this paper we formulate and prove an analogue of the Langlands--Rapoport conjecture for Igusa varieties of Hodge type, building off the work of Kisin \cite{MR3630089} in the case of Shimura varieties. We then use this description of the points to derive a point-counting formula for Igusa varieties of Hodge type, generalizing the fomula in PEL type of Shin \cite{MR2484281}, by drawing on the techniques of Kisin--Shin--Zhu \cite{KSZ}.
\end{abstract}

\tableofcontents
\newpage

\section{Introduction}

\subsection{Context}
\label{sec: history and context}

In this paper we investigate the representations appearing in the cohomology of Igusa varieties, with a view towards applications in the cohomology of Shimura varieties and Langlands program.

A great deal of inspiration comes from the Langlands--Kottwitz method, pioneered by Langlands \cite{MR0354617, MR0437500, langlands_1977, MR546619, langlands_1979} and developed further by Kottwitz in \cite{MR1044820, MR1124982}, which uses geometric and group-theoretic techniques to obtain a trace formula for the cohomology of Shimura varieties that can be compared to the automorphic trace formula, which comparison eventually allows us to relate Galois and automorphic representations. The case treated in \cite{MR1044820,MR1124982} is that of PEL type and hyperspecial level at $p$. In this case Shimura varieties have good reduction and a good moduli structure which makes this approach feasible.

To see ramified representations we must go beyond hyperspecial level at $p$, and the resulting Shimura varieties have bad reduction. An approach in the case of modular curves (i.e., $\GL_2$) was described in Deligne's letter to Piatetski-Shapiro. This approach was extended to some simple Shimura varieties by Harris--Taylor \cite{MR1876802}, where the role of Igusa varieties became clear; and further developed by Mantovan \cite{MR2074715,MR2169874} and Shin \cite{MR2484281, MR2684263, MR2800722, JSTOR41426412}. In short, Mantovan's formula \cite[Thm 22]{MR2169874} allows us to express the cohomology of Shimura varieties in terms of that of Igusa varieties and Rapoport--Zink spaces, with the bad reduction going to the Rapoport--Zink space and the remaining global information going to the Igusa variety. Then a Langlands--Kottwitz style analysis of Igusa varieties \cite{MR2484281, MR2684263} allows us to draw conclusions about Shimura varieties \cite{MR2800722} and Rapoport--Zink spaces \cite{JSTOR41426412}.

Beyond PEL type, the construction of integral models of Shimura varieties no longer guarantees a good moduli structure, so more work is needed to get a good description of points in the special fiber. In particular, the Langlands--Rapoport conjecture describes the points on the special fiber of more general Shimura varieties in a way that is suitable for counting points. This conjecture has essentially been proven by Kisin \cite{MR3630089} for Shimura varieties of abelian type, and the subsequent point-counting work carried out by Kisin--Shin--Zhu \cite{KSZ}.

Igusa varieties and Mantovan's formula have been generalized to Hodge type by Hamacher and Hamacher-Kim \cite{MR3994311, hamacherkim}. The goal of the present work is to derive a trace formula for the cohomology of Igusa varieties of Hodge type, analogous to those for Shimura varieties given in \cite{MR1044820,MR1124982,KSZ} and generalizing the formula for Igusa varieties of PEL type \cite{MR2484281}. This provides an important missing tool for developing our understanding of Shimura varieties of Hodge type, and fits in well with many other recent works in the same thread (cf. \S\ref{sec: applications}).

\subsection{Methods}

This work is descended from the work of Shin \cite{MR2484281} counting points on Igusa varieties of PEL type together, as well as the work of Kisin \cite{MR3630089} and Kisin--Shin--Zhu \cite{KSZ} counting points on Shimura varieties of Hodge (and abelian) type. Our methods draw heavily from these sources, and readers who are not intimate with these sources may prefer to read the author's thesis \cite{thesis}, which provides more detail on how their methods are adapted to the present situation.

As in the case of Shimura varieties, Igusa varieties of Hodge type do not have a good moduli structure suitable for counting points. Thus we need a better description of the points on our Igusa variety.

Our first main theorem is an analogue of the Langlands--Rapoport conjecture for Igusa varieties of Hodge type. This is the subject of \S\ref{igusa LR section}.

Recall (e.g., \cite[Conj. 3.3.7]{MR3630089}) that the Langlands--Rapoport conjecture for Shimura varieties predicts a Frobenius-Hecke equivariant bijection
\[
\sS_{K_p}(G,X)(\oFFp)\overset{\sim}{\longrightarrow}\coprod_{[\phi]}I_\phi(\QQ)\backslash X^p(\phi)\times X_p(\phi),
\]
where $\sS_{K_p}(G,X)$ is the special fiber of the Shimura variety associated to a datum $(G,X)$ at infinite level away from $p$; on the right hand side, the objects are defined in terms of Galois gerbs, but intuitively the disjoint union over $[\phi]$ represents the different isogeny classes on the Shimura variety, the sets $X^p(\phi)$ and $X_p(\phi)$ represent away-from-$p$ isogenies and $p$-power isogenies respectively, and the group $I_\phi(\QQ)$ represents self-isogenies.

The main difference between Shimura varieties and Igusa varieties is the structure at $p$, where Igusa varieties fix an isomorphism class of $p$-divisible groups and add the data of a trivialization. Thus to formulate an analogue of the Langlands--Rapoport conjecture for Igusa varieties, we expect this difference to reflect in the set $X_p(\phi)$.

In the case of Shimura varieties we have
\[
X_p(\phi)\cong X_v(b)=\{g\in G(L)/G(\OO_L):gb\sigma(g)^{-1}\in G(\OO_L)v(p)G(\OO_L)\}
\]
(cf. \ref{ADLV}), where $v$ is a cocharacter of $G$ arising from the Shimura datum, and $b\in G(L)$ is an element essentially recording the Frobenius on the isocrystal associated to a chosen point on the Shimura variety (here $L=\breve{\QQ}_p$ is the completion of the maximal unramified extension of $\QQ_p$). Our chosen point determines a Dieudonn\'e module inside this isocrystal; intuitively, choosing an element $g\in G(L)/G(\OO_L)$ corresponds to transforming this Dieudonn\'e module into another lattice inside the isocrystal, and the condition $gb\sigma(g)^{-1}\in G(\OO_L)v(p)G(\OO_L)$ ensures that this lattice is again a Dieudonn\'e module.

For Igusa varieties we replace this $X_v(b)$ by
\[
J_b(\QQ_p)=\{g\in G(L):gb\sigma(g)^{-1}=b\}
\]
(cf. \ref{Jb}). Intuitively, replacing the condition $gb\sigma(g)^{-1}\in G(\OO_L)v(p)G(\OO_L)$ by the condition $gb\sigma(g)^{-1}=b$ corresponds to fixing an isomorphism class of $p$-divisible group or Dieudonn\'e module (rather than fixing an isogeny class, i.e. isocrystal); and replacing $G(L)/G(\OO_L)$ by $G(L)$ corresponds to adding a trivialization of the $p$-divisible group (i.e. choosing a basis rather than simply a lattice). We have no need to modify the term $X^p(\phi)\cong G(\AA_f^p)$ away from $p$.

The other change in our Langlands--Rapoport conjecture for Igusa varieties is to define a notion of \emph{$\bfb$-admissible morphism} (Definition \ref{b-admissible morphism}) to replace the admissible morphisms appearing in the Langlands--Rapoport conjecture. Since our Igusa variety fixes an isomorphism class of $p$-divisible groups, in particular it fixes an isogeny class, and therefore it lies over a single Newton stratum of the Shimura variety. Restricting to $\bfb$-admissible morphisms corresponds to restricting to isogeny classes in the $\bfb$-stratum.

Indeed, one of the main ideas of the proof (undertaken in \S\ref{isogeny classes on igusa varieties}) is to relate isogeny classes on the Igusa variety and the Shimura variety. Namely, we show that taking the preimage of a Shimura isogeny class along the natural map $\Ig\to\Sh$ gives a bijection between the set of Igusa isogeny classes and the set of Shimura isogeny classes contained in the $\bfb$-stratum (Proposition \ref{isogeny class comparison}).

This relation allows us to use the methods of \cite{MR3630089} (namely Kottwitz triples and their refinements) to establish a bijection between the set of Igusa isogeny classes and the set of conjugacy classes of $\bfb$-admissible morphisms of Galois gerbs, as well as bijections between each individual isogeny class and its parametrizing set. Note that in \cite{MR3630089} the Langlands--Rapoport conjecture is proven only up to possibly twisting the action of $I_\phi(\QQ)$ on $X^p(\phi)\times X_p(\phi)$ by an element $\tau\in I_\phi^\ad(\AA_f)$. A crucial part of \cite{KSZ} is to show that this twist can be taken to satisfy hypotheses that allows us nonetheless to arrive at the expected point-counting formula (namely taking $\tau$ to be tori-rational and to lie in a distinguished space $\Gamma(\cH)_0$, cf. \ref{tau-twists section} or \cite[2.6]{KSZ}). As we use their methods, the same ambiguity of the twist will appear in our case, but in our case too it will not interfere with point-counting.

Thus we arrive at our first main theorem, an analogue of the Langlands--Rapoport conjecture for Igusa varieties of Hodge type.

\begin{theorem*}[\ref{igusa LR}]
There exists a tori-rational element $\tau\in \Gamma(\cH)_0$ which admits a $G(\AA_f^p)\times J_b(\QQ_p)$-equivariant bijection
\[
\Ig_\Sigma(\oFFp)
\overset{\sim}{\longrightarrow}
\coprod_{[\phi]} I_\phi(\QQ)\backslash G(\AA_f^p)\times J_b(\QQ_p),
\]
where the disjoint union ranges over conjugacy classes of $\bfb$-admissible morphisms $\phi:\fQ\to\fG_G$, and the action of $I_\phi(\QQ)$ on $G(\AA_f^p)\times J_b(\QQ_p)$ is twisted by $\tau$.
\end{theorem*}

Next we put this theorem to work as our basis for point-counting, to derive the trace formula for cohomology of Igusa varieties of Hodge type. This is the subject of \S\ref{sec: point-counting}.

This falls into two steps. First, we interpret the appropriate class of test functions as correspondences on our Igusa varieties, and use Fujiwara's trace formula to convert the problem of computing traces of the action on cohomology to the problem of computing fixed points of these correspondences. We can use our first main theorem above to describe the fixed points of these correspondences, resulting in a preliminary form of our point-counting formula parametrized by Galois-gerb-theoretic LR pairs $(\phi,\varepsilon)$ (cf. \ref{def: LR pair}). This is undertaken in \S\S\ref{sec: acceptable functions and fujiwara}-\ref{preliminary point-counting section}.

For comparison with the automorphic trace formula, the second step is to re-parametrize our point-counting formula in the more group-theoretic terms of Kottwitz parameters (cf. \ref{kottwitz parameter}). For this we adapt the techniques of \cite[\S3]{KSZ}. The theory required is quite analogous, but the relevant class of LR pairs is different; instead of their $p^n$-admissible pairs, we define notions of \emph{$\bfb$-admissible} and \emph{acceptable} pairs (Definitions \ref{def: LR pair} and \ref{def: acceptable LR pair}). Then we need to re-work a substantial part of the theory under these new hypotheses. Fortunately it is possible to prove essentially the same results, though the arguments often require different techniques. This is undertaken in \S\S\ref{LR/KP section}-\ref{sec: final point counting}.

In the end we arrive at our second main theorem, the point-counting formula for Igusa varieties of Hodge type.

\begin{theorem*}[\ref{point counting formula}]
For any acceptable function $f\in C^\infty_c(G(\AA_f^p)\times J_b(\QQ_p))$, we have
\begin{align*}
\tr(f\mid H_c(\Ig_\Sigma,&\sL_\xi))
=
\sum_{\gamma_0\in\Sigma_{\Rell}(G)\phantom{W}}\sum_{(a,[b_0])\in\KP(\gamma_0)}
\\
&\frac{\lvert\Sha_G(\QQ,G_{\gamma_0}^\circ)\rvert}{\lvert(G_{\gamma_0}/G_{\gamma_0}^\circ)(\QQ)\rvert}
\vol\left(I^\circ_{\fc}(\QQ)\backslash I^\circ_{\fc}(\AA_f)\right) O_{\gamma\times\delta}^{G(\AA_f^p)\times J_b(\QQ_p)}(f)\tr(\xi(\gamma_0))
\end{align*}
where $I_\fc$ is the inner form of $G_{\gamma_0}^\circ$ associated to the Kottwitz parameter $\fc=(\gamma_0,a,[b_0])$ as in \ref{par: I_c}, and $\gamma,\delta$ are the elements belonging to the classical Kottwitz parameter $(\gamma_0,\gamma,\delta)$ associated to $\fc$ as in \ref{par: KP -> classical KP}.
\end{theorem*}

\subsection{Applications}
\label{sec: applications}

For applications, it is necessary to stabilize our point-counting formula. This is work in progress by Bertoloni Meli and Shin \cite{BMS}.

As described in \S\ref{sec: history and context}, we expect our formula to be useful in combination with Mantovan's formula (due to Hamacher-Kim \cite{hamacherkim} in Hodge type) to investigate the cohomology of Shimura varieties and Rapoport--Zink spaces, as has been done to great effect in \cite{MR1876802,MR2800722,JSTOR41426412} for the case of PEL type. A particular example is the recent work of Bertoloni Meli and Nguyen \cite{meli2021kottwitz} who prove the Kottwitz conjecture for a certain class of unitary similitude groups; their method uses Mantovan's formula and the point counting formula for Igusa varieties of PEL type, and therefore our formula in Hodge type is an important ingredient in extending their results.

Another promising application is to generalize the results of Caraiani-Scholze on cohomology of Shimura varieties \cite{MR3702677, CS19non-cpt}. A crucial part of their approach is to use Hodge-Tate period map $\pi_{\rm HT}:\Sh\to\Fl$ to transfer problems from a Shimura variety to an associated flag variety. They realize the fibers of $\pi_{\rm HT}$ as essentially Igusa varieties, and use the point-counting formula for Igusa varieties of PEL type \cite{MR2484281, MR2684263} to describe the cohomology of those fibers. Thus our second main theorem above is one of the crucial ingredients needed to generalize their arguments to Hodge type.

Even more recently, Kret and Shin \cite{KS2021h0} have given a description of the $H^0$ cohomology of Igusa varieties in terms of automorphic representations by combining our point-counting formula with automorphic trace formula techniques.

\subsection{Acknowledgments}

The author owes many thanks to Sug Woo Shin for his patient and thoughtful advising, without which this work would not have been possible. Thanks also to Xinyi Yuan, Martin Olsson, Ken Ribet, and all of the faculty in and around number theory for creating a stimulating environment in which to study. Many thanks to Alex Youcis, Alexander Bertoloni Meli, Rahul Dalal, Ian Gleason, Zixin Jiang, Dong Gyu Lim, and Koji Shimizu for helpful conversations and support. This work was generously funded by NSF award numbers 1752814 
and 1646385. 

\section{Background}

\subsection{Isocrystals}

\paragraph{}
\label{isocrystal definition}
Let $L=\breve{\QQ}_p$ the completion of the maximal unramified extension of $\QQ_p$ and $\sigma$ the lift of Frobenius on $L$ (coming from $\breve{\ZZ}_p=W(\oFFp)$). An \emph{isocrystal} over $\oFFp$ is a finite-dimensional vector space $V$ over $L$ equipped with a $\sigma$-semilinear bijection $F:V\to V$, which we call its Frobenius map. A morphism of isocrystals $f:(V_1,F_1)\to(V_2,F_2)$ is a linear map $f:V_1\to V_2$ intertwining their Frobenius maps, i.e. $f\circ F_1=F_2\circ f$.

For $G$ a reductive group, an \emph{isocrystal with $G$-structure} is an exact faithful tensor functor
\[
\Rep_{\QQ_p}(G)\to\Isoc
\]
where $\Isoc$ is the category of isocrystals \cite[Def 3.3]{MR1411570}.

\paragraph{}
\label{B(G)}
Giving an isocrystal structure $F$ on $L^n$ is the same as choosing an element $b\in \GL_n(L)$, via $F=b\sigma$. The isocrystal structures defined by $b_0,b_1\in \GL_n(L)$ are isomorphic precisely when $b_0,b_1$ are $\sigma$-conjugate, i.e. $b_1=gb_0\sigma(g)^{-1}$ for some $g\in \GL_n(L)$. This gives an identification between isomorphism classes of isocrystals and $\sigma$-conjugacy classes in $\GL_n(L)$.

More generally, we can associate an isocrystal with $G$-structure to an element $b\in G(L)$ by setting ``$F=b\sigma$''; namely, the isocrystal with $G$-structure defined by the functor
\begin{align*}
\Rep_{\QQ_p}(G) & \to \Isoc \\
(V,\rho) & \mapsto (V\otimes_{\QQ_p}L, \rho(b)(\id_V\otimes\sigma)).
\end{align*}
This association identifies the set of isomorphism classes of isocrystals with $G$-structure with the set of $\sigma$-conjugacy classes in $G(L)$. We denote this common set by $B(G)$, and we write $[b]\in B(G)$ for the $\sigma$-conjugacy class of an element $b\in G(L)$. If $G$ is connected, then in fact every element of $B(G)$ has a representative in $G(\QQ_{p^r})$ for some finite unramified extension $\QQ_{p^r}$ of $\QQ_p$ \cite[4.3]{MR809866}. Given a cocharacter $\mu$ of $G$, there is a distinguished finite subset $B(G,\mu)$ of \emph{$\mu$-admissible} classes, defined in \cite[\S6]{MR1485921}.

\paragraph{}
\label{slopes}
The Dieudonn\'e-Manin classification gives a concrete description of the category of isocrystals. Given $\lambda=\frac{r}{s}$ a rational number in lowest terms ($s>0$), we can define an isocrystal $E_\lambda=L\langle F\rangle/(F^s-p^r)$ with Frobenius given by left multiplication by $F$; here $L\langle F\rangle$ is the twisted polynomial ring where $Fx=\sigma(x)F$ for $x\in L$. The Dieudonn\'e-Manin classification states that the category of isocrystals is semi-simple, and the simple objects are $E_\lambda$ for $\lambda\in\QQ$. In other words, an isocrystal is determined by a finite set of rational numbers $\lambda$ with multiplicities. These rational numbers are called its \emph{slopes}, and the decomposition of an isocrystal into a direct sum $V=\bigoplus_\lambda V_\lambda$ of subspaces of slope $\lambda$ (so $V_\lambda\cong E_\lambda^{\oplus r}$ for some $r$) is called the \emph{slope decomposition}.

We can also associate slopes to an isocrystal with $G$-structure, in a slightly different form. Let $\DD$ be the pro-torus with character group $X^*(\DD)=\QQ$. Then an isocrystal with slope decomposition $V=\bigoplus_\lambda V_\lambda$ produces a fractional cocharacter $\DD\to\GL(V)$ over $L$ defined by $\DD$ acting on $V_\lambda$ by the character $\lambda\in \QQ=X^*(\DD)$.

Now let $b\in G(L)$, defining an isocrystal with $G$-structure. Given a representation $(V,\rho)\in \Rep_{\QQ_p}(G)$, this isocrystal with $G$-structure produces an isocrystal on $V_L=V\otimes_{\QQ_p}L$, and therefore a fractional cocharacter $\nu_\rho:\DD\to\GL(V_L)$. The \emph{slope homomorphism} of $b$ is the unique fractional cocharacter $\nu_b:\DD\to G$ over $L$ satisfying $\nu_\rho=\rho\circ\nu_b$ for all $p$-adic representations $\rho$ of $G$.

Alternatively, $\nu_b$ can be defined (cf. \cite[4.3]{MR809866}) as the unique element of $\Hom_L(\DD,G)$ for which there exists an $n>0$ and $c\in G(L)$ such that
\begin{itemize}
\item $n\nu_b\in \Hom_L(\GG_m,G)$,
\item $\Int(c)\circ n\nu_b$ is defined over a finite unramified extension $\QQ_{p^n}$ of $\QQ_p$, and
\item $c(b\sigma)^nc^{-1}=c\cdot n\nu_b(\pi)\cdot c^{-1}\cdot\sigma^n$ (considered in $G(L)\rtimes\langle \sigma\rangle$).
\end{itemize}

From the definition we can see the slope homomorphism transforms nicely under the action of $\sigma$ and conjugation by $G(L)$:
\begin{itemize}
\item $\nu_{\sigma(b)}=\sigma(\nu_b)$;
\item $\nu_{gb\sigma(g)^{-1}}=\Int(g)\circ\nu_b$.
\end{itemize}

The following lemma states that, to check if an element $g\in G(L)$ commutes with $\nu_b$, it suffices to check on a single faithful representation.

\begin{lemma}
\label{slope decomposition}
Let $b\in G(L)$, defining an isocrystal with $G$-structure. Let $\rho:G\to\GL(V)$ be a faithful $p$-adic representation, and $(V_L,\rho(b)\sigma)$ the associated isocrystal. If $g\in G(L)$ (acting via $\rho(g)$) preserves the slope decomposition of this isocrystal, then $g$ commutes with the slope homomorphism $\nu_b$.
\end{lemma}
\begin{proof}
Suppose $g\in G(L)$ preserves the slope decomposition $V_L=\bigoplus_\lambda V_\lambda$. Consider the two fractional cocharacters $\nu_b$ and $\Int(g)\circ\nu_b$ of $G$; we want to show they are equal.

Composing with $\rho$, we get two fractional cocharacters of $\GL(V_L)$,
\[
\nu_\rho=\rho\circ\nu_b\qquad\text{ and }\qquad \Int(\rho(g))\circ\nu_\rho=\rho\circ\Int(g)\circ\nu_b.
\]
The fractional cocharacter $\nu_\rho:\DD\to\GL(V_L)$ is defined by $\DD$ acting on $V_\lambda$ by the character $\lambda\in\QQ=X^*(\DD)$. By assumption $\rho(g)$ preserves the slope decomposition, so we see that $\Int(\rho(g))\circ\nu_\rho$ also acts on $V_\lambda$ by the character $\lambda$. Since these two actions are the same, we find $\nu_\rho=\Int(\rho(g))\circ\nu_\rho$; and by the monomorphism property of $\rho$ this implies $\nu_b=\Int(g)\circ\nu_b$, as desired.
\end{proof}

\subsection{Acceptable Elements of $J_b(\QQ_p)$}
\label{acceptable elements section}

\paragraph{}
\label{Jb}
For $b\in G(L)$, define an algebraic group $J_b$ (or $J_b^G$ when it is helpful to specify the group) over $\QQ_p$ by defining its points for a $\QQ_p$-algebra $R$ by
\[
J_b(R)=\{g\in G(R\otimes_{\QQ_p}L):gb\sigma(g)^{-1}=b\},
\]
and define
\[
M_b=\text{centralizer in $G$ of }\nu_b.
\]
If necessary we can change $b$ inside its $\sigma$-conjugacy class to ensure that $M_b$ is defined over $\QQ_p$ (this is always possible if $G$ is quasi-split, cf. \cite[Prop 6.2]{MR809866}). Then $J_b$ is the automorphism group of the isocrystal with $G$-structure defined by $b$ (indeed, the condition $gb\sigma(g)^{-1}=b$ precisely means that $g$ commutes with $b\sigma$), and furthermore $J_b$ is an inner form of $M_b$. Changing $b$ by $\sigma$-conjugation in $G(L)$ does not essentially change the situation: if $b_0=gb_1\sigma(g)^{-1}$, then we have a canonical isomorphism
\begin{align*}
J_{b_1} & \overset{\sim}{\longrightarrow}J_{b_0} \\
x & \longmapsto gxg^{-1}
\end{align*}
and $M_{b_0}=\Int(g) M_{b_1}$.

\paragraph{}
\label{fixed isocrystal}

Choose a faithful representation $V$ of $G$. Then our isocrystal with $G$-structure associated to $b$ produces an isocrystal $(V_L,b\sigma)$ (we abuse notation by suppressing the map $\rho$, e.g. writing $b\sigma$ rather than $\rho(b)\sigma$ for the Frobenius). The group $J_b(\QQ_p)$ acts on this isocrystal by linear automorphisms, via its natural inclusion in $G(L)$---the defining condition $gb\sigma(g)^{-1}=b$ exactly means that $g$ commutes with $b\sigma$. Write $V_L=\bigoplus_i V_{\lambda_i}$ for the slope decomposition of our isocrystal, with slopes in decreasing order $\lambda_1>\lambda_2>\cdots>\lambda_r$.

\begin{definition}
\label{acceptable element}
Define an element $\delta\in J_b(\QQ_p)$ to be \emph{acceptable} (or say $\delta$ is \emph{acceptable with respect to $b$}) if, regarding $\delta=(\delta_i)\in\prod \GL(V_{\lambda_i})$, any eigenvalues $e_i$ of $\delta_i$ and $e_j$ of $\delta_j$ with $i<j$ (i.e. $\lambda_i>\lambda_j$) satisfy $v_p(e_i)<v_p(e_j)$.
\end{definition}

\paragraph{}
\label{acceptable elements are plentiful}
An important example of an acceptable element is defined as follows. Choose an integer $s$ so that $s\lambda_i$ is an integer for all slopes $\lambda_i$ of our isocrystal. Define an element of $J_b(\QQ_p)$, formally written as $\fr^{s}$, by acting on $V_{\lambda_i}$ by $p^{s\lambda_i}$. For the inverse of this element we write $\fr^{-s}$, defined by acting on $V_{\lambda_i}$ by $p^{-s\lambda_i}$. Any element of $J_b(\QQ_p)$ becomes acceptable after multiplying by a sufficiently high power of $\fr^{-s}$---note that $\fr^{-s}$ decreases $v_p(e_i)$ by $s\lambda_i$, and as the slopes are strictly decreasing a high enough power will satisfy Definition \ref{acceptable element}.

\begin{lemma}
\label{not root argument}
Let $\varepsilon\in G(\QQ_p^\ur)$ semi-simple. Suppose $\varepsilon$ commutes with $b\sigma$, and is therefore an element of $J_b(\QQ_p)$, and furthermore is acceptable. Then $G_\varepsilon\subset M_b$.
\end{lemma}
\begin{proof}
We continue to work with our fixed isocrystal of \ref{fixed isocrystal}, with slope decomposition
\[
V_L=\bigoplus_i V_{\lambda_i}.
\]
Since $\varepsilon$ is semi-simple, its action on $V_L$ is diagonalizable, and since it commutes with $b\sigma$ it preserves the slope components $V_{\lambda_i}$. Thus each slope component has a basis of eigenvectors for the action of $\varepsilon$. The acceptable condition implies that $\varepsilon$ has different eigenvalues on different slope components, so in fact each slope component is a direct sum of full eigenspaces of $\varepsilon$.

Now, suppose $x\in G_\varepsilon$, so $x$ commutes with $\varepsilon$. Then $x$ must preserve the eigenspaces of $\varepsilon$. Since the slope components are direct sums of eigenspaces of $\varepsilon$, we see that $x$ preserves the slope decomposition. By Lemma \ref{slope decomposition} this implies $x\in M_b$.
\end{proof}

\begin{lemma}
\label{same nu}
Let $b_0,b_1$ be $\sigma$-conjugate elements of $G(L)$. Suppose that there is a semi-simple element $\varepsilon\in G(L)$ such that $\varepsilon$ lies in both $J_{b_0}(\QQ_p)$ and $J_{b_1}(\QQ_p)$ and furthermore is acceptable with respect to both $b_0$ and $b_1$. Then $\nu_{b_0}=\nu_{b_1}$.
\end{lemma}
\begin{proof}
We continue to work with our fixed representation $V$ of \ref{fixed isocrystal}, and write $(V_L,b_0\sigma)$ and $(V_L,b_1\sigma)$ for the isocrystals produced from $V$ by $b_0$ and $b_1$. These two isocrystal structures induce two slope decompositions
\[
V_L=\bigoplus_i V_{\lambda_i,0} \qquad\text{ and }\qquad V_L=\bigoplus_i V_{\lambda_i,1}.
\]
Since $b_0,b_1$ are $\sigma$-conjugate, we see these two isocrystal structures are abstractly isomorphic. Taking $g\in G(L)$ with $b_1=gb_0\sigma(g)^{-1}$, the slope decompositions are related by $V_{\lambda_i,1}=g\cdot V_{\lambda_i,0}$. In particular, for each slope $\lambda_i$, the slope components $V_{\lambda_i,0}$ and $V_{\lambda_i,1}$ have the same dimension.

Now consider the element $\varepsilon$, which we assume to commute with $b_0\sigma,b_1\sigma$ so that it acts on both isocrystal structures, and furthermore is acceptable for both. As in the proof of \ref{not root argument}, the semi-simple and acceptable conditions imply that each slope component $V_{\lambda_\bullet,\bullet}$ is a direct sum of full eigenspaces of $\varepsilon$. Combining the fact that $\dim V_{\lambda_i,0}=\dim V_{\lambda_i,1}$ with the fact that (by the acceptable condition) the eigenspaces appearing in the various slope components must be ordered by the $p$-adic valuation of the corresponding eigenvalue, this shows that $V_{\lambda_i,0}$ and $V_{\lambda_i,1}$ must consist of the same eigenspaces; that is, $V_{\lambda_i,0}=V_{\lambda_i,1}$ for all slopes $\lambda_i$.

Now we've shown $g\cdot V_{\lambda_i,0}=V_{\lambda_i,1}=V_{\lambda_i,0}$; that is, $g$ preserves the slope decomposition with respect to $b_0$, and by Lemma \ref{slope decomposition} this means $g$ centralizes $\nu_{b_0}$. On the other hand, we chose $g$ to satisfy $b_1=gb_0\sigma(g)^{-1}$, so as in \ref{slopes} we have
\[
\nu_{b_1}=\nu_{gb_0\sigma(g)^{-1}}=\Int(g)\circ \nu_{b_0}.
\]
Thus $\nu_{b_0}=\nu_{b_1}$, as desired.
\end{proof}

\paragraph{}
\label{S_b}
We now define a submonoid $S_b\subset J_b(\QQ_p)$ which will play a role in defining the group actions on our Igusa varieties. Let $\delta\in J_b(\QQ_p)$, and suppose that $\delta^{-1}$ is an isogeny. Regard $\delta=(\delta_i)\in\prod\GL(V_{\lambda_i})$ as in \ref{fixed isocrystal}--\ref{acceptable element}, and for each $i$ let $e_i(\delta)$ and $f_i(\delta)$ be the minimal and maximal (respectively) integers such that
\[
\ker p^{f_i(\delta)}\subset\ker\delta_i^{-1}\subset\ker p^{e_i(\delta)}.
\]
Then $S_b$ is the submonoid of $J_b(\QQ_p)$ defined by
\[
S_b=\{\delta\in J_b(\QQ_p):\delta^{-1}\text{ an isogeny, and } f_{i-1}(\delta)\ge e_i(\delta)\text{ for all }i\}.
\]
For example, $S_b$ contains $p^{-1}$, since multiplication by $p$ is an isogeny and $f_{i-1}(p^{-1})=e_i(p^{-1})=1$ for all $i$. Also $S_b$ contains $\fr^{-s}$, as $\fr^s$ is an isogeny and $f_{i-1}(\fr^{-s})=s\lambda_{i-1}>e_i(\fr^{-s})=s\lambda_i$ for all $i$. In fact, $J_b(\QQ_p)$ is generated as a monoid by $S_b$ together with $p$ and $\fr^s$---in other words, any element of $J_b(\QQ_p)$ can be translated into $S_b$ by multiplying by high enough powers of $p^{-1}$ and $\fr^{-s}$.

\subsection{$p$-divisible Groups}

In this section we briefly recall some definitions related to $p$-divisible groups that will be essential for defining and working with Igusa varieties.

\paragraph{}
\label{p-divisible groups basics}
A \emph{$p$-divisible group} over a scheme $S$ is an fppf sheaf of abelian groups $\sG$ on $S$ such that
\begin{itemize}
\item $\sG[p^n]=\ker(\sG\overset{p^n}{\longrightarrow}\sG)$ is a finite locally free group scheme for all $n$,
\item $\sG=\dlim_n\sG[p^n]$, and
\item $\sG\overset{p}{\longrightarrow}\sG$ is an epimorphism.
\end{itemize}
An \emph{isogeny} of $p$-divisible groups $f:\sG_1\to\sG_2$ is an epimorphism whose kernel is a finite locally free group scheme. A \emph{quasi-isogeny} is a global section $f$ of $\Hom_S(\sG_1,\sG_2)\otimes_\ZZ\QQ$ which Zariski-locally admits an integer $n$ so that $p^nf$ is an isogeny.

The examples we will be most concerned with are $p$-divisible groups arising from abelian varieties: if $A$ is an abelian variety, then the limit $A[p^\infty]=\dlim A[p^n]$ over the natural inclusion maps $A[p^n]\hookrightarrow A[p^{n+1}]$ forms a $p$-divisible group.

We write $\sG\mapsto\DD(\sG)$ for the contravariant Dieudonn\'e module functor, which gives a contravariant equivalence of categories between the category of $p$-divisible groups over $\oFFp$ and the category of Dieudonn\'e modules (e.g. \cite{MR0344261}). By composing with the functor from Dieudonn\'e modules to isocrystals, we get a contravariant functor $\sG\to\VV(\sG)$.

While the Dieudonn\'e module records a $p$-divisible group up to isomorphism, the isocrystal records a $p$-divisible group up to isogeny. To be precise, consider the \emph{isogeny category} of $p$-divisible groups over $S$, whose objects are $p$-divisible groups and whose morphisms are global sections of $\Hom_S(\sG_1,\sG_2)\otimes_\ZZ\QQ$; in this category quasi-isogenies are isomorphisms. Then the functor $\sG\to\VV(\sG)$ is an equivalence of categories between the isogeny category of $p$-divisible groups over $\oFFp$ and the category of isocrystals.

A $p$-divisible group can be equipped with $G$ structure; in our case this will usually consist of a set of tensors on $\DD(\sG)$. Then the associations above produce a Dieudonn\'e module with $G$-structure and an isocrystal with $G$-structure. In particular a $p$-divisible group with $G$-structure has an associated isocrystal with $G$-structure, whose isomorphism class is recorded by a class $[b]\in B(G)$, and we say the $p$-divisible group is of \emph{type $[b]$}.

\paragraph{}
\label{completely slope divisible}
Say that a $p$-divisible group is \emph{isoclinic} if it has only a single slope (possibly with multiplicity). A \emph{slope filtration} for a $p$-divisible group $\sG$ with slopes $\lambda_1>\cdots>\lambda_r$ is a filtration
\[
0=\sG_0\subset\cdots\subset\sG_r=\sG
\]
such that each successive quotient $\sG_i/\sG_{i-1}$ is isoclinic of slope $\lambda_i$. If it exists, it is unique. A slope filtration always exists for a $p$-divisible group over a field of positive characteristic, and the filtration splits canonically if the field is perfect \cite{MR0417192}.

We say $\sG$ is \emph{completely slope divisible} if it has a slope filtration such that for each successive quotient $X$ of slope $\lambda=\frac{a}{b}$, the quasi-isogeny $p^{-a}\Frob^b:X\to X^{(p^b)}$ is an isogeny. Over $\oFFp$, this is equivalent to being a direct sum of isoclinic $p$-divisible groups defined over finite fields \cite{MR1938119}. For a class $[b]\in B(G)$, there is a completely slope divisible $p$-divisible group with $G$-structure of type $[b]$ exactly when $[b]$ has a representative in $G(\QQ_{p^r})$ for some $r$; as in \ref{B(G)}, this is the case for all classes in $B(G)$ if $G_{\QQ_p}$ is connected.

\subsection{Igusa Varieties of Siegel Type}
\label{siegel type}

In this section we review the case of Siegel type, as it is required for understanding Hodge type. Many objects will be decorated with a tick $\bullet'$ to distinguish them from the analogous objects of Hodge type introduced in \S\ref{hodge type}.

\paragraph{}
Let $V$ be a $\ZZ$-module of finite rank $2r$ and $\psi$ a symplectic pairing on $V$. For any ring $R$, we write $V_R=V\otimes_\ZZ R$. This data gives rise to a Shimura datum of Siegel type $(\GSp,S^\pm)$ consisting of the group of symplectic similitudes
\[
\GSp(R)=\{g\in \GL(V_R): \psi(gx, gy)=c(g)\psi(x,y) \text{ for some }c(g)\in R^\times\}
\]
and the Siegel double space $S^\pm$, defined to be the set of complex structures $J$ on $V_\RR$ such that the form $(x,y)\mapsto\psi(x,Jy)$ is symmetric and positive or negative definite. We can regard such a complex structure on $V_\RR$ also as a Hodge structure on $V$, or also (in the usual manner of a Shimura datum) as a homomorphism $\bbS\to \GSp_\RR$, where $\bbS$ denotes the Deligne torus $\bbS=\Res_{\CC/\RR}\GG_m$.

We can consider $\GSp$ as a reductive group over $\ZZ_{(p)}$, acting on $V_{\ZZ_{(p)}}$, and let $K_p'=\GSp(\ZZ_p)\subset \GSp(\QQ_p)$ be the corresponding hyperspecial subgroup. We write $K'=K'_pK'^p$, where $K'^p\subset \GSp(\AA_f^p)$ is a sufficiently small compact open subgroup.

The Shimura variety $\Sh_{K'}(\GSp,S^\pm)$ associated to this Shimura datum is defined over $\QQ$, and has a canonical integral model which can be defined by a moduli problem. Consider the category of abelian schemes up to prime-to-p isogeny, whose objects are abelian schemes, and whose morphisms are $\Hom(A,B)\otimes_\ZZ\ZZ_{(p)}$ (where $\Hom(A,B)$ denotes the usual homomorphisms of abelian schemes). Isomorphisms in this category are called prime-to-p quasi-isogenies. Our moduli problem has the form
\begin{align*}
\textbf{Schemes}_{/\ZZ_{(p)}} & \longrightarrow \textbf{Sets} \\
X & \longmapsto \{(A,\lambda,\eta^p_{K'})\}/\sim
\end{align*}
where
\begin{itemize}
\item $A$ is an abelian scheme over $X$ up to prime-to-p isogeny;
\item $\lambda$ is a weak polarization of $A$, i.e. a prime-to-p quasi-isogeny $\lambda:A\to A^\vee$ modulo scaling by $\ZZ_{(p)}^\times$, some multiple of which is a polarization;
\item $\eta^p_{K'}\in\Gamma(X,\underline{\Isom}(V_{\AA_f^p},\hat{V}^p(A))/K'^p)$ is a $K'^p$-level structure, where we regard $\hat{T}^p(A)=\ilim_{p\not\mid n}A[n]$ and $\hat{V}^p(A)=\hat{T}^p(A)\otimes_\ZZ\QQ$ as \'etale sheaves on $X$, and define $\underline{\Isom}(V_{\AA_f^p},\hat{V}^p(A))$ to be the \'etale sheaf of isomorphisms compatible with the pairings induced by $\psi$ and $\lambda$ up to $\AA_f^p{}^\times$-scalar; and
\item two triples are equivalent $(A_1,\lambda_1,\eta^p_{K',1})\sim(A_2,\lambda_2,\eta^p_{K',2})$ if there is a prime-to-p quasi-isogeny $A_1\to A_2$ sending $\lambda_1\to\lambda_2$ and $\eta^p_{K',1}$ to $\eta^p_{K',2}$. 
\end{itemize}
If $K'^p$ is sufficiently small, this moduli problem is represented by a smooth $\ZZ_{(p)}$-scheme $\sS_{K'}(\GSp,S^\pm)$, which is the canonical integral model of $\Sh_{K'}(\GSp,S^\pm)$. By virtue of the moduli structure, it carries a universal polarized abelian scheme $\cA'\to\sS_{K'}(\GSp,S^\pm)$. Write
\[
\osS_{K'}(\GSp,S^\pm)=\sS_{K'}(\GSp,S^\pm)\otimes_{\ZZ_{(p)}}\FF_p
\]
for the special fiber.

\paragraph{}
\label{siegel newton stratification}
The universal polarized abelian scheme gives rise to a universal polarized $p$-divisible group $(\cA'[p^\infty],\lambda)$ and hence isocrystal with $\GSp$-structure over $\osS_{K'}(\GSp,S^\pm)$. Restricting to a geometric point $\overline{x}\to\osS_{K'}(\GSp,S^\pm)$ gives an isocrystal with $\GSp$-structure over $\overline{x}$, which is classified by an element $\bfb_x\in B(\GSp)$ depending only on the topological point $x$ underlying $\overline{x}$. For each $\bfb\in B(\GSp)$, we define the \emph{Newton stratum}
\[
\osS^{(\bfb)}_{K'}(\GSp,S^\pm)=\{x\in \osS_{K'}(\GSp,S^\pm):\bfb_x=\bfb\}
\]
to be the locus in $\osS_{K'}(\GSp,S^\pm)$ where the universal $p$-divisible group is of type $\bfb$. It is a locally closed subset, which we promote to a subscheme by taking the reduced subscheme structure.

Fixing a class $\bfb\in B(\GSp)$ whose Newton stratum is non-empty, let $(\Sigma,\lambda_\Sigma)$ be a $p$-divisible group with $\GSp$-structure of type $\bfb$, i.e.
\begin{itemize}
\item $\Sigma$ a $p$-divisible group over $\oFFp$ and
\item $\lambda_\Sigma$ a polarization of $\Sigma$, such that
\item there is an isomorphism $\DD(\Sigma)\overset{\sim}{\to}V_{\OO_L}$ preserving the pairings induced by $\lambda_\Sigma$ and $\psi$, and taking the Frobenius on $\DD(\Sigma)$ to an endomorphism $b\sigma$ on $V_{\OO_L}$ with $b\in \GSp(\OO_L)$ belonging to the class $\bfb\in B(\GSp)$.
\end{itemize}
Note choosing a different isomorphism changes $b$ by $\sigma$-conjugacy in $\GSp(\OO_L)$, so its class in $B(\GSp)$ is well-defined.

Define the \emph{central leaf} corresponding to $(\Sigma,\lambda_\Sigma)$ by
\[
C'_{\Sigma,K'}=\{x\in\osS^{(\bfb)}_{K'}(\GSp,S^\pm):(\cA[p^\infty]_x,\lambda_x)\otimes_{k(x)}\overline{k(x)}\cong(\Sigma,\lambda_\Sigma)\otimes_{\oFFp}\overline{k(x)}\}.
\]
This is a closed subset of the Newton stratum, and when equipped with the reduced subscheme structure, is smooth \cite[Prop 1]{MR2169874}.

\paragraph{}
\label{siegel igusa varieties}
Now we assume that $\Sigma$ is completely slope divisible (which we can do since $\GSp$ is connected, cf. \ref{completely slope divisible}). Then the universal $p$-divisible group $\cA'[p^\infty]$ over $C_\Sigma$, being isomorphic to $\Sigma$ (over each geometric generic point of $C_\Sigma$), is also completely slope divisible. Let $\cA'[p^\infty]^{(i)}$ be the successive quotients of the slope filtration, and define $\cA'[p^\infty]^{\rm sp}=\bigoplus_i\cA'[p^\infty]^{(i)}$ to be the associated split $p$-divisible group, which inherits a polarization from $\cA'[p^\infty]$. Denote its $p^m$-torsion by $\cA'[p^m]^{\rm sp}$.

The \emph{level-$m$ Igusa variety} of Siegel type $\Ig'_{\Sigma,K',m}$ is a smooth $\oFFp$-scheme, finite \'etale and Galois over $C'_{\Sigma,K'}$, defined by the moduli problem
\begin{align*}
\Ig'_{\Sigma,K',m}(X)
=
\Big\{(A,\lambda,\eta^p_{K'},j_m):\;&(A,\lambda,\eta^p_{K'})\in C'_{\Sigma,K'}(X),
\\
&j_m:\Sigma[p^m]\times_{\oFFp}X\overset{\sim}{\to}\cA'[p^m]^{\rm sp}\times_{C_\Sigma}X\Big\},
\end{align*}
where $j_m$ is an isomorphism preserving polarizations up to $(\ZZ/p^m)^\times$-scalar, and extending \'etale locally to any higher level $m'\ge m$. Essentially we are adding level structure over $C'_{\Sigma,K'}$ in the form of a trivialization $j_m$ of $\cA'[p^m]$.

Let $\Ig'_{\Sigma,K'}=\ilim_m\Ig'_{\Sigma,K',m}$ be the Igusa variety at infinite $m$-level, and $\sJ'_{\Sigma,K'}=\Ig'^{(p^{-\infty})}_{\Sigma,K'}$ its perfection. Then $\sJ'_{\Sigma,K'}$ is the moduli space over $C'_{\Sigma,K'}$ parametrizing isomorphisms
\[
j:\Sigma\times_{\oFFp}C'_{\Sigma,K'}\to\cA'[p^\infty]
\]
preserving polarizations up to scaling. Note that the slope filtration splits canonically over a perfect base \cite[4.1]{newtonstratification}, so we no longer need to impose the splitting on $\cA'[p^\infty]$.

The group $\GSp(\AA_f^p)$ acts on the system $\Ig'_{\Sigma,K'}$, varying $K'$, by acting on the level structure $\eta_{K'}^p$. This action is inherited from the Siegel modular variety, since it happens away from $p$ (i.e. it does not interact with the Igusa level structure). The system $\Ig'_{\Sigma,K'}$ also receives an action of the submonoid $S_b\subset J_b^{\GSp}(\QQ_p)$ of \ref{S_b}, which extends to an action of the full group $J_b^{\GSp}(\QQ_p)$ on the perfections $\sJ'_{\Sigma,K'}$ and on \'etale cohomology.

\subsection{Igusa Varieties of Hodge Type}
\label{hodge type}

\paragraph{}
\label{par: hodge embedding}
Let $(G,X)$ be a Shimura datum of Hodge type: $G$ is a connected reductive $\QQ$-group, which we further assume to be unramified at $p$; $X$ is a $G(\RR)$-conjugacy class of homomorphisms $h:\bbS\to G_{\RR}$; and there exists a closed embedding $G\hookrightarrow\GSp$ which sends $X$ to $S^\pm$. Denote by $E$ the reflex field of $(G,X)$. 
We permanently fix a Hodge embedding $G\hookrightarrow\GSp$.

As $G$ is unramified, it has a reductive model $G_{\ZZ_{(p)}}$ over $\ZZ_{(p)}$ and corresponding hyperspecial subgroup $K_p=G_{\ZZ_{(p)}}(\ZZ_p)\subset G(\QQ_p)$. As in \cite[1.3.3]{MR3630089}, there is a $\ZZ_{(p)}$-lattice $V_{\ZZ_{(p)}}\subset V_{\QQ}$ such that the embedding $G\hookrightarrow\GSp$ is induced by an embedding $G_{\ZZ_{(p)}}\hookrightarrow\GL(V_{\ZZ_{(p)}})$. Enlarging our symplectic space $V$ if necessary, $\psi$ induces a perfect pairing on $V_{\ZZ{(p)}}$, so we can define a hyperspecial subgroup $K_p'=\GSp(V_{\ZZ_{(p)}})(\ZZ_p)\subset \GSp(\QQ_p)$ which is compatible in the sense that the embedding $G\hookrightarrow\GSp$ takes $K_p$ into $K_p'$. For any compact open $K^p\subset G(\AA_f^p)$ there is a $K'^p\subset\GSp(\AA_f^p)$ so that $K=K_pK^p\subset K'=K'_pK'^p$ and the natural map
\[
\Sh_K(G,X)\to\Sh_{K'}(\GSp,S^\pm)
\]
is a closed embedding.

In the case of Hodge type, our integral models will not be defined by a moduli problem. Instead, consider the composition
\[
\Sh_K(G,X)\to\Sh_{K'}(\GSp,S^\pm)\to\sS_{K'}(\GSp,S^\pm),
\]
and let $\sS_K(G,X)$ be the closure of $\Sh_K(G,X)$ in $\sS_{K'}(\GSp,S^\pm)\otimes_{\ZZ_{(p)}}\OO_{E,(p)}$ (where $\OO_{E,(p)}=\OO_E\otimes_\ZZ \ZZ_{(p)}$). Then $\sS_K(G,X)$ is the canonical integral model of $\Sh_K(G,X)$. (By \cite{MR2669706}, the normalization of $\sS_K(G,X)$ is the canonical integral model, and recently it has been shown \cite{xu2020normalization} that the normalization is unnecessary).

Pulling back the universal abelian scheme $\sA'\to\sS_{K'}(\GSp,S^\pm)$ along the map $\sS_K(G,X)\to\sS_{K'}(\GSp,S^\pm)$ we obtain a universal abelian scheme $\cA\to\sS_K(G,X)$ and $p$-divisible group $\cA[p^\infty]$.

\paragraph{}
\label{tensors}
For a vector space or module $W$, we write $W^\otimes$ be the direct sum of all finite combinations of tensor powers, duals, and symmetric and exterior powers of $W$. Since we include duals, we can identify $W^\otimes$ with $(W^\vee)^\otimes$.

As in \cite[2.3.2]{MR2669706}, our group $G$ has a reductive model $G_{\ZZ_{(p)}}$ over $\ZZ_{(p)}$ defined as a subgroup of $\GL(V_{\ZZ_{(p)}})$ by a finite collection of tensors $\{s_\alpha\}\subset V_{\ZZ_{(p)}}^\otimes$. Fix such a collection of tensors $\{s_\alpha\}$.

Our $G$-structures essentially amount to transporting the tensors $\{s_\alpha\}$ to all relevant spaces. Following \cite[1.3.6-10]{MR3630089}, to each point $x\in\sS_K(G,X)(S)$ we associate a finite set of tensors
\[
\{s_{\alpha,\ell,x}\}\subset H^1_\et(\cA_x,\QQ_\ell)^\otimes\cong V_\ell(\cA_x)^\otimes, \qquad\ell\neq p
\]
(here we use the insensitivity of $\bullet^\otimes$ to dualizing); and to each point $x\in\sS_K(G,X)(\oFFp)$ a finite set of tensors
\[
\{s_{\alpha,0,x}\}\subset H^1_\crys(\cA_x/\OO_L)^\otimes\cong \DD(\cA_x[p^\infty])^\otimes
\]
(defined even over $\ZZ_{p^r}$ for $r$ sufficiently large). These tensors are compatible with the original tensors $\{s_\alpha\}$ in the following way. As in \cite[3.2.4]{MR2669706}, for $x\in\sS_K(G,X)(S)$ with associated Siegel data $(\cA_x,\lambda,\eta_{K'}^p)$, the section
\[
\eta_{K'}^p\in\Gamma(S,\underline{\Isom}(V_{\AA_f^p},\hat{V}^p(\cA_x))/K'^p)
\]
can be promoted to a section
\[
\eta_K^p\in\Gamma(S,\underline{\Isom}(V_{\AA_f^p},\hat{V}^p(\cA_x))/K^p),
\]
and this isomorphism $\eta_K^p$ takes $s_\alpha$ to $s_{\alpha,\ell,x}$. At $p$, there is an isomorphism
\[
V^\vee_{\ZZ_{p^r}}\overset{\sim}{\longrightarrow}\DD(\cA_x[p^\infty])(\ZZ_{p^r})
\]
taking $s_\alpha$ to $s_{\alpha,0,x}$.

\begin{proposition}[{\cite[Cor 1.3.11]{MR3630089}}]
\label{tensors distinguish points}
If $x,x'\in\osS_K(G,X)(\oFFp)$ lie over the same point of $\osS_{K'}(\GSp,S^\pm)(\oFFp)$, then $x=x'$ if and only if $s_{\alpha,0,x}=s_{\alpha,0,x'}$ for all $\alpha$.
\end{proposition}

\paragraph{}
\label{par: fix Qpbar embedding}
Fix an embedding $\overline{\QQ}\hookrightarrow\overline{\QQ}_p$, and let $v$ be a prime of $E$ over $p$ determined by this embedding. Denote the residue field by $k(v)$, and let $\osS_K(G,X)=\sS_K(G,X)\otimes_{\OO_{E,(p)}}k(v)$. We denote again by $\cA\to\osS_K(G,X)$ the pullback of the abelian scheme $\cA\to\sS_K(G,X)$.

By \cite{2017arXiv170206611L}, the $p$-divisible group $\cA[p^\infty]$ equipped with polarization and tensors on $\DD(\cA_x[p^\infty])$ gives rise to an isocrystal with $G$-structure over $\osS_K(G,X)$ in the sense of \cite{MR1411570}. It therefore gives rise, as in \ref{siegel newton stratification}, to a Newton stratification
\[
\osS_K^{(\bfb)}(G,X)=\{x\in\osS_K(G,X):\bfb_x=\bfb\}
\]
parametrized by classes $\bfb\in B(G)$. As in Siegel type, these are locally closed subsets, which we equip with the reduced subscheme structure. The $\bfb$-stratum is non-empty exactly when $\bfb\in B(G,\mu^{-1})$, where $\mu$ is a member of the conjugacy class of cocharacters arising from the Shimura datum $(G,X)$ \cite{KMPS}.

\paragraph{}
\label{hodge central leaf}
Fix a class $\bfb\in B(G,\mu^{-1})$ (corresponding to a non-empty Newton stratum) and let $(\Sigma, \lambda_\Sigma, \{s_{\alpha, \Sigma}\})$ be a $p$-divisible group with $G$-structure over $\oFFp$ of type $\bfb$, namely
\begin{itemize}
\item $\Sigma$ a $p$-divisible group over $\oFFp$,
\item $\lambda_\Sigma$ a polarization of $\Sigma$, and
\item $\{s_{\alpha,\Sigma}\}\subset\DD(\Sigma)^\otimes$ a collection of tensors, such that
\item there is an isomorphism $\DD(\Sigma)\to V_{\OO_L}$ preserving the pairings induced by $\lambda_\Sigma$ and $\psi$, taking $s_{\alpha,\Sigma}$ to $s_{\alpha}$, and taking the Frobenius on $D(\Sigma)$ to an endomorphism $b\sigma$ on $V_{\OO_L}$ with $b\in G(\OO_L)$ belonging to the fixed class $\bfb\in B(G)$.
\end{itemize}
Changing the isomorphism changes $b$ by $\sigma$-conjugacy in $G(\OO_L)$, so its class in $B(G)$ is well-defined.

As this data is required for the definition of Igusa varieties of Hodge type, we fix these choices of $\bfb$ and $(\Sigma,\lambda_\Sigma,\{s_{\alpha,\Sigma}\})$ for the remainder of the paper. We fix also an isomorphism $\VV(\Sigma)\to V_L$, and therefore an element $b\in G(L)$ giving the Frobenius $b\sigma$ on $\VV(\Sigma)$, which is a representative of our fixed class $\bfb$.

Define the central leaf corresponding to $(\Sigma,\lambda_\Sigma,\{s_{\alpha,\Sigma}\})$ by
\begin{align*}
C_{\Sigma,K}
=
\Big\{x\in \osS&_K^{(\bfb)}(G,X):
\\
&(\cA_x[p^\infty],\lambda_x,\{s_{\alpha,0,x}\})\otimes_{k(x)}\overline{k(x)}\cong(\Sigma,\lambda_\Sigma,\{s_{\alpha,\Sigma}\})\otimes_{\oFFp}\overline{k(x)}\Big\}.
\end{align*}
This is a closed subset of the Newton stratum, and when equipped with the reduced subscheme structure, is smooth \cite[Prop 2.5]{MR3994311}.

\paragraph{}
\label{hodge igusa varieties}
Following \cite{MR3994311}, we define the perfect infinite level Igusa variety of Hodge type
\[
\sJ_{\Sigma,K}\subset(\Ig'_{\Sigma,K'}\times_{C'_{\Sigma,K'}}C_{\Sigma,K})^{(p^{-\infty})}
\]
to be the locus of points where $j^*(s_{\alpha,0,x})=s_{\alpha,\Sigma}$. That is, $\sJ_{\Sigma,K}$ parametrizes isomorphisms $\Sigma\otimes_{\oFFp}C_{\Sigma,K}\overset{\sim}{\to}\cA[p^\infty]$ preserving polarizations (up to scaling) and tensors. Define (un-perfected) infinite level and finite level Igusa varieties of Hodge type by
\begin{align*}
\Ig_{\Sigma,K}=\im(\sJ_{\Sigma,K}\to\Ig'_{\Sigma,K'}\times_{C'_{\Sigma,K'}}C_{\Sigma,K}), \\
\Ig_{\Sigma,K,m}=\im(\sJ_{\Sigma,K}\to\Ig'_{\Sigma,K',m}\times_{C'_{\Sigma,K'}}C_{\Sigma,K}).
\end{align*}

We will work primarily with the Igusa variety with infinite level at $m$ and infinite level away from $p$, which we call
\[
\Ig_\Sigma=\ilim_{K^p}\Ig_{\Sigma,K_pK^p}
\]
(and analogously $\Ig'_{\Sigma}$ for the Siegel version). Since perfection does not change $\oFFp$-points nor \'etale cohomology, for these purposes it is essentially similar to work with the perfect Igusa variety. An important example is that the natural map $\sJ_{\Sigma,K}\to\sJ'_{\Sigma,K'}$ is a closed embedding \cite[Prop 4.10]{MR3994311}, so we can regard $\Ig_\Sigma(\oFFp)$ as a subset of $\Ig'_{\Sigma}(\oFFp)$.

The Igusa variety $\Ig_\Sigma$ is not an honest moduli space, but we can nonetheless attach useful data to its points, which we will refer to as partial moduli data. A point $x\in\Ig_\Sigma(\oFFp)$ has associated data
\[
(\cA_x,\lambda_x,\eta^p,\{s_{\alpha,0,x}\}, j)
\]
where $\cA_x,\lambda_x,\eta^p,\{s_{\alpha,0,x}\}$ is the data for the image of $x$ in $\osS_{K_p}(G,X)(\oFFp)$---on account of passing to infinite level away from $p$, the level structure $\eta^p$ is now a full trivialization $V_{\AA_f^p}\overset{\sim}{\to}\hat{V}^p(\cA_x)$ sending $s_\alpha$ to $(s_{\alpha,\ell,x})_\ell$---and
\[
j:\Sigma\overset{\sim}{\to}\cA_x[p^\infty]
\]
is the Igusa level structure attached to the image of $x$ in the Siegel Igusa variety $\Ig'_{\Sigma}(\oFFp)$, an isomorphism of $p$-divisible groups over $\oFFp$ respecting polarizations up to scaling and sending $s_{\alpha,\Sigma}$ to $s_{\alpha,0,x}$.

Consider two sets of data to be equivalent if there is a prime-to-p isogeny between the abelian varieties sending one set of data to the other. With this equivalence, points are distinguished by their partial moduli data.

\paragraph{}
The Igusa variety $\Ig_{\Sigma}$ (or system $\Ig_{\Sigma,K,m}$) receives an action of $G(\AA_f^p)$ inherited from the Shimura variety, and a commuting action of the submonoid $S_b\subset J_b(\QQ_p)$, which extends to an action of the full group $J_b(\QQ_p)$ on the perfection and on \'etale cohomology. (Here $b$ is a representative of our fixed class $\bfb\in B(G)$).

Let $\xi$ be a finite-dimensional representation of $G$, and $\sL_\xi$ the system of sheaves on $\Ig_{\Sigma,K,m}$ defined by $\xi$, as in \cite[\S6]{MR1124982}. Define
\begin{gather}
\begin{split}
\label{H_c(Ig)}
H_c^i(\Ig_\Sigma,\sL_\xi)=\dlim_{K^p,m}H^i_c(\Ig_{\Sigma,K,m},\sL_\xi), \\
H_c(\Ig_\Sigma,\sL_\xi)=\sum_i(-1)^iH_c^i(\Ig_\Sigma,\sL_\xi),
\end{split}
\end{gather}
the latter as an element of $\Groth(G(\AA_f^p)\times J_b(\QQ_p))$, where $H^i_c$ denotes \'etale cohomology with compact supports. This is the representation that will be described by our eventual counting point formula.

We can describe the action of $G(\AA_f^p)\times J_b(\QQ_p)$ on the points $\Ig_\Sigma(\oFFp)$ as follows. Let $x\in\Ig_\Sigma(\oFFp)$ be a point with associated partial moduli data $(\cA_x,\lambda_x,\eta^p,\{s_{\alpha,0,x}\},j)$. Note that we can regard $J_b(\QQ_p)$ as the group of self-quasi-isogenies of $\Sigma$.

The action of $G(\AA_f^p)$ is inherited from the Shimura variety, and as there, it acts on the level structure $\eta^p$: for $g^p\in G(\AA_f^p)$, the data associated to $x\cdot g^p$ is
\[
(\cA_x,\lambda_x,\eta^p\circ g^p, \{s_{\alpha,0,x}\},j).
\]

To describe the action of $g_p\in J_b(\QQ_p)$, regard it as a quasi-isogeny $g_p:\Sigma\to\Sigma$, and choose $m\ge0$ such that $p^mg_p^{-1}:\Sigma\to\Sigma$ is an isogeny. The Igusa level structure $j:\Sigma\overset{\sim}{\to}\cA_x[p^\infty]$ allows us to transfer this to $\cA_x$. The data associated to $x\cdot g_p$ is
\[
(\cA_x/j(\ker p^mg_p^{-1}), g_p^*\lambda_x, \eta^p, \{s_{\alpha,0,x}\}, g_p^*j)
\]
where $g_p^*\lambda_x$ is the induced polarization; we can take the same level structure $\eta^p$ because $\cA_x$ is unchanged away from $p$, and the same tensors $\{s_{\alpha,0,x}\}$ because $J_b(\QQ_p)$, being a subgroup of $G(L)$, preserves tensors; and
\[
g_p^*j:\Sigma\overset{\sim}{\longrightarrow}\cA_{x\cdot g_p}[p^\infty]=\cA_x[p^\infty]/j(\ker p^mg_p^{-1})
\]
is the unique map making the following diagram commute.
\[
\xymatrix{
\Sigma \ar[r]^j \ar[d]_{p^mg_p^{-1}} & \cA_x[p^\infty] \ar[d] \\
\Sigma \ar[r]_{g_p^*j\hspace{4.5em}} & \cA_{x}[p^\infty]/j(\ker p^m g_p^{-1})
}
\]
Note that the choice of $m$ does not matter because multiplication by $p^k$ induces an isomorphism $A/\ker p^k\to A$, and this gives an equivalence between moduli data for different choices of $m$.

\paragraph{}
\label{alternative partial moduli data}
There is an alternative partial moduli description which admits a simpler description of the group actions, but has the downside that it makes the map to $\osS_{K_p}(G,X)(\oFFp)$ more opaque. By \cite[Lem 7.1]{MR2484281}, we have the following moduli description of the Siegel Igusa variety:
\[
\Ig'_\Sigma(\oFFp)=\{(A,\lambda,\eta^p,j)\}/\sim
\]
where
\begin{itemize}
\item $A$ is an abelian variety over $\oFFp$,
\item $\lambda$ is a polarization of $A$,
\item $\eta^p:V_{\AA_f^p}\overset{\sim}{\to}\hat{V}^p(A)$ is an isomorphism preserving the pairings induced by $psi$ and $\lambda$ up to scaling,
\item $j:\Sigma\to A[p^\infty]$ is a quasi-isogeny preserving polarizations up to scaling, and
\item two tuples are equivalent if there is an isogeny $A_1\to A_2$ sending $\lambda_1$ to a scalar multiple of $\lambda_2$, and sending $\eta^p_1$ to $\eta^p_2$ and $j_1$ to $j_2$. (It is equivalent to replace ``isogeny'' here with ``quasi-isogeny'').
\end{itemize}
Note the difference that we allow $j$ to be a quasi-isogeny rather than an isomorphism, and equivalence requires only an isogeny $A_1\to A_2$, rather than a prime-to-$p$ isogeny.

Under this moduli description, $\Ig'_\Sigma(\oFFp)$ has a right action of $\GSp(\AA_f^p)\times J_{b}^{\GSp}(\QQ_p)$ (where we write $b$ again for the image of $b$ in $\GSp(L)$) described by
\[
(g^p,g_p):(A,\lambda,\eta^p,j)\mapsto (A,\lambda,\eta^p\circ g^p,j\circ g_p).
\]
As noted in \ref{hodge igusa varieties}, we can regard $\Ig_\Sigma(\oFFp)\subset\Ig'_{\Sigma}(\oFFp)$, and furthermore this is compatible with the actions of $G(\AA_f^p)\times J_b^G(\QQ_p)\subset \GSp(\AA_f^p)\times J_b^{\GSp}(\QQ_p)$. Thus each point of $\Ig_\Sigma(\oFFp)$ can be associated data $(A,\lambda,\eta^p,j)$ as above, with distinct points having distinct data, and we can write the action of $G(\AA_f^p)\times J_b^G(\QQ_p)$ in a precisely similar way.

\subsection{Galois Gerbs}
\label{section: galois gerbs}

In this section we review Galois gerbs. We refer to \S3 of \cite{MR3630089} and \S2 of \cite{KSZ} for details omitted here.

\paragraph{}
A \emph{$k'/k$-Galois gerb} is a linear algebraic group $G$ over $k'$ and an extension of topological groups (giving $G(k')$ the discrete topology)
\[
1\to G(k')\to\fG\to\Gal(k'/k)\to 1
\]
satisfying certain technical conditions \cite[3.1.1]{MR3630089}. We will use the name of the extension $\fG$ also to refer to the whole data, and we say that $\fG^\Delta=G$ is the \emph{kernel} of $\fG$. We sometimes refer to a $\overline{k}/k$-Galois gerb as simply a Galois gerb over $k$.

A $k'/k$-Galois gerb $\fG$ induces a $\overline{k}/k$-Galois gerb via pullback by $\Gal(\overline{k}/k)\to\Gal(k'/k)$ and pushout by $G(k')\to G(\overline{k})$. Similarly, for any place $v$ of $\QQ$, a Galois gerb $\fG$ over $\QQ$ induces a Galois gerb $\fG(v)$ over $\QQ_v$ by pulling back by $\Gal(\overline{\QQ}_v/\QQ_v)\to\Gal(\overline{\QQ}/\QQ)$ and pushing out by $G(\overline{\QQ})\to G(\overline{\QQ}_v)$.

An important example will be the \emph{neutral $k'/k$-Galois gerb} attached to a group $G$ over $k$, defined to be the semi-direct product $\fG_G=G(k')\rtimes\Gal(k'/k)$, where the action of $\Gal(k'/k)$ on $G(k')$ is given by the $k$-structure on $G$.

A morphism of $k'/k$-Galois gerbs $f:\fG_1\to\fG_2$ is a continuous group homomorphism which fits into a commuting diagram
\[
\xymatrix{
1 \ar[r] & G_1(k') \ar[d]_{f^\Delta(k')} \ar[r] & \fG_1 \ar[d]_{f} \ar[r] & \Gal(k'/k) \ar[d]_{\id} \ar[r] & 1 \\
1 \ar[r] & G_2(k') \ar[r] & \fG_2 \ar[r] & \Gal(k'/k) \ar[r] & 1
}
\]
where $f^\Delta(k')$ is induced by a map of algebraic groups $f^\Delta:G_1\to G_2$. If $\fG_1,\fG_2$ are Galois gerbs over $\QQ$, then $f$ induces a morphism $f(v):\fG_1(v)\to\fG_2(v)$ of Galois gerbs over $\QQ_v$. Two morphism $\fG_1\to\fG_2$ are \emph{conjugate} if they are related by conjugation by an element of $G_2(k')$. We denote the conjugacy class of a morphism $f$ by $[f]$.

Let $f:\fG_1\to\fG_2$ be a morphism. If $R$ is a $k$-algebra, we can consider the pushouts $\fG_{1,R}, \fG_{2,R}$ of $\fG_1,\fG_2$ by the map $G_1(k')\to G_1(k'\otimes_kR)$, and $f$ induces a map $\fG_{1,R}\to\fG_{2,R}$; we define the automorphism group of the morphism $f$ by
\[
I_f(R)=\{g\in G_2(k'\otimes_kR):\Int(g)\circ f=f\}.
\]
This defines a group scheme $I_f$ over $k$.

In the case that $\fG_1$ is any Galois gerb and $\fG_2=\fG_G$ is the neutral Galois gerb attached to a linear algebraic group $G$, we have the following lemma.

\begin{lemma}[{\cite[Lem 3.1.2]{MR3630089}}]
\label{lem: gerb kernel facts}
Let $f:\fG_1\to\fG_G$ be a map of $k'/k$-Galois gerbs.
\begin{itemize}
\item The map $I_{f,k'}\to G_{k'}$ given by
\[
I_{f,k'}(R)\longrightarrow G(k'\otimes_k R)\to G(R)
\]
(where $R$ is a $k'$-algebra, and $k'\otimes_k R\to R$ is the multiplication map) identifies $I_{f,k'}$ with the centralizer $Z_G(f^\Delta)$ in $G_{k'}$.\
\item The set of morphisms $f':\fG_1\to\fG_G$ with $f'^\Delta=f^\Delta$ is in bijection with $Z^1(\Gal(k'/k),I_f(k'))$, via the map sending $e\in Z^1(\Gal(k'/k),I_f(k'))$ to the morphism $ef$ defined such that, if $f(q)=g\rtimes\rho$, we have $ef(q)=e_\rho g\rtimes\rho$. Furthermore, $ef$ is conjugate to $e'f$ exactly when $e$ is cohomologous to $e'$.
\end{itemize}
\end{lemma}

\paragraph{}
\label{quasi-motivic gerb}
There is a distinguished Galois gerb over $\QQ$ called the \emph{quasi-motivic Galois gerb}, and denoted $\fQ$, which plays a central role in point counting. We refer to \cite[3.1]{MR3630089} for full details, but we review here the essential properties we will need.

For $L/\QQ$ a finite Galois extension, define
\[
Q^L=(Res_{L(\infty)/\QQ}\GG_m \times \Res_{L(p)/\QQ}\GG_m)/\GG_m,
\]
where the action of $\GG_m$ is the diagonal action, and $L(\infty)=L\cap\RR$ and $L(p)=L\cap\QQ_p$. This group is equipped with cocharacters $\nu(p)^L$ over $\QQ_p$ and $\nu(\infty)^L$ over $\RR$ defined by
\[
\nu(v)^L:\GG_m\to\Res_{L(v)/\QQ}\GG_m\to Q^L.
\]
For $L'/L$ Galois there is a natural map $Q^{L'}\to Q^L$, and the limit is a pro-torus $Q=\ilim_L Q^L$ over $\QQ$ equipped with a fractional cocharacter $\nu(p):\DD\to Q$ over $\QQ_p$ and cocharacter $\nu(\infty):\GG_m\to Q$ over $\RR$. The kernel of the quasi-motivic Galois gerb is this pro-torus $\fQ^\Delta=Q$.

The quasi-motivic Galois gerb $\fQ$ has a $v$-adic realization for any place $v$ of $\QQ$ in the form of a morphism $\zeta_v:\fG_v\to\fQ(v)$ from a distinguished Galois gerb over $\QQ_v$.

At $\ell\neq p,\infty$
\[
\fG_\ell=\Gal(\overline{\QQ}_\ell/\QQ_\ell)
\]
is the trivial Galois gerb.

At $p$, we have
\[
\fG_p=\ilim_L \fG_p^L
\]
where $L$ runs over finite Galois extensions of $\QQ_p$, and $\fG_p^L$ is the $L/\QQ_p$-gerb (induced to $\overline{\QQ}_p$) with kernel $\fG_p^{L,\Delta}=\GG_m$ given by the fundamental class in $H^2(\Gal(L/\QQ_p),L^\times)$. The kernel of $\fG_p$ is $\fG_p^\Delta=\DD$, the pro-torus with character group $\QQ$.

At $\infty$,
\[
1\to\CC^\times\to\fG_\infty\to\Gal(\CC/\RR)\to 1
\]
is the extension corresponding to the fundamental class in $H^2(\Gal(\CC/\RR),\CC^\times)$.

At $p$ and $\infty$ we have $\zeta_p^\Delta=\nu(p)$ and $\zeta_\infty^\Delta=\nu(\infty)$; these cocharacters will play a role in some later arguments.

The quasi-motivic Galois gerb is also equipped with a distinguished morphism $\psi:\fQ\to\fG_{\Res_{\overline{\QQ}/\QQ}\GG_m}$ which allows us to construct morphisms to neutral Gerbs from cocharacters in the following way. Let $T$ be a torus over $\QQ$ and $\mu$ a cocharacter of $T$ defined over a Galois extension $L/\QQ$. Then $\mu$ induces a map $\Res_{L/\QQ}\GG_m\to T$ which further induces a morphism of Galois gerbs $\fG_{\Res_{L/\QQ}}\to\fG_T$. We define a morphism $\psi_\mu:\fQ\to\fG_T$ by the composition
\[
\psi_\mu:\fQ\overset{\psi}{\longrightarrow}\fG_{\Res_{\overline{\QQ}/\QQ}}\longrightarrow\fG_{\Res_{L/\QQ}}\longrightarrow \fG_T.
\]

\paragraph{}
Recall that $\fG_p$ was defined as a limit of Galois gerbs $\fG_p^L$, where $L$ runs over finite Galois extensions of $\QQ_p$. If we restrict $L$ to run unramified extensions, we can define a $\QQ_p^\ur/\QQ_p$-Galois gerb $\fD$ with kernel $\fD^\Delta=\DD$, which becomes $\fQ_p$ when induced to $\overline{\QQ}_p$. Writing $\sigma\in\Gal(\QQ_p^\ur/\QQ_p)$ for the Frobenius, there is a distinguished element $d_\sigma\in \fD$ lying over $\sigma$ and such that $d_\sigma^n$ maps to $p^{-1}\in \GG_m=\fG_p^{\QQ_{p^n},\Delta}$ under the projection to $\fG_p^{\QQ_{p^n}}$.

\begin{definition}
\label{unramified morphism}
A morphism $\theta:\fG_p\to\fG_G$ is \emph{unramified} if it is induced by a morphism $\theta^\ur:\fD\to\fG_G^\ur$, where $\fG_G^\ur$ is the neutral $\QQ_p^\ur/\QQ_p$-gerb attached to $G$. For $\theta$ an unramified morphism, we define an element $b_\theta\in G(\QQ_p^\ur)$ by $\theta^\ur(d_\sigma)=b_\theta\rtimes\sigma$.
\end{definition}

\paragraph{}
\label{unramified paragraph}
If $G$ is connected, every morphism $f:\fG_p\to\fG_G$ is conjugate to an unramified morphism \cite[Lem 2.2.4 (i)]{KSZ}. If $\theta$ and $\theta'$ are unramified morphisms conjugate to $f$, then $b_\theta$ and $b_{\theta'}$ are related by $\sigma$-conjugacy in $G(\QQ_p^\ur)$, so we can associate to $f$ a well-defined class $[b_\theta]\in B(G)$.

\begin{lemma}[{\cite[Prop 2.2.6]{KSZ}}]
\label{unramified morphism facts}
Let $G$ be a connected linear algebraic $\QQ_p$-group, $\theta:\fG_p\to\fG_G$ an unramified morphism, and $\nu$ the fractional cocharacter $\theta^{\ur,\Delta}:\DD_{\QQ_p^\ur}\to G_{\QQ_p^\ur}$. Then
\begin{itemize}
\item $\nu=-\nu_{b_\theta}$, where $\nu_{b_\theta}$ is the slope homomorphism of $b_\theta$, and
\item there are natural $\QQ_p$-isomorphisms $J_{b_\theta}\cong I_{\theta^\ur}\cong I_\theta$.
\end{itemize}
\end{lemma}

\paragraph{}
The morphisms $\fQ\to\fG_G$ that will be used in our point-counting are required to satisfy an admissibility condition. For $\ell\neq p,\infty$, let $\xi_\ell:\fG_\ell\to\fG_G(\ell)$ be the map sending $\rho\mapsto 1\rtimes \rho$.

\begin{definition}
\label{admissible morphism}
A morphism $\phi:\fQ\to\fG_G$ is \emph{admissible} if
\begin{itemize}
\item for $\ell\neq p,\infty$, the morphism $\phi(\ell)\circ\zeta_\ell:\fG_\ell\to\fG_G(\ell)$ is conjugate to $\xi_\ell$;
\item at $p$, the morphism $\phi(p)\circ\zeta_p:\fG_p\to\fG_G(p)$ is conjugate to an unramified morphism $\theta$ such that $b_\theta\in G(\ZZ_p^\ur)\mu(p)^{-1}G(\ZZ_p^\ur)$;
\end{itemize}
as well as satisfying a global condition \cite[3.3.6 (1)]{MR3630089}.
\end{definition}

An important example of admissible morphisms is given by the morphisms $\psi_\mu$ constructed in \ref{quasi-motivic gerb}. Namely, if $T\subset G$ is a torus and $\mu\in X_*(T)$, we have an admissible morphism for $G$ defined by
\[
i\circ\psi_\mu:\fQ\overset{\psi_\mu}{\longrightarrow}\fG_T\overset{i}{\longrightarrow}\fG_G.
\]

\paragraph{}
\label{par: twisting admissible morphisms}
As in Lemma \ref{lem: gerb kernel facts} we can twist a morphism $\phi:\fQ\to\fG_G$ by a cocycle $e\in Z^1(\QQ,I_\phi)$.

Define
\[
\Sha_G^S(\QQ,I_\phi)\subset H^1(\QQ,I_\phi)
\]
to be the subset of classes which are trivial in $H^1(\QQ,G)$ as well as $H^1(\QQ_v,I_\phi)$ for $v\in S$. For our purposes $S$ will be $\{\infty\}$ or $\{p,\infty\}$ or $\{\text{all places of $\QQ$}\}$. In the case that $S$ is all places of $\QQ$, we also write $\Sha_G(\QQ,I_\phi)$ for $\Sha_G^S(\QQ,I_\phi)$. See \cite[1.2.5]{KSZ} for more details.

\begin{proposition}[{\cite[Prop 2.6.12]{KSZ}}]
\label{prop: admissible twist}
If $\phi$ is an admissible morphism and $e\in Z^1(\QQ,I_\phi)$, then $e\phi$ is admissible exactly when $e$ lies in $\Sha_G^\infty(\QQ,I_\phi)$.
\end{proposition}

\paragraph{}
\label{par: X^p(phi) and cocycles}
For an admissible morphism $\phi$ we can define a set
\[
X^p(\phi)=\{x=(x_\ell)\in G(\overline{\AA}_f^p):\Int(x_\ell)\circ\xi_\ell=\phi(\ell)\circ\zeta_\ell\}.
\]
It is non-empty by the admissible condition for $\ell\neq p,\infty$, and furthermore is a $G(\AA_f^p)$-torsor under the natural right action (note that $I_{\xi_\ell}(\QQ_\ell)=G(\QQ_\ell)$).

Let $x\in X^p(\phi)$, and define a cocycle $\zeta_\phi^{p,\infty}:\Gal(\overline{\QQ}/\QQ)\to G(\overline{\AA}_f^p)$ by $\rho\mapsto x\rho(x)^{-1}$. This does not depend on the choice of $x$, because any other choice $x'$ is related by $x'=xg$ for some $g\in G(\AA_f^p)$, and since $g$ is rational we have $x'\rho(x')^{-1}=xg\rho(xg)^{-1}=xgg^{-1}\rho(x)^{-1}=x\rho(x)^{-1}$.

Define cocycles $\zeta_{\phi,\ell}:\Gal(\overline{\QQ}_\ell/\QQ_\ell)\to G(\overline{\QQ}_\ell)$ by
\[
\zeta_{\phi,\ell}:\Gal(\overline{\QQ}_\ell/\QQ_\ell)\to\Gal(\overline{\QQ}/\QQ)\overset{\zeta_\phi^{p,\infty}}{\longrightarrow}G(\overline{\AA}_f^p)\to G(\overline{\QQ}_\ell),
\]
where the map $G(\overline{\AA}_f^p)\to G(\overline{\QQ}_\ell)$ is induced by the map $\overline{\AA}_f^p\to\overline{\QQ}_\ell$ sending a point to its $\ell$-component. These cocycles record the component of $\phi$ at $\ell$ in the sense that $(\phi(\ell)\circ\zeta_\ell)(\rho)=\zeta_{\phi,\ell}(\rho)\rtimes\rho$ for all $\rho\in\Gal(\overline{\QQ}_\ell/\QQ_\ell)$.

\section{Langlands--Rapoport Conjecture for Igusa Varieties of Hodge Type}
\label{igusa LR section}

\subsection{Isogeny Classes on Shimura Varieties of Hodge Type}
\label{isogeny classes on shimura varieties}

\paragraph{}
\label{ADLV}
Following \cite[1.4.1]{MR3630089} we define a cocharacter $v$ of $G$ and an element $b\in G(L)$. These are needed to define the affine Deligne-Lusztig variety $X_v(b)$, which records the $p$-part of an isogeny class on the Shimura variety.

Let $x\in \osS_{K}(G,X)(\FF_{p^r})$ for $r$ divisible by the residue degree of $v$ (our fixed prime of $E$ lying over $p$ as in \ref{par: fix Qpbar embedding}). As in \ref{tensors}, there is an isomorphism
\[
V_{\ZZ_{p^r}}^\vee\overset{\sim}{\longrightarrow}\DD(\cA_x[p^\infty])(\ZZ_{p^r})
\]
taking $s_\alpha$ to $s_{\alpha,0,x}$. This allows us to identify
\begin{equation}
\label{eq: G acting on D}
G_{\ZZ_{(p)}}\otimes_{\ZZ_{(p)}}\ZZ_{p^r}\overset{\sim}{\to}G_{0,x}\subset\GL(\DD(\cA_x[p^\infty])(\ZZ_{p^r})),
\end{equation}
where $G_{0,x}$ is the subgroup fixing the tensors $s_{\alpha,0,x}$. Under this identification, the action of Frobenius is given by $b\sigma$, where $b\in G(\QQ_{p^r})$ and $\sigma$ is the lift of Frobenius on $\ZZ_{p^r}$. The element $b$ is well-defined up to $\sigma$-conjugacy by $G(\ZZ_{p^r})$.

Fix a Borel $B$ and maximal torus $T$ in $G_{\ZZ_{(p)}}$, let $\mu\in X_*(T)$ be the dominant cocharacter in the conjugacy class of $\mu_h$ for $h\in X$, and let $v=\sigma(\mu^{-1})$. Then we have $b\in G(\ZZ_{p^r})v(p)G(\ZZ_{p^r})$.

Define the \emph{affine Deligne-Lusztig variety}
\[
X_v(b)=\{g\in G(L)/G(\OO_L):g^{-1}b\sigma(g)\in G(\OO_L)v(p)G(\OO_L)\}.
\]
We consider $X_v(b)$ as a set, and equip it with a Frobenius operator
\[
\Phi(g)=(b\sigma)^rg=b\cdot\sigma(b)\cdots\sigma^{r-1}(b)\cdot\sigma^r(g).
\]
Note that it also carries an action of $J_b(\QQ_p)$ by left multiplication, as $J_b(\QQ_p)$ is the group of elements of $G(L)$ which $\sigma$-centralize $b$.

\paragraph{}
\label{ADLV map description}
Following \cite[1.4.2]{MR3630089} we define a map $X_v(b)\to\osS_{K_p}(G,X)(\oFFp)$ as follows. Choose a base point $x\in \osS_{K_p}(G,X)(\oFFp)$, with associated $p$-divisible group $\cA_x[p^\infty]$. For $g\in X_v(b)$, the lattice $g\cdot\DD(\cA_x[p^\infty])\subset \VV(\cA_x[p^\infty])$ is again a Dieudonn\'e module, and corresponds to a $p$-divisible group $\sG_{gx}$ equipped with a quasi-isogeny $\cA_x[p^\infty]\to\sG_{gx}$.

Let $\cA_{gx}$ be the corresponding abelian variety equipped with the polarization and level structure induced from $\cA_x$. Sending $g\mapsto \cA_{gx}$ with polarization and level structure defines a map $X_v(b)\to\osS_{K'_p}(\GSp,S^\pm)(\oFFp)$. By \cite[Prop 1.4.4]{MR3630089} this map has a unique lift to a map $i_x:X_v(b)\to\osS_{K_p}(G,X)(\oFFp)$ satisfying $s_{\alpha,0,x}=s_{\alpha,0,i_x(g)}\in \DD(\cA_{gx}[p^\infty])$. Extending by the action of $G(\AA_f^p)$, we get a map
\begin{equation}
\label{shimura isogeny class map}
i_x:G(\AA_f^p)\times X_v(b)\longrightarrow\osS_{K_p}(G,X)(\oFFp)
\end{equation}
which is equivariant for the action of $G(\AA_f^p)$ and intertwines the action of $\Phi$ on $X_v(b)$ with the action of geometric $p^r$-Frobenius on $\osS_{K_p}(G,X)$.

\begin{definition}
For $x\in\osS_{K_p}(G,X)(\oFFp)$, the \emph{isogeny class of $x$}, denoted $\sI_x^\Sh$, is the image of the map \eqref{shimura isogeny class map}.
\end{definition}

\paragraph{}
\label{automorphism group at x}
Define $\Aut_{\QQ}(\cA_x)$ to be the algebraic group over $\QQ$ with points
\[
\Aut_{\QQ}(\cA_x)(R)=(\End_{\QQ}(\cA_x)\otimes_\QQ R)^\times,
\]
and define $I_x\subset\Aut_{\QQ}(\cA_x)$ to be the subgroup preserving the polarization of $\cA_x$ up to scaling and fixing the tensors $s_{\alpha,\ell,x}$ $(\ell\neq p)$ and $s_{\alpha,0,x}$. 

The level structure $\eta^p:V_{\AA_f^p}\overset{\sim}{\longrightarrow}\hat{V}^p(\cA_x)$ away from $p$ identifies the tensors $s_\alpha$ and $(s_{\alpha,\ell,x})_{\ell\neq p}$, and therefore identifies $G(\AA_f^p)$ with the subgroup of $\GL(\hat{V}^p(\cA_x))$ fixing $(s_{\alpha,\ell,x})_{\ell\neq p}$. Thus the embedding $\Aut_{\QQ}(\cA_x)(\QQ)\hookrightarrow\GL(\hat{V}^p(\cA_x))$ induces an embedding $I_x(\QQ)\hookrightarrow G(\AA_f^p)$, canonical up to conjugation by $G(\AA_f^p)$.

Similarly, our above identification \eqref{eq: G acting on D} allows us to identify $J_b(\QQ_p)$ with the subgroup of $\GL(\VV(\cA_x[p^\infty]))$ fixing the tensors $s_{\alpha,0,x}$ and commuting with the Frobenius. Thus the embedding $\Aut_{\QQ}(\cA)(\QQ)\hookrightarrow\GL(\VV(\cA_x[p^\infty]))$ induces an embedding $I_x(\QQ)\hookrightarrow J_b(\QQ_p)$, canonical up to conjugation by $J_b(\QQ_p)$.

Thus we have an embedding
\[
I_x(\QQ)\hookrightarrow G(\AA_f^p)\times J_b(\QQ_p),
\]
canonical up to conjugation. We fix such a choice of embedding.

Through this embedding, $I_x(\QQ)$ acts on the set $G(\AA_f^p)\times X_v(b)$. By \cite[Prop 2.1.3]{MR3630089}, the map \eqref{shimura isogeny class map} induces an injective map
\begin{equation}
i_x:I_x(\QQ)\backslash G(\AA_f^p)\times X_v(b)\hookrightarrow\osS_{K_p}(G,X)(\oFFp).
\end{equation}
Thus the isogeny class of a point $x\in\osS_{K_p}(G,X)(\oFFp)$ is parametrized by the set $I_x(\QQ)\backslash G(\AA_f^p)\times X_v(b)$.

We can also give a description of isogeny classes in terms of the partial moduli structure, which (in addition to being useful) gives a very plain relation to isogenies of the moduli data.

\begin{proposition}[{\cite[Prop 1.4.15]{MR3630089}}]
\label{isogeny moduli description}
Two points $x,x'\in\osS_{K_p}(G,X)(\oFFp)$ lie in the same isogeny class exactly when there is a quasi-isogeny $\cA_x\to\cA_{x'}$ preserving polarizations up to scaling and such that the induced maps $\DD(\cA_{x'}[p^\infty])\to\DD(\cA_x[p^\infty])$ and $\hat{V}^p(A_x)\to\hat{V}^p(\cA_{x'})$ send $s_{\alpha,0,x'}$ to $s_{\alpha,0,x}$ and $s_{\alpha,\ell,x}$ to $s_{\alpha,\ell,x'}$.
\end{proposition}

\subsection{Isogeny Classes on Igusa Varieties of Hodge Type}
\label{isogeny classes on igusa varieties}

\paragraph*{}
Recall from \ref{automorphism group at x} the group $I_x$ and embedding $I_x(\QQ)\hookrightarrow G(\AA_f^p)\times J_b(\QQ_p)$. We begin by giving the parametrization (and definition) of an isogeny class in $\Ig_\Sigma(\oFFp)$.

\begin{lemma}
\label{ICmap}
For $x\in \Ig_\Sigma(\oFFp)$, the map
\begin{align*}
i_x:I_x(\QQ)\backslash G(\AA_f^p)\times J_b(\QQ_p) & \to\Ig_\Sigma(\oFFp) \\
(g^p,g_p)& \mapsto x\cdot (g^p,g_p)
\end{align*}
is well-defined and injective.
\end{lemma}
\begin{proof}
To show that this map is well-defined and injective is to show that $I_x(\QQ)$ is the stabilizer of $x$ under the action of $G(\AA_f^p)\times J_b(\QQ_p)$. We use the partial moduli structure described in \ref{alternative partial moduli data}, as it provides a simpler description of the group action. Let $(A,\lambda,\eta^p,j)$ be the data associated to $x$. Then the data associated to $x\cdot(g^p,g_p)$ is $(A,\lambda,\eta^p\circ g^p,j\circ g_p)$.

If $x=x\cdot(g^p,g_p)$, then by the partial moduli interpretation there is a quasi-isogeny $\theta:A\to A$ preserving $\lambda$ up to $\QQ^\times$-scaling and sending $\eta^p$ to $\eta^p\circ g^p$ and sending $j$ to $j\circ g_p$. This implies that $\theta$ acts as $g^p$ on $\hat{V}^p(A)$, and therefore preserves the tensors $s_{\alpha,\ell,x}$; also $\theta$ acts as $g_p$ on $\VV(\cA[p^\infty])$, and therefore preserves the tensors $s_{\alpha,0,x}$. Thus we can regard $\theta$ as an element of $I_x(\QQ)$, which is identified with $(g^p,g_p)$ under our embedding; i.e. $(g^p,g_p)=\theta\in I_x(\QQ)$.

Conversely, suppose that $(g^p,g_p)=\theta\in I_x(\QQ)$. Then $\theta:A\to A$ is a quasi-isogeny preserving the polarization up to $\QQ^\times$ and sending $\eta^p$ to $\eta^p\circ g^p$ and sending $j$ to $j\circ g_p$. Thus $\theta$ gives an equivalence between $(A,\lambda,\eta^p,j)$ and $(A,\lambda,\eta^p\circ g^p,j\circ g_p)$ in the moduli description, showing $x=x\cdot(g^p,g_p)$.
\end{proof}

\begin{definition}
\label{igusa isogeny classes}
For $x\in\Ig_\Sigma(\oFFp)$, the \emph{isogeny class of $x$}, denoted $\sI_x^\Ig$, is the image of the map of \ref{ICmap}. We may also speak of an \emph{isogeny class} $\sI^{\Ig}$ if we do not specify $x$.
\end{definition}

The next two lemmas relate isogeny classes on the Igusa variety and the Shimura variety. The essential relationship is stated in Proposition \ref{isogeny class comparison}.

\begin{lemma}
\label{pullback}
Let $x\in\Ig_\Sigma(\oFFp)$, and $x'\in\osS_{K_p}(G,X)(\oFFp)$ the image of $x$ under the natural map. The isogeny class maps for $x$ and $x'$ fit into a pullback diagram as below, where the left vertical map is defined by $(g^p,g_p)\mapsto (g^p,g_p)$.
\[
\xymatrix{
I_x(\QQ)\backslash G(\AA_f^p)\times J_b(\QQ_p) \ar[r]^{\hspace{3em}i_x} \ar[d] & \Ig_\Sigma(\oFFp) \ar[d] \\
I_{x'}(\QQ)\backslash G(\AA_f^p)\times X_v(b) \ar[r]_{\hspace{2em}i_{x'}} & \osS_{K_p}(G,X)(\oFFp)
}
\]
That is, an isogeny class on the Igusa variety is the preimage of an isogeny class on the Shimura variety.
\end{lemma}
\begin{proof}
Recall from \ref{ADLV} that we have chosen $b$ so that $b\in G(\ZZ_{p^r})v(p)G(\ZZ_{p^r})$, which ensures that the coset of $1\in G(L)$ is in $X_v(b)$. By $(g^p,g_p)\mapsto(g^p,g_p)$ we mean mapping $G(\AA_f^p)\to G(\AA_f^p)$ by the identity, and mapping $J_b(\QQ_p)\to X_v(b)$ by sending the identity to the identity coset and extending by the action of $J_b(\QQ_p)$ on $X_v(b)$ by left multiplication.

First we show that this diagram commutes. For this we can ignore the quotients by $I_x(\QQ)$. The whole diagram is $G(\AA_f^p)$-equivariant, so it suffices to check commutativity for elements of $J_b(\QQ_p)$. Let $(1,g_p)\in G(\AA_f^p)\times J_b(\QQ_p)$. We use the partial moduli description of \ref{hodge igusa varieties}, because it provides a simpler description of the map to the Shimura variety. Write $(A_x,\lambda,\eta^p,\{s_{\alpha,0,x}\},j)$ for the data at $x$. Then the data at $x\cdot g_p$, i.e. the image of $(1,g_p)$ in $\Ig_\Sigma(\oFFp)$, is
\[
(A_x/j(\ker p^mg_p^{-1}),g_p^*\lambda,\{s_{\alpha,0,x}\},\eta^p, g_p^*j),
\]
and the data at the image in $\osS_{K_p}(G,X)(\oFFp)$ (following the right vertical arrow of the diagram) is the same, simply forgetting $g_p^*j$.

Going the other way around the diagram, the image in $G(\AA_f^p)\times X_v(b)$ is $(1,g_p)$. The image in $\osS_{K_p}(G,X)(\oFFp)$ (following the bottom horizontal arrow) is defined as in \ref{ADLV map description} by taking the $p$-divisible group $\sG_{g_px}$ associated to the Dieudonn\'e module $g_p\cdot\DD(\cA_x[p^\infty])$, with quasi-isogeny $\cA_x[p^\infty]\to\sG_{g_px}$ induced by the isomorphism $g_p^{-1}:g_p\cdot\VV(\cA_x[p^\infty])\overset{\sim}{\to}\VV(\cA_x[p^\infty])$, then the corresponding abelian variety $A_{g_px}$ with induced polarization and level structure, and the same tensors $s_{\alpha,0,x}$ (as usual note $g_p$ preserves tensors).

Since $g_p$ is (the image of) an element of $J_b(\QQ_p)$, which consists of \emph{self}-quasi-isogenies, we see that $\sG_{g_px}$ is isomorphic to $\cA_x[p^\infty]$, and the quasi-isogeny $\cA_x[p^\infty]\to\sG_{g_px}$ corresponds via $j$ to $g_p^{-1}$. Thus, taking $m$ large enough that $p^mg_p^{-1}$ is an isogeny, we can identify $\sG_{g_px}$ with $\cA_x[p^\infty]/j(\ker p^mg_p^{-1})$ and $\cA_{g_px}$ with $\cA_x/j(\ker p^mg_p^{-1})$, with the induced polarization and away-from-$p$ level structure, the same tensors, and the Igusa level structure $g_p^*j$. This matches the data produced by traversing the diagram the other way, and we see that the diagram commutes.

To show the diagram is a pullback, let $x_1\in\Ig_\Sigma(\oFFp)$ be any point whose image $x_1'$ in $\osS_{K_p}(G,X)(\oFFp)$ is in the isogeny class of $x'$. We want to show that $x_1$ is in the isogeny class of $x$.

Since $x_1'$ and $x'$ lie in the same isogeny class, they are related by a pair $(g^p,g_0)\in G(\AA_f^p)\times X_v(b)$. In particular, the $p$-divisible groups are related by a quasi-isogeny $\cA_{x}[p^\infty]\to\cA_{x_1}[p^\infty]$ corresponding to the isomorphism
\[
\VV(\cA_{x_1}[p^\infty])\cong g_0\cdot\DD(\cA_x[p^\infty])\otimes_{\ZZ_p}\QQ_p\overset{\sim}{\to}\VV(\cA_x[p^\infty]).
\]
Using the Igusa level structures $j$ at $x$ and $j_1$ at $x_1$, we can translate this to a quasi-isogeny $\Sigma\to\Sigma$ given by an element $g_p^{-1}\in J_b(\QQ_p)$ (i.e. we define $g_p$ to be the inverse of this quasi-isogeny). We claim that $x_1$ is related to $x$ by $(g^p,g_p)\in G(\AA_f^p)\times J_b(\QQ_p)$.

Indeed, $g_p$ maps to $g_0\in X_v(b)$ because by construction they send $\DD(\cA_x[p^\infty])$ to the same lattice $g_0\cdot\DD(\cA_x[p^\infty])=g_p\cdot\DD(\cA_x[p^\infty])$ in $\VV(\cA_x[p^\infty])$. Thus $x_1$ and $x\cdot (g^p,g_p)$ have the same image in $\osS_{K_p}(G,X)(\oFFp)$, so it only remains to show they have the same Igusa level structure. This is also essentially by construction: $g_p$ was defined to make the left-hand diagram commute, and the Igusa level structure $g_p^*j$ at $x\cdot(g^p,g_p)$ is defined to make the right-hand diagram commute.
\[
\xymatrix{
\Sigma \ar[d]_{g_p^{-1}} \ar[r]^{j} & \cA_x[p^\infty] \ar[d] & & \Sigma \ar[r]^{j} \ar[d]_{p^mg_p^{-1}} & \cA_x[p^\infty] \ar[d] \\
\Sigma \ar[r]_{j_1} & \cA_{x_1}[p^\infty] & &  \Sigma \ar[r]_{g_p^*j\hspace{4.5em}} & \cA_x[p^\infty]/j(\ker p^m g_p^{-1})
}
\]
But $\cA_{x_1}[p^\infty]=\cA_x[p^\infty]/j(\ker p^mg_p^{-1})$, because $x_1$ and $x\cdot(g^p,g_p)$ have the same associated abelian variety. Thus (after adjusting the vertical arrows by $p^m$ in the first diagram to make them isogenies), we see that $j_1$ and $g_p^*j$ both make the same diagram commute---and since there is a unique isomorphism making the diagram commute, we conclude $j_1=g_p^*j$ as desired.
\end{proof}

Recall that in \ref{hodge central leaf} we have fixed a class $\bfb\in B(G)$, which specifies the isogeny class of our fixed $p$-divisible group $\Sigma$, and therefore the Newton stratum over which our Igusa variety $\Ig_\Sigma$ lies.

\begin{lemma}
\label{intersection}
Each isogeny class in $\osS_{K_p}(G,X)(\oFFp)$ is contained in a single Newton stratum. The isogeny classes in $\osS_{K_p}(G,X)(\oFFp)$ which give rise to a non-empty isogeny class in $\Ig_\Sigma(\oFFp)$ are precisely those contained in the $\bfb$-stratum $\osS_{K_p}^{(\bfb)}(G,X)(\oFFp)$.
\end{lemma}
\begin{proof}
By the moduli description of isogeny classes on the Shimura variety in Proposition \ref{isogeny moduli description}, if two points $x_1,x_2\in\osS_{K_p}(G,X)(\oFFp)$ lie in the same isogeny class, there is a quasi-isogeny $\cA_{x_1}\to\cA_{x_2}$ such that the induced map $\DD(\cA_{x_2}[p^\infty])\to\DD(\cA_{x_1}[p^\infty])$ takes $s_{\alpha,0,x_2}$ to $s_{\alpha,0,x_1}$. Thus it induces an isomorphism between the isocrystals with $G$-structure at $x_1$ and $x_2$, which shows that $x_1$ and $x_2$ lie in the same Newton stratum. This verifies the first claim.

Since each isogeny class on the Igusa variety is the preimage of an isogeny class on the Shimura variety (Proposition \ref{pullback}), the isogeny classes in the Shuimura variety $\osS_{K_p}(G,X)(\oFFp)$ which give rise to a non-empty isogeny class in $\Ig_\Sigma(\oFFp)$ are precisely those which intersect the central leaf $C_\Sigma$. Thus to prove the claim, it suffices to show that every isogeny class in $\osS_{K_p}^{(\bfb)}(G,X)(\oFFp)$ intersects the central leaf.

Let $x_1\in\osS_{K_p}^{(\bfb)}(G,X)(\oFFp)$ and $x_2\in C_\Sigma(\oFFp)$. The point $x_2$ is auxiliary; we use it to show $x_1$ is isogenous to a point in $C_\Sigma$, but this point may not be $x_2$.

Choosing $\tilde{x}_1\in \Sh_{K_p}(G,X)(L)$ specializing to $x_1$, we have isomorphisms
\[
\xymatrix{
\DD(\cA_{x_1}[p^\infty])\ar@{<->}[r]
& H^1_\crys(\cA_{x_1}/\OO_L) \ar@{<->}[d] & \\
& H^1_\et(\cA_{\tilde{x}_1,\overline{L}},\ZZ_p)\otimes_{\ZZ_p}\OO_L \ar@{<->}[r]
& V_{\ZZ_p}^*\otimes_{\ZZ_p}\OO_L,
}
\]
where the vertical arrow is given by \cite[1.3.7(2)]{MR3630089}. We use this isomorphism (and the same for $x_2$) to identify $\DD(\cA_{x_1}[p^\infty])$ and $\DD(\cA_{x_2}[p^\infty])$ with $V_{\ZZ_p}^*\otimes_{\ZZ_p}\OO_L$, so that we can name isogenies and Frobenius by elements of $G(L)$.

Now, since $x_1$ and $x_2$ lie in the same Newton stratum, we can choose an isogeny $\cA_{x_1}[p^\infty]\to\cA_{x_2}[p^\infty]$. Using the above identifications, this isogeny is given by an element $g^{-1}\in G(L)$ (i.e. we define $g$ to be the inverse of this isogeny), and
\[
\DD(\cA_{x_2}[p^\infty])=g\cdot\DD(\cA_{x_1}[p^\infty])\subset\DD(\cA_{x_1}[p^\infty])\otimes_{\ZZ_p}\QQ_p.
\]
This relates the Frobenius $b_1\sigma$ on $\DD(\cA_{x_1}[p^\infty])$ and $b_2\sigma$ on $\DD(\cA_{x_2}[p^\infty])$, regarding $b_1,b_2\in G(L)$, by
\[
g^{-1}b_2\sigma(g)=b_1.
\]
By \cite[1.1.12]{MR3630089}, we have $b_1\in G(\OO_L)v_0(p)G(\OO_L)$, where $v_0=\sigma(\mu_0^{-1})$ and $\mu_0^{-1}$ is a $G_{\OO_L}$-valued cocharacter giving the filtration on $\DD(\cA_{x_1}[p^\infty])\otimes_{\OO_L}\oFFp$. By \cite[1.3.7(3)]{MR3630089}, $\mu_0$ is conjugate to $\mu_h$ for $h\in X$, and therefore conjugate to the cocharacter $\mu$ specified in \ref{ADLV}. This implies that $v_0$ is conjugate to the cocharacter $v$ of \ref{ADLV}, and since both are $G_{\OO_L}$-valued, they are conjugate by an element of $G(\OO_L)$. Thus $G(\OO_L)v_0(p)G(\OO_L)=G(\OO_L)v(p)G(\OO_L)$. Furthermore we can choose our basepoint $x$ of \ref{ADLV} to be $x_2$, so that the element $b\in G(L)$ arising from this basepoint is $b_2$.

Combining these facts, we get
\[
g^{-1}b\sigma(g)=b_1\in G(\OO_L)v_0(p)G(\OO_L)=G(\OO_L)v(p)G(\OO_L),
\]
which shows that $g$ defines an element of $X_v(b)$. Furthermore, the image of $g$ under the isogeny class map associated to $x_1$ is a point with Dieudonn\'e module $g\cdot\DD(\cA_{x_1}[p^\infty])\cong\DD(\cA_{x_2}[p^\infty])$. This image is a point in the central leaf which is isogenous to $x_1$, as desired.
\end{proof}

Finally we summarize the results of this section in the following proposition.

\begin{proposition}
\label{isogeny class comparison}
There is a canonical bijection between isogeny classes on $\Ig_\Sigma(\oFFp)$ and isogeny classes on $\osS_{K_p}(G,X)(\oFFp)$ contained in the $\bfb$-stratum, given by taking preimage under the map $\Ig_\Sigma(\oFFp)\to\osS_{K_p}(G,X)(\oFFp)$.
\end{proposition}
\begin{proof}
This follows from Lemma \ref{pullback} and Lemma \ref{intersection}.
\end{proof}

\subsection{$\bfb$-admissible Morphisms of Galois Gerbs}
\label{b-admissible morphisms of galois gerbs}

Recall from \ref{hodge central leaf} that we have fixed a class $\bfb\in B(G)$ specifying the Newton stratum our Igusa variety lies over.

\begin{definition}
\label{b-admissible morphism}
A morphism $\phi:\fQ\to\fG_G$ is \emph{$\bfb$-admissible} if it is admissible as in Definition \ref{admissible morphism}, with a refined condition at $p$: that $\phi(p)\circ\zeta_p:\fG_p\to\fG_G(p)$ is conjugate to an unramified morphism $\theta$ with $[b_\theta]=\bfb$.
\end{definition}

This property is preserved under conjugation by $G(\overline{\QQ})$, so we have a well-defined notion of $\bfb$-admissibility for a conjugacy class $[\phi]$ of admissible morphisms.

\paragraph{}
\label{AA_f embedding}
If $\phi:\fQ\to\fG_G$ is $\bfb$-admissible, we define a morphism $I_\phi(\AA_f)\to G(\AA_f^p)\times J_b(\QQ_p)$ as follows. Recall from \S\ref{section: galois gerbs} the set
\[
X^p(\phi)=\{g=(g_\ell)\in G(\overline{\AA}_f^p):\Int(g_\ell)\circ\xi_\ell=\phi(\ell)\circ\zeta_\ell\}
\]
is a $G(\AA_f^p)$-torsor, and therefore the natural left-multiplication action of $I_\phi(\AA_f^p)$ on $X^p(\phi)$ gives a morphism $I_\phi(\AA_f^p)\to G(\AA_f^p)$ well-defined up to $G(\AA_f^p)$-conjugacy. At $p$, let $\theta:\fG_p\to\fG_G(p)$ be an unramified morphism conjugate to $\phi(p)\circ\zeta_p$. Then we have maps
\[
I_\phi\hookrightarrow I_{\phi(p)\circ\zeta_p}\overset{\sim}{\longrightarrow}I_\theta\overset{\sim}{\longrightarrow}J_{b_\theta}\overset{\sim}{\longrightarrow}J_b,
\]
the first isomorphism because $\theta$ is conjugate to $\phi(p)\circ\zeta_p$, the second from Lemma \ref{unramified morphism facts}, and the third because $b_\theta$ is $\sigma$-conjugate to $b$ (by the $\bfb$-admissible condition). This gives a map $I_\phi(\QQ_p)\to J_b(\QQ_p)$ well-defined up to $J_b(\QQ_p)$-conjugacy.

We have produced a morphism
\[
I_\phi(\AA_f)\to G(\AA_f^p)\times J_b(\QQ_p).
\]
Define
\begin{equation}
\label{S^Ig(phi)}
S^{\Ig}(\phi)=I_\phi(\QQ)\backslash G(\AA_f^p)\times J_b(\QQ_p),
\end{equation}
where the action of $I_\phi(\QQ)$ is given by the map above. This is the set that we will use to parametrize the isogeny class $\sI^\Ig$ corresponding to a conjugacy class $[\phi]$---although we will only identify the action of $I_\phi(\QQ)$ with $I_x(\QQ)$ up to a twist, as described in \S\ref{tau-twists section}.

\subsection{Kottwitz Triples}
\label{kottwitz triples}

To make the connection between isogeny classes and admissible morphisms, we import the technique of \cite{MR3630089}.

\paragraph{}
A \emph{Kottwitz triple of level $r$} is a triple $\fk=(\gamma_0,\gamma,\delta)$ consisting of
\begin{itemize}
\item $\gamma_0\in G(\QQ)$ a semi-simple element which is elliptic in $G(\RR)$,
\item $\gamma=(\gamma_\ell)_{\ell\neq p}\in G(\AA_f^p)$ conjugate to $\gamma_0$ in $G(\overline{\AA}_f^p)$, and
\item $\delta\in G(\QQ_{p^r})$ such that $\gamma_0$ is conjugate to $\gamma_p=\delta\sigma(\delta)\cdots\sigma^{r-1}(\delta)$ in $G(\overline{\QQ}_p)$;
\end{itemize}
this data is required to satisfy the further condition (to be explained presently) that
\begin{itemize}
\item[$(*)$] there is an inner twist $I$ of $I_0$ over $\QQ$ with $I\otimes_{\QQ}\RR$ anisotropic mod center, and $I\otimes_{\QQ}\QQ_v$ is isomorphic to $I_v$ as inner twists of $I_0$ for all finite places $v$ of $\QQ$.
\end{itemize}
Here $I_0$ is the centralizer of $\gamma_0^n$ in $G$, and $I_\ell$ for $\ell\neq p$ is the centralizer of $\gamma_\ell^n$ in $G_{\QQ_\ell}$, and $I_p$ is a $\QQ_p$-group defined on points by
\[
I_p=\{g\in G(W(\FF_{p^n})\otimes_{\ZZ_p}R):g^{-1}\delta\sigma(g)=\delta\};
\]
all of these groups stabilize for $n$ sufficiently large, and we mean to take the stabilized group.

A \emph{Kottwitz triple} is an equivalence class of Kottwitz triples of various level, where we take the smallest equivalence relation so that
\begin{itemize}
\item two triples $(\gamma_0,\gamma,\delta)$ and $(\gamma_0',\gamma',\delta')$ of the same level $r$ are equivalent if $\gamma_0,\gamma_0'$ are conjugate in $G(\overline{\QQ})$ and $\gamma,\gamma'$ are conjugate in $G(\AA_f^p)$ and $\delta, \delta'$ are $\sigma$-conjugate in $\QQ_{p^r}$, and
\item a triple $(\gamma_0,\gamma,\delta)$ of level $r$ is equivalent to the triple $(\gamma_0^m,\gamma^m,\delta)$ of level $rm$.
\end{itemize}

Define a Kottwitz triple $(\gamma_0,\gamma,\delta)$ to be \emph{$\bfb$-admissible} if the $\sigma$-conjugacy class of $\delta$ is $\bfb$. This is clearly seen to be preserved under the equivalences above.

A \emph{refined Kottwitz triple} is a tuple $\tilde{\fk}=(\gamma_0,\gamma,\delta,I,\iota)$ where $\fk=(\gamma_0,\gamma,\delta)$ is a Kottwitz triple (we say $\tilde{\fk}$ is a refinement of $\fk$), and $I$ is a group as in condition $(*)$ above, and $\iota:I\otimes_\QQ\AA_f\to I_{\AA_f^p}\times I_p$ is an isomorphism of inner twists of $I_0$, where $I_{\AA_f^p}$ is the centralizer of $\gamma^n$ in $G_{\AA_f^p}$ for $n$ sufficiently large.

We consider two refined Kottwitz triples $(\gamma_0,\gamma,\delta,I,\iota)$ and $(\gamma_0',\gamma',\delta',I',\iota')$ to be equivalent if
\begin{itemize}
\item $(\gamma_0,\gamma,\delta)$ is equivalent to $(\gamma_0',\gamma',\delta')$ as Kottwitz triples, so that we can (and do) identify the groups $I_{\AA_f^p}\times I_p$ in each case; and
\item there is an isomorphism $I\to I'$ over $\QQ$ as inner twists of $I_0$ which intertwines $\iota,\iota'$.
\end{itemize}

The condition $(*)$ determines the inner twist $I$ uniquely up to conjugation by $I(\QQ)$, so the last condition above is equivalent to requiring that $\iota,\iota'$ are intertwined up to conjugation by $I(\QQ)$.

If $\tilde{\fk}$ is a refined Kottwitz triple, then the isomorphism $\iota:I\otimes_{\QQ}\AA_f\to I_{\AA_f^p}\times I_p$ gives an injection
\[
I(\AA_f)\hookrightarrow G(\AA_f^p)\times J_\delta(\QQ_p),
\]
as $I_{\AA_f^p}$ is a subgroup of $G_{\AA_f^p}$ and $I_p$ is a subgroup of $J_\delta$. For $\tilde{\fk}$ a $\bfb$-admissible Kottwitz triple, we define
\begin{equation}
\label{S^Ig(KT)}
S^\Ig(\tilde{\fk})=I(\QQ)\backslash G(\AA_f^p)\times J_b(\QQ_p),
\end{equation}
where we have used $\bfb$-admissibility to replace $J_\delta$ with the isomorphic group $J_b$. This is the intermediate set that allows us to connect an isogeny class $\sI^\Ig$ with our parametrizing set $S^\Ig(\phi)$.

\paragraph{}
\label{associated KT}
Let $x\in\Ig_\Sigma(\oFFp)$. We recall how to associate a refined Kottwitz triple to the isogeny class $\sI^\Ig_x$, following \cite[4.4.6]{MR3630089}. The same construction associates a refined Kottwitz triple to an isogeny class $\sI^\Sh$ on the Shimura variety, and since it is the same construction the triples match for isogeny classes $\sI^\Ig$ and $\sI^\Sh$ matched by Proposition \ref{isogeny class comparison}.

The level structure $\eta^p$ at $x$ identifies the group $G_{\QQ_\ell}$ with the subgroup of $\GL(H^1_\et(\cA_x,\QQ_\ell))$ fixing the tensors $\{s_{\alpha,\ell,x}\}\subset H^1_\et(\cA_x,\QQ_\ell)^\otimes$, and this allows us to write the geometric Frobenius on $H^1_\et(\cA_x,\QQ_\ell)$ as an element $\gamma_\ell\in G(\QQ_\ell)$. Let $\gamma=(\gamma_\ell)\in G(\AA_f^p)$.

At $p$, we similarly have an isomorphism
\[
V^\vee_{\ZZ_{p^r}}\overset{\sim}{\longrightarrow}\DD(\cA_x[p^\infty])(\ZZ_{p^r})
\]
which identifies $G_{\QQ_p}$ with the subgroup of $\GL(\DD(\cA_x[p^\infty])(\ZZ_{p^r}))$ fixing the tensors $\{s_{\alpha,0,x}\}$, and allows us to write the Frobenius on $\DD(\cA_x[p^\infty])$ as $\delta\sigma$ for some $\delta\in G(\QQ_{p^r})$.

With this choice of $\gamma$ and $\delta$, \cite[Cor 2.3.1]{MR3630089} states that there is an element $\gamma_0\in G(\QQ)$ that makes $(\gamma_0,\gamma,\delta)$ a Kottwitz triple, which we denote $\fk(x)$.

Furthermore if we take $I=I_x\subset\Aut_{\QQ}(\cA_x)$ then the identifications above give an isomorphism $\iota:I_x\otimes_\QQ\AA_f\to I_{\AA_f^p}\times I_p$ by taking $\ell$-adic and crystalline cohomology
\[
\Aut_\QQ(\cA_x)\to\GL(H^1_\et(\cA_x,\QQ_\ell))\times\GL(H^1_\crys(\cA_x/\ZZ_{p^r})).
\]
This completes the refined Kottwitz triple
\[
\tilde{\fk}(x)=(\gamma_0,\gamma,\delta,I,\iota)
\]
associated to the isogeny class $\sI^\Ig_x$.

Let $\phi:\fQ\to\fG_G$ be an admissible morphism. We refer to \cite[4.5.1]{MR3630089} for the construction of a Kottwitz triple $\fk(\phi)$ associated to the conjugacy class $[\phi]$, except to note that the element $\delta$ appearing in this Kottwitz triple is a $\sigma$-conjugate of the element $b_\theta$ produced by an unramified morphism $\theta$ conjugate to $\phi(p)\circ\zeta_p$ as in Definition \ref{admissible morphism}. This Kottwitz triple has a natural refinement $\tilde{\fk}(\phi)$ taking $I=I_\phi$.

We now establish a number of simple compatibilities between isogeny classes, $\bfb$-admissible morphisms, and their associated Kottwitz triples.

\begin{lemma}
\label{b-admissible iff b-stratum}
Let $\sI^\Sh$ be an isogeny class in $\osS_{K_p}(G,X)(\oFFp)$. The associated Kottwitz triple $\fk(\sI^\Sh)$ is $\bfb$-admissible if and only if $\sI^\Sh$ is contained in the $\bfb$-stratum.
\end{lemma}
\begin{proof}
The Kottwitz triple $\fk(\sI^\Sh)=(\gamma_0,\gamma,\delta)$ is $\bfb$-admissible if the $\sigma$-conjugacy class of $\delta$ is $\bfb$---but $\delta$ arises from the Frobenius on the Dieudonn\'e module at some point in $\sI^\Sh$, and so the $\sigma$-conjugacy class of $\delta$ records the Newton stratum in which the isogeny class lies.
\end{proof}

\begin{lemma}
\label{SIg(IC)=SIg(KT)}
For any $x\in\Ig_\Sigma(\oFFp)$, there is a $G(\AA_f^p)\times J_b(\QQ_p)$-equivariant bijection
\[
\sI^\Ig_x\cong S^{\Ig}(\tilde{\fk}(x)).
\]
\end{lemma}
\begin{proof}
Compare Definition \ref{igusa isogeny classes} with \eqref{S^Ig(KT)}, noting that the refined Kottwitz triple $\tilde{\fk}(x)$ is $\bfb$-admissible and has $I=I_x$.
\end{proof}

\begin{lemma}
\label{b-admissible phi iff b-admissible KT}
An admissible morphism $\phi:\fQ\to\fG_G$ is $\bfb$-admissible if and only if the associated Kottwitz triple $\fk(\phi)=(\gamma_0,\gamma,\delta)$ is $\bfb$-admissible.
\end{lemma}
\begin{proof}
The $\bfb$-admissibility of an admissible morphism $\phi$ depends on the $\sigma$-conjugacy class $[b_\theta]$ produced by an unramified morphism $\theta$ conjugate to $\phi(p)\circ\zeta_p$, while for a Kottwitz triple $(\gamma_0,\gamma,\delta)$ it depends on the $\sigma$-conjugacy class of $\delta$. But $\delta$ is a $\sigma$-conjugate of $b_\theta$, so these conditions are equivalent.
\end{proof}

\begin{lemma}
\label{SIg(AM)=SIg(KT)}
For any $\bfb$-admissible morphism $\phi:\fQ\to\fG_G$ and $\tau\in I_\phi^\ad(\AA_f^p)$, there is a $G(\AA_f^p)\times J_b(\QQ_p)$-equivariant bijection
\[
S^{\Ig}(\phi)
\cong S^{\Ig}(\tilde{\fk}(\phi)).
\]
\end{lemma}
\begin{proof}
Compare \eqref{S^Ig(phi,tau)} with \eqref{S^Ig(KT)}, and note that $\tilde{\fk}(\phi,\tau)$ has $I=I_\phi$.
\end{proof}

\subsection{$\tau$-twists}
\label{tau-twists section}

Our goal is to parametrize each isogeny class $\sI^\Ig\subset\Ig_\Sigma(\oFFp)$ by a set $S^\Ig(\phi)=I_\phi(\QQ)\backslash G(\AA_f^p)\times J_b(\QQ_p)$. However, we will only identify the action of $I_x(\QQ)$ with $I_\phi(\QQ)$ up to twisting by an element $\tau\in I_\phi^\ad(\AA_f)$. We develop the necessary theory in this section, following \cite[\S2.6]{KSZ}.

\paragraph{}
Let $\phi$ be an admissible morphism, and define
\begin{gather*}
\cH(\phi)=I_\phi(\AA_f)\backslash I_\phi^\ad(\AA_f)/I_\phi^\ad(\QQ), \\
\fE^p(\phi)=I_\phi(\AA_f^p)\backslash I_\phi^\ad(\AA_f^p).
\end{gather*}
These sets can be given the structure of abelian groups (by comparison with certain abelianized cohomology groups, see \cite[2.6.13, Lemma 2.6.14]{KSZ}). By weak approximation the natural inclusion $I_\phi^\ad(\AA_f^p)\to I_\phi^\ad(\AA_f)$ induces a surjection $\fE^p(\phi)\to\cH(\phi)$, which allows us to lift an element of $\cH(\phi)$ to $I_\phi^\ad(\AA_f^p)$, or further $I_\phi(\overline{\AA}_f^p)$. We will often do this implicitly, writing $\tau\in\Gamma(\cH)$ and $\tau(\phi)\in I_\phi^\ad(\AA_f^p)$, when the ambiguity in the lift is harmless.

Let $\AM$ be the set of admissible morphisms. Since we want to consider assignments of an element of $\cH(\phi)$ for all $\phi$ simultaneously, it is convenient to consider $\cH(\phi)$ as the stalks of a sheaf $\cH$ on $\AM$ (regarded as a discrete topological space), and similarly $\fE^p(\phi)$ the stalks of a sheaf $\fE^p$. Let $\Gamma(\cH)$ and $\Gamma(\fE^p)$ be the global sections of these sheaves, so an element $\tau\in\Gamma(\cH)$ assigns to each admissible morphism $\phi$ an element $\tau(\phi)\in\cH(\phi)$, and similarly for $\fE^p$. The surjections $\fE^p(\phi)\to\cH(\phi)$ produce surjections $\fE^p\to\cH$ and $\Gamma(\fE^p)\to\Gamma(\cH)$.

Define an equivalence relation $\phi_1\approx\phi_2$ if $\phi_1^\Delta$ is conjugate to $\phi_2^\Delta$ by $G(\overline{\QQ})$. If $\phi_1\approx\phi_2$, then there are canonical isomorphisms
\begin{gather*}
\Comp_{\phi_1,\phi_2}:\cH(\phi_1)\to\cH(\phi_2) \\
\Comp_{\phi_1,\phi_2}^{\fE^p}:\fE^p(\phi_1)\to\fE^p(\phi_2)
\end{gather*}
satisfying the relations $\Comp_{\phi_2,\phi_3}\circ\Comp_{\phi_1,\phi_2}=\Comp_{\phi_1,\phi_3}$ and $\Comp_{\phi_1,\phi_1}=\id_{\cH(\phi_1)}$, and similarly for $\fE^p$. These isomorphisms show that $\cH$ and $\fE^p$ are pulled back from sheaves $\cH/\approx$ and $\fE^p/\approx$ on $\AM/\approx$, under the natural quotients
\[
\AM\to\AM/\text{conj}\to\AM/\approx.
\]
Write $\cH/\text{conj}$ and $\fE^p/\text{conj}$ for the intermediate pullbacks to $\AM/\text{conj}$, the set of admissible morphisms up to conjugacy.

\paragraph{}
\label{par: tori-rational and sha-compatible}
Let $\Gamma(\cH)_0$ be the set of global sections of $\cH$ that descend to $\AM/\approx$, and $\Gamma(\cH)_1$ those that descend to $\AM/\text{conj}$, so we have $\Gamma(\cH)_0\subset\Gamma(\cH)_1\subset\Gamma(\cH)$. Define $\Gamma(\fE^p)_0\subset\Gamma(\fE^p)_1\subset\Gamma(\fE^p)$ similarly. The surjection $\Gamma(\fE^p)\to\Gamma(\cH)$ induces a surjection $\Gamma(\fE^p)_0\to\Gamma(\cH)_0$.

There is one further technical definition we will need, namely the notion of \emph{tori-rationality} of an element of $\Gamma(\cH)$ or $\Gamma(\fE^p)$. For this we refer to \cite[Def 2.6.19]{KSZ}. We will also need the fact \cite[Lemma 2.6.20]{KSZ} that an element of $\Gamma(\cH)$ is tori-rational exactly when one (equivalently, every) lift to $\Gamma(\fE^p)$ is tori-rational.

Let $\tau\in \Gamma(\cH)_1$. Define
\begin{equation}
\label{S^Ig(phi,tau)}
S^{\Ig}_\tau(\phi)=I_\phi(\QQ)\backslash G(\AA_f^p)\times J_b(\QQ_p),
\end{equation}
as in \eqref{S^Ig(phi)}, except the action of $I_\phi(\QQ)$ composed with the action of $\tau(\phi)$. Likewise define $\tilde{\fk}(\phi,\tau)$ to be the refined Kottwitz triple obtained from $\tilde{\fk}(\phi)$ replacing $\iota:I_\phi\otimes_\QQ\AA_f\overset{\sim}{\to}I_{\AA_f^p}\times I_p$ by its composition with the action of $\tau$. By taking $\tau\in\Gamma(\cH)_1$ we ensure that these definitions only depend on the conjugacy class of $\phi$. Then we have an immediate analogue of Lemma \ref{SIg(AM)=SIg(KT)}.

\begin{lemma}
\label{SIg(AM,tau)=SIg(KT,sigma)}
Let $\tau\in\Gamma(\cH)$. For any $\bfb$-admissible morphism $\phi:\fQ\to\fG_G$, there is a $G(\AA_f^p)\times J_b(\QQ_p)$-equivariant bijection
\[
S^{\Ig}_\tau(\phi)
\cong S^{\Ig}(\tilde{\fk}(\phi,\tau)).
\]
\end{lemma}

\subsection{Langlands--Rapoport-$\tau$ Conjecture for Igusa Varieties of Hodge Type}
\label{statement and proof}

\begin{proposition}
\label{IC-AM bijection}
There exists a bijection
\[
\left\{\parbox{8em}{isogeny classes in $\Ig_\Sigma(\oFFp)$}\right\}
\longleftrightarrow
\left\{\parbox{11em}{conjugacy classes of $\bfb$-admissible morphisms $\phi$}\right\},
\]
compatible with (unrefined) Kottwitz triples. Furthermore, there is a tori-rational element $\tau\in\Gamma(\cH)_0$ making such a bijection compatible with refined Kottwitz triples, in the sense that if $\sI_x$  maps to $[\phi]$, then $\tilde{\fk}(x)$ is equivalent to $\tilde{\fk}(\phi,\tau(\phi))$.
\end{proposition}
\begin{proof}
In \cite{MR3630089}, Kisin constructs a bijection
\[
\left\{\parbox{8em}{isogeny classes in $\osS_{K_p}(G,X)(\oFFp)$}\right\}
\longleftrightarrow
\left\{\parbox{6em}{conj. classes of admissible morphisms $\phi$}\right\},
\]
compatible with Kottwitz triples (that is, corresponding elements give rise to equivalent Kottwitz triples). This bijection induces a bijection in our case, as we now explain.

By Proposition \ref{isogeny class comparison} we can consider the set of isogeny classes in $\Ig_\Sigma(\oFFp)$ as a subset of the set of isogeny classes in $\osS_{K_p}(G,X)(\oFFp)$, namely the subset of isogeny classes contained in the $\bfb$-stratum. By Lemma \ref{b-admissible iff b-stratum}, this subset is characterized as those isogeny classes whose associated Kottwitz triple is $\bfb$-admissible.

On the other side, the $\bfb$-admissible condition on admissible morphisms $\fQ\to\fG_G$ is a strengthening of the admissible condition, so we can regard the set of conjugacy classes of $\bfb$-admissible morphisms as a subset of the set of conjugacy classes of admissible morphisms. By Lemma \ref{b-admissible phi iff b-admissible KT}, this subset is characterized as those conjugacy classes whose associated Kottwitz triple is $\bfb$-admissible.

Because the bijection above is compatible with Kottwitz triples, it induces a bijection between the subsets on each side corresponding to $\bfb$-admissible Kottwitz triples, and this gives our desired bijection.

It remains to settle the last claim regarding the $\tau$-twist. The bijection between Shimura isogeny classes and conjugacy classes of admissible morphisms constructed in \cite{MR3630089} is not known to preserve refined Kottwitz triples, and so we can only identify the refined Kottwitz triples on each side up to a $\tau$-twist. That is, if $\sI_x$ maps to $[\phi]$, then $\tilde{\fk}(x)$ is equivalent to $\tilde{\fk}(\phi,\tau(\phi))$ for some $\tau\in\Gamma(\cH)$.

In \cite{KSZ} it is conjectured that it should be possible to take $\tau$ to be tori-rational and to lie in $\Gamma(\cH)_0$ (part of their Conjecture 3.4.4), and one of their main results \cite[Thm 6.3.6]{KSZ} is to prove this conjecture in the case of abelian type. By the same argument as above, their result carries over to our case, which completes the proof.
\end{proof}

\begin{theorem}
\label{igusa LR}
There exists a tori-rational element $\tau\in \Gamma(\cH)_0$ admitting a $G(\AA_f^p)\times J_b(\QQ_p)$-equivariant bijection
\[
\Ig_\Sigma(\oFFp)\overset{\sim}{\longrightarrow}\coprod_{[\phi]} S^{\Ig}_\tau(\phi),
\]
where the disjoint union ranges over conjugacy classes of $\bfb$-admissible morphisms $\phi:\fQ\to\fG_G$.
\end{theorem}
\begin{proof}
By Proposition \ref{IC-AM bijection}, the conjugacy classes $[\phi]$ parametrizing the disjoint union are in bijection with the set of isogeny classes in $\Ig_\Sigma(\oFFp)$. Furthermore there is a tori-rational element $\tau\in\Gamma(\cH)_0$ such that for corresponding $\sI_x$ and $[\phi]$, the refined Kottwitz triples $\tilde{\fk}(x)$ and $\tilde{\fk}(\phi,\tau(\phi))$ are equivalent.

Combining the equivalence of these refined Kottwitz triples with Lemma \ref{SIg(IC)=SIg(KT)} and Lemma \ref{SIg(AM,tau)=SIg(KT,sigma)}, we have $G(\AA_f^p)\times J_b(\QQ_p)$-equivariant bijections
\[
\sI_x\cong S^\Ig(\tilde{\fk}(x))\cong S^\Ig(\tilde{\fk}(\phi,\tau))\cong S^\Ig_\tau(\phi).
\]
Thus we have produced bijections between the set of isogeny classes and the set of conjugacy classes of $\bfb$-admissible morphisms, and bijections between each isogeny class and the set $S^\Ig_\tau(\phi)$ produced by its corresponding admissible morphism. Together, these produce a bijection $\Ig_\Sigma(\oFFp)\overset{\sim}{\to}\coprod_{[\phi]}S^\Ig_\tau(\phi)$ as desired.
\end{proof}

\section{Point-Counting Formula for Igusa Varieties of Hodge Type}
\label{sec: point-counting}

We continue to use the notation of the previous sections. In particular, we have fixed a Shimura datum $(G,X)$ of Hodge type, a class $\bfb\in B(G,\mu^{-1})$, and representative $b\in G(L)$ of this class, which data defines an Igusa variety $\Ig_\Sigma$.

\subsection{Acceptable Functions and Fujiwara's Trace Formula}
\label{sec: acceptable functions and fujiwara}

\paragraph{}
\label{trace preparations}
The (right) action of $G(\AA_f^p)\times S_b$ on $\Ig_\Sigma$ extends to an (left) action of $G(\AA_f^p)\times J_b(\QQ_p)$ on $H_c^i(\Ig_\Sigma,\sL_\xi)$ as in \eqref{H_c(Ig)}. Let $f\in C^\infty_c(G(\AA_f^p)\times J_b(\QQ_p))$, a smooth (i.e. locally constant) compactly supported function on $G(\AA_f^p)\times J_b(\QQ_p)$. We define an operator on $H_c^i(\Ig_\Sigma,\sL_\xi)$, also called $f$ by abuse of notation, by
\[
v\mapsto\int_{G(\AA_f^p)\times J_b(\QQ_p)}f(x)x\cdot v\;dx.
\]
This operator has finite rank, so we can write
\[
\tr(f\mid H_c(\Ig_\Sigma,\sL_\xi))=\sum_i(-1)^i\tr(f\mid H_c^i(\Ig_\Sigma,\sL_\xi)).
\]

Any $f\in C^\infty_c(G(\AA_f^p)\times J_b(\QQ_p))$ is a finite linear combination of indicator functions $\one_{UgU}$ for $U\subset G(\AA_f^p)\times J_b(\QQ_p)$ compact open and $g\in G(\AA_f^p)\times J_b(\QQ_p)$, so by linearity of trace it suffices to consider $f=\one_{UgU}$. Writing $UgU=\coprod_i g_igU$ for some finite collection $g_i\in U$, the action of $\one_{UgU}$ on $H_c^i(\Ig_\Sigma,\sL_\xi)$ is given by
\begin{align*}
\int_{G(\AA_f^p)\times J_b(\QQ_p)}\one_{UgU}(x)x\cdot v\;dx
& = \sum_i g_ig\int_U x\cdot v\;dx.
\end{align*}
Now $\int_U x\cdot v\;dx$ is $\vol(U)$ times projection to $H_c^i(\Ig_\Sigma,\sL_\xi)^U$, which realizes the action of $\one_{UgU}$ as $\vol(U)$ times the operator
\[
v\mapsto \sum_i g_ig\cdot v\qquad\text{ on }\qquad H_c^i(\Ig_\Sigma,\sL_\xi)^U.
\]
We call this operator $v\mapsto\sum_ig_ig\cdot v$ the \emph{double coset operator $[UgU]$}.

Define
\[
U_p(m)=\ker(\Aut(\Sigma,\lambda_\Sigma,\{s_{\alpha,\Sigma}\})\to\Aut(\Sigma[p^m],\lambda_\sigma,\{s_{\alpha,\Sigma}\}))\subset J_b(\QQ_p).
\]
These subsets form a neighborhood basis of the identity in $J_b(\QQ_p)$, so we can assume $U=U^p\times U_p(m)$ for $U^p\subset G(\AA_f^p)$ compact open. Then
\[
H_c^i(\Ig_\Sigma,\sL_\xi)^{U}=H_c^i(\Ig_{\Sigma,U^p,m},\sL_\xi).
\]
In particular, the action of $\one_{UgU}$ is $\vol(U)$ times the double coset action $[UgU]$, and we have (taking the alternating sum over $i$)
\begin{equation}
\label{indicator to double coset}
\tr(\one_{UgU}\mid H_c(\Ig_\Sigma,\sL_\xi))=\vol(U)\tr([UgU]\mid H_c(\Ig_{\Sigma,U^p,m},\sL_\xi)).
\end{equation}

\paragraph{}
Now assume $g=g^p\times g_p\in G(\AA_f^p)\times S_b$ (recall $S_b$ from \ref{S_b}), so that we can consider the action of $g$ on finite-level Igusa varieties. Then the double coset action $[UgU]$ is induced by the following correspondence.
\[
\xymatrix{
& \Ig_{\Sigma}(\oFFp)/(U\cap gUg^{-1}) \ar[dl]_{[\cdot 1]} \ar[dr]^{[\cdot g]} & \\
\Ig_{\Sigma}(\oFFp)/U & & \Ig_{\Sigma}(\oFFp)/U
}
\]
Note that this is a set-theoretic correspondence rather than an algebro-geometric correspondence, but we can obtain an algebro-geometric correspondence which induces the same action after by increasing the level of the top Igusa variety. This introduces a constant factor into the resulting formulas, but this factor is reabsorbed by defining the set of fixed points in terms of the set-theoretic correspondence as well.

\begin{definition}
\label{Fix(UgU)}
Define the fixed point set of the above correspondence by
\[
\Fix(UgU)=\{x\in \Ig_\Sigma(\oFFp)/(U\cap gUg^{-1}) : x=xg\text{ in } \Ig_\Sigma(\oFFp)/U\}.
\]
\end{definition}

Recall from \S\ref{acceptable elements section} the notion of an acceptable element of $J_b(\QQ_p)$.

\begin{definition}
\label{acceptable function}
A function $f\in C_c^\infty(G(\AA_f^p)\times J_b(\QQ_p))$ is \emph{acceptable} if
\begin{enumerate}
\item \label{acceptable elements part} for all $(g,\delta)\in \supp f$, we have $\delta\in S_b$ and $\delta$ is acceptable;
\item \label{fujiwara part} there is a sufficiently small compact open subgroup $U=U^p\times U_p(m)\subset G(\AA_f^p)\times J_b(\QQ_p)$ and a finite subset $I\subset G(\AA_f^p)\times J_b(\QQ_p)$ such that $f=\sum_{g\in I}\one_{UgU}$; and for each term in this sum, we have
\begin{enumerate}
\item $\Fix(UgU)$ is finite, and
\item the trace of the correspondence on cohomology is given by Fujiwara's formula:
\begin{equation}
\label{fujiwara equation}
\tr([UgU]\mid H_c(\Ig_{\Sigma,U^p,m},\sL_\xi))=\sum_{x\in \Fix(UgU)}\tr([UgU]\mid (\sL_\xi)_x).
\end{equation}
\end{enumerate}
\end{enumerate}
\end{definition}

The next two lemmas are slightly adjusted versions of Lemmas 6.3 and 6.4 of \cite{MR2484281}, which verify that acceptable functions are a sufficient class of test functions to determine a representation. The proofs are precisely analogous to those in \cite{MR2484281}.

\begin{lemma}[cf. {\cite[Lemma 6.3]{MR2484281}}]\label{fujiwara}
For any $f\in C^\infty_c(G(\AA_f^p)\times J_b(\QQ_p))$, the function $f^{(m,n)}$ defined by $f^{(m,n)}(x):=f(x\cdot p^m(\fr^s)^n)$ is an acceptable function for sufficiently large $m,n$.
\end{lemma}

Write $\Groth(G(\AA_f^p)\times J_b(\QQ_p))$ for the Grothendieck group of admissible representations of $G(\AA_f^p)\times J_b(\QQ_p)$.

\begin{lemma}[cf. {\cite[Lemma 6.4]{MR2484281}}]\label{acceptable functions suffice}
If $\Pi_1,\Pi_2\in\Groth(G(\AA_f^p)\times J_b(\QQ_p))$ satisfy $\tr(f\mid \Pi_1)=\tr(f\mid \Pi_2)$ for all acceptable functions $f$, then $\Pi_1=\Pi_2$ as elements of the Grothendieck group.
\end{lemma}

\subsection{Preliminary Point-Counting}
\label{preliminary point-counting section}

Let $f$ be an acceptable function; our goal is a preliminary formula for $\tr(f\mid H_c(\Ig_\Sigma,\sL_\xi))$. We assume for the moment that $f$ is of the form $\one_{UgU}$ (for $U=U^p\times U_p(m)$) and satisfies condition \ref{fujiwara part} of Definition \ref{acceptable function}. Combining \eqref{indicator to double coset} and \eqref{fujiwara equation}, we find
\begin{equation}
\label{fujiwara'd trace}
\tr(\one_{UgU}\mid H_c(\Ig_\Sigma,\sL_\xi))
=\vol(U)\sum_{x\in\Fix(UgU)}\tr([UgU]\mid(\sL_\xi)_x).
\end{equation}

We proceed to analyze the set $\Fix(UgU)$ (Definition \ref{Fix(UgU)}) parametrizing the sum. 

From Theorem \ref{igusa LR} we see that for any compact open $U\subset G(\AA_f^p)\times J_b(\QQ_p)$ we have a bijection
\begin{equation}
\label{LR mod U}
\Ig_\Sigma(\oFFp)/U\overset{\sim}{\longrightarrow}\coprod_{[\phi]} I_\phi(\QQ)\backslash G(\AA_f^p)\times J_b(\QQ_p)/U,
\end{equation}
where $[\phi]$ ranges over conjugacy classes of $\bfb$-admissible morphisms, and with the action of $I_\phi(\QQ)$ possibly twisted by a tori-rational element $\tau\in\Gamma(\cH)_0$. Recall from \ref{tau-twists section} that we can lift $\tau(\phi)\in \cH(\phi)$ to an element of $I_\phi^\ad(\AA_f^p)$, which we will also call $\tau(\phi)$ by abuse of notation.

Consider $X=G(\AA_f^p)\times J_b(\QQ_p)/U$ and $Y=G(\AA_f^p)\times J_b(\QQ_p)/(U\cap gUg^{-1})$ sets with action by $I=I_\phi(\QQ)$, and maps $a,b:X\to Y$ given by $a:x\mapsto x\mod U$ and $b:x\mapsto xg\mod U$. Then the set $(I\backslash Y)^{a=b}$ on which the maps $a,b:I\backslash Y\to I\backslash X$ agree is on the one hand
\[
\{x\in I_\phi(\QQ)\backslash G(\AA_f^p)\times J_b(\QQ_p)/(U\cap gUg^{-1}) : x=xg\mod U\},
\]
the $[\phi]$-component of $\Fix(UgU)$ as in \eqref{LR mod U}; and on the other hand
\[
\coprod_\varepsilon I_{\phi,\varepsilon}(\QQ)\backslash \{x\in G(\AA_f^p)\times J_b(\QQ_p)/(U\cap gUg^{-1}) : \varepsilon x=xg \mod U\}
\]
by applying Milne's combinatorial lemma \cite[Lemma 5.3]{MR1155229} in the setup above (here we are taking $C={1}$ in the notation of \cite[Lemma 5.3]{MR1155229}, which satisfies the necessary hypotheses by Lemmas \ref{MCL stabilizer} and \ref{MCL center}).

Thus we can write
\begin{align}
\begin{split}
\label{Fix(UgU) after MCL}
&\Fix(UgU)= \\
&\coprod_{[\phi]} \coprod_\varepsilon I_{\phi,\varepsilon}(\QQ)\backslash\{x\in G(\AA_f^p)\times J_b(\QQ_p)/(U\cap gUg^{-1}) : \varepsilon x=xg\mod U\}.
\end{split}
\end{align}

Proceeding as in \cite[Lemma 7.4]{MR2484281}, we obtain the first form of the counting point formula.

\begin{proposition}
\label{preliminary formula}
Let $\tau\in\Gamma(\cH)_0$ a tori-rational element satisfying Theorem \ref{igusa LR}, which we may lift to a tori-rational element of $\Gamma(\fE^p)_0$ (still called $\tau$ by abuse of notation). For any acceptable function $f\in C_c^\infty(G(\AA_f^p)\times J_b(\QQ_p))$, we have
\[
\tr(f\mid H_c(\Ig_\Sigma,\sL_\xi))=\sum_{[\phi]}\sum_\varepsilon \frac{\vol\left(I^\circ_{\phi,\varepsilon}(\QQ)\backslash G_{b,\gamma\times\delta}^\circ\right)}{[I_{\phi,\varepsilon}(\QQ):I^\circ_{\phi,\varepsilon}(\QQ)]} O_{\gamma\times\delta}^{G(\AA_f^p)\times J_b(\QQ_p)}(f)\tr(\xi(\varepsilon)),
\]
where $\phi$ ranges over conjugacy classes of $\bfb$-admissible morphisms $\phi:\fQ\to\fG_G$ and $\varepsilon$ ranges over conjugacy classes in $I_\phi(\QQ)$, and $\gamma\times\delta$ is the image of $\Int(\tau(\phi))\varepsilon$ in $G(\AA_f^p)\times J_b(\QQ_p)$.
\end{proposition}

The remainder of the section consists of technical results deferred from earlier in the section.

\begin{lemma}\label{MCL stabilizer}
For $U$ a sufficiently small compact open subgroup of $G(\AA_f^p)\times J_b(\QQ_p)$ and $\phi$ an admissible morphism, the stabilizer in $I_\phi(\QQ)$ of any element of $G(\AA_f^p)\times J_b(\QQ_p)/U$ is $Z_G(\QQ)\cap U$. The same is true if the action of $I_\phi(\QQ)$ is twisted by an element $\tau(\phi)\in I_{\phi}^\ad(\AA_f^p)$.
\end{lemma}
\begin{proof}
The proof of \cite[Lemma 3.7.2(i)]{KSZ} carries over to our case, the only difference being the component at $p$---and in our case $\varepsilon\in I_\phi(\QQ)$ stabilizing an element $\varepsilon x_p=x_p\mod U_p$ implies that $\varepsilon$ is contained in the compact subgroup $x_pU_px_p^{-1}$ of $J_b(\QQ_p)$ (or of $G(\overline{\QQ}_p)$ via the embedding $J_b\to G$ over $\overline{\QQ}_p$), providing the necessary ingredient for the proof.
\end{proof}

\begin{lemma}\label{MCL center}
For a sufficiently small compact open subgroup $U\subset G(\AA_f^p)\times J_b(\QQ_p)$ we have $Z_G(\QQ)\cap U=\{1\}$.
\end{lemma}
\begin{proof}
Since $G$ is part of a Shimura datum of Hodge type, $Z_G^\circ$ satisfies the Serre condition---this is (equivalent to) the condition that $Z_G^\circ$ is isogenous over $\QQ$ to a torus $T^+\times T^-$ where $T^+$ is split over $\QQ$ and $T^-$ is compact over $\RR$. This implies that $Z_G^\circ(\QQ)$ is discrete in $Z_G^\circ(\AA_f)$ (e.g. \cite[Lemma 1.5.5]{KSZ}), and via 
\[
Z_G^\circ(\AA_f)\hookrightarrow G(\AA_f^p)\times J_b(\QQ_p)
\]
we see $Z_G^\circ(\QQ)$ is discrete in $G(\AA_f^p)\times J_b(\QQ_p)$ (the embedding $Z_G\hookrightarrow J_b$ coming from the fact that $J_b$ is an inner form of a Levi subgroup of $G$). Since $[Z_G(\QQ):Z_G^\circ(\QQ)]$ is finite, $Z_G(\QQ)$ is also discrete in $G(\AA_f^p)\times J_b(\QQ_p)$. Thus any sufficiently small compact open subgroup $U$ will intersect $Z_G(\QQ)$ trivially.
\end{proof}

\subsection{LR Pairs and Kottwitz Parameters}
\label{LR/KP section}

\begin{definition}
\label{def: LR pair}
An \emph{LR pair} is a pair $(\phi,\varepsilon)$ where $\phi:\fQ\to\fG_G$ is a morphism of Galois gerbs and $\varepsilon\in I_\phi(\QQ)$. The element $\varepsilon$ can also be regarded as an element of $G(\overline{\QQ})$ via $I_\phi(\QQ)\subset G(\overline{\QQ})$. Two LR pairs $(\phi_1,\varepsilon_1)$ and $(\phi_2,\varepsilon_2)$ are \emph{conjugate} if there is an element $g\in G(\overline{\QQ})$ which conjugates $\phi_1$ to $\phi_2$ and $\varepsilon_1$ to $\varepsilon_2$.

An LR pair $(\phi,\varepsilon)$ is \emph{semi-admissible} if $\phi$ is admissible, and \emph{$\bfb$-admissible} if $\phi$ is $\bfb$-admissible. We denote the set of LR pairs by $\LRP$, the subset of semi-admissible pairs by $\LRP_\sa$, and the subset of $\bfb$-admissible pairs by $\LRP_\ba$.
\end{definition}

Note that if $(\phi,\varepsilon)$ is a semi-admissible LR pair then $\varepsilon$ is semi-simple, as $I_\phi/Z_G$ is compact over $\RR$ (cf. \cite[3.1.2]{KSZ}).

\paragraph{}
\label{p-adic realization}
Let $(\phi,\varepsilon)$ be a semi-admissible LR pair. As in \ref{unramified paragraph}, the morphism $\phi(p)\circ\zeta_p:\fG_p\to\fG_G(p)$ is conjugate by some $g\in G(\overline{\QQ}_p)$ to an unramified morphism $\theta:\fG_p\to\fG_G(p)$, which defines an element $b_\theta\in G(\QQ_p^\ur)$. Since $\varepsilon$ commutes with $\phi$, its conjugate $\varepsilon'=g\varepsilon g^{-1}$ commutes with $\theta$. Then for any $\rho\in \Gal(\overline{\QQ}_p/\QQ_p^\ur)$ we have
\[
1\rtimes\rho=\theta(\rho)=\Int(\varepsilon')\circ\theta(\rho)=\Int(\varepsilon')(1\rtimes\rho)=\varepsilon'\rho(\varepsilon')^{-1}\rtimes\rho.
\]
This shows that $\varepsilon'$ is fixed by all $\rho\in\Gal(\overline{\QQ}_p/\QQ_p^\ur)$, and is therefore an element of $G(\QQ_p^\ur)$ (a priori only being an element of $G(\overline{\QQ}_p)$). Furthermore, since $\varepsilon'$ commutes with $\theta^\ur$, it must $\sigma$-centralize $b_\theta$, and we can regard it as an element of $J_{b_\theta}(\QQ_p)$.

A pair $(b_\theta,\varepsilon')$ arising from $(\phi,\varepsilon)$ in this way is said to be a \emph{$p$-adic realization} of $(\phi,\varepsilon)$. (For comparison, the set of $p$-adic realizations is called $\text{cls}_p(\phi,\varepsilon)$ in \cite[3.1.4]{KSZ}). Note that all $p$-adic realizations of an LR pair are conjugate in the following sense: if $g_1,g_2\in G(\overline{\QQ}_p)$ give rise to $p$-adic realizations $(b_1,\varepsilon_1)$ and $(b_2,\varepsilon_2)$ respectively, then
\[
b_2=(g_2g_1^{-1})b_1\sigma(g_2g_1^{-1})^{-1}
\]
and conjugation by $g_2g_1^{-1}$ gives an isomorphism $J_{b_1}\to J_{b_2}$ sending $\varepsilon_1$ to $\varepsilon_2$.

\begin{definition}
\label{def: acceptable LR pair}
Let $(\phi,\varepsilon)$ be a semi-admissible LR pair with $(b_\theta,\varepsilon')$ a $p$-adic realization. We can regard $\varepsilon'$ as an element of $J_{b_\theta}(\QQ_p)$. Define $(\phi,\varepsilon)$ to be \emph{acceptable} if $\varepsilon'$ is acceptable as an element of $J_{b_\theta}(\QQ_p)$. Note that this does not depend on the choice of $p$-adic realization, as in the paragraph just above.
\end{definition}

\paragraph{}
\label{par: special point datum}

A \emph{special point datum} is a triple $(T,h_T,i)$ where
\begin{itemize}
\item $T$ is a torus,
\item $h_T:\bbS\to T_\RR$ is a morphism from the Deligne torus, and
\item $i:T\to G$ is an embedding realizing $T$ as a maximal torus of $G$ defined over $\QQ$, and sending $h_T$ into $X$.
\end{itemize}

A special point datum induces an LR pair in the following way. Recall from \ref{quasi-motivic gerb} that the quasi-motivic Galois gerb is equipped with a distinguished morphism $\psi:\fQ\to\fG_{\Res_{\overline{\QQ}/\QQ}\GG_m}$. Composing this with the map $\fG_{\Res_{\overline{\QQ}/\QQ}\GG_m}\to\fG_T$ induced by the cocharacter $\mu_{h_T}$ and the map $\fG_T\to\fG_G$ induced by the inclusion $i$, we obtain an admissible morphism $\phi=i\circ\psi_{\mu_{h_T}}:\fQ\to\fG_G$.

Furthermore, in this setup $T(\QQ)$ as a subgroup of $G(\overline{\QQ})$ lies inside $I_\phi(\QQ)$, so any $\varepsilon\in T(\QQ)$ makes a semi-admissible LR pair $(\phi,\varepsilon)$. Such a pair is called \emph{very special}; an LR pair conjugate to a very special pair is called \emph{special}.

It is a fact \cite[3.3.9]{KSZ} that every semi-admissible pair is special.

\begin{definition}
\label{def: gg}
An LR pair $(\phi,\varepsilon)$ is \emph{gg} (abbreviated from \emph{g\"unstig gelegen}, German for ``well-positioned'') if 
\begin{itemize}
\item $\phi^\Delta$ is defined over $\QQ$;
\item $\varepsilon$ lies in $G(\QQ)$ (under the inclusion $I_{\phi}(\QQ)\subset G(\overline{\QQ})$) and is semi-simple and elliptic in $G(\RR)$; and
\item for any $\rho\in\Gal(\overline{\QQ}/\QQ)$, letting $q_\rho\in \fQ$ a lift of $\rho$ and $\phi(q_\rho)=g_\rho\rtimes\rho$, we have $g_\rho\in G_{\varepsilon}^\circ$.
\end{itemize}
\end{definition}

Very special pairs are gg, and consequently every semi-admissible LR pair is conjugate to a gg pair. We denote the set of gg LR pairs by $\LRP^\gg\subset \LRP$.

The following definition is inspired by \cite[Def 10.1]{MR2484281} and \cite[Def 1.6.2]{KSZ}.

\begin{definition}
\label{def: classical KP}
A \emph{classical Kottwitz parameter} of \emph{type $\bfb$} is a triple $(\gamma_0,\gamma,\delta)$ where
\begin{itemize}
\item $\gamma_0\in G(\QQ)$ is semi-simple and elliptic in $G(\RR)$,
\item $\gamma=(\gamma_\ell)_\ell\in G(\AA_f^p)$ such that $\gamma_\ell$ is stably conjugate to $\gamma_0$ in $G(\QQ_\ell)$, and
\item $\delta\in J_b(\QQ_p)$ is acceptable, and conjugate to $\gamma_0$ in $G(\overline{\QQ}_p)$ under the embedding $J_b(\QQ_p)\to G(\overline{\QQ}_p)$.
\end{itemize}

We say that $(\gamma_0,\gamma,\delta)$ and $(\gamma_0',\gamma',\delta')$ are \emph{equivalent} if $\gamma_0$ is stably conjugate to $\gamma_0'$ in $G(\QQ)$, $\gamma$ is conjugate to $\gamma'$ in $G(\AA_f^p)$, and $\delta$ is conjugate to $\delta'$ in $J_b(\QQ_p)$.
\end{definition}

In order to handle the case that $G_\der$ is not simply connected, we need a (closely related but) more abstract treatment.

\begin{definition}[cf. {\cite[Def 1.6.4]{KSZ}}]
\label{kottwitz parameter}
A \emph{Kottwitz parameter} is a triple $\fc=(\gamma_0,a,[b_0])$ where
\begin{enumerate}
\item $\gamma_0\in G(\QQ)$ is semi-simple and elliptic in $G(\RR)$, and we write $I_0=G_{\gamma_0}^\circ$;
\item $a$ is an element of
\[
\fD(I_0,G;\AA_f^p)=\ker\left(H^1(I_0,\AA_f^p)\to H^1(G,\AA_f^p)\right);
\]
\item $[b_0]\in B(I_0)$; and
\item \label{KP0} the image of $[b_0]$ under the Kottwitz map $B(I_0)\to B(G)\overset{\kappa}{\to}\pi_1(G)_{\Gamma_p}$ is equal to the image of $\mu^{-1}$, where $\mu$ is the cocharacter induced by an(y) element of $X$.
\end{enumerate}
\end{definition}

\paragraph{}
\label{par: isomorphism of KP}
An isomorphism of Kottwitz parameters is essentially given by conjugation by $G(\overline{\QQ})$. We make this precise as follows. Let $(\gamma_0,a,[b_0])$ be a Kottwitz parameter, $I_0=G_{\gamma_0}^\circ$, and $u\in G(\overline{\QQ})$ an element such that $\Int(u)\gamma_0=\gamma_0'$ is again in $G(\QQ)$ and $u^{-1}\rho(u)\in I_0$ for all $\rho\in \Gal(\overline{\QQ}/\QQ)$. Write $I_0'=G_{\gamma_0'}^\circ$.

To relate the away-from-$p$ parts, we consider the bijection
\begin{align*}
u_*:\fD(I_0,G;\AA_f^p) & \to\fD(I_0',G;\AA_f^p) \\
e_\rho & \mapsto ue_\rho\rho(u)^{-1}
\end{align*}
induced by the element $u$.

To relate the $p$-parts, we construct a bijection $u_*:B(I_0)\to B(I_0')$. The cocycle $\rho\mapsto u^{-1}\rho(u)\in Z^1(\QQ_p,I_0)$ is trivial in $H^1(\breve{\QQ}_p,I_0)$ by the Steinberg vanishing theorem. That is, we can find $d\in I_0(\overline{\breve{\QQ}}_p)$ so that $u^{-1}\rho(u)=d^{-1}\rho(d)$ for all $\rho\in\Gal(\overline{\breve{\QQ}}_p/\breve{\QQ}_p)$. Then we have $ud^{-1}=\rho(ud^{-1})$ for all such $\rho$, so $u_0:=ud^{-1}$ lies in $G(\breve{\QQ}_p)$.

Since $d$ commutes with $\gamma_0$, we have $u_0\gamma_0 u_0^{-1}=\gamma_0'$, and thus $u_0$ induces a bijection
\begin{align*}
u_*:B(I_0) & \to B(I_0') \\
[b] & \mapsto [u_0b\sigma(u_0)^{-1}].
\end{align*}
This bijection is independent of $d$ (and therefore deserves the name $u_*$) because any other choice of $d$ is related by an element of $I_0(\breve{\QQ}_p)$ and therefore its $\sigma$-conjugation of $b$ does not change the class in $B(I_0)$.

Now with the above setup, an \emph{isomorphism} between Kottwitz parameters $(\gamma_0,a,[b_0])$ and $(\gamma_0',a',[b_0'])$ is an element $u\in G(\overline{\QQ})$ such that
\begin{itemize}
\item $\Int(u)\gamma_0=\gamma_0'$ and $u^{-1}\rho(u)\in I_0$ for all $\rho\in \Gal(\overline{\QQ}/\QQ)$ (i.e. $u$ stably conjugates $\gamma_0$ to $\gamma_0'$),
\item the bijection $u_*:\fD(I_0,G;\AA_f^p)\to\fD(I_0',G;\AA_f^p)$ sends $a$ to $a'$, and
\item the bijection $u_*:B(I_0)\to B(I_0')$ sends $[b_0]$ to $[b_0']$.
\end{itemize}

\paragraph{}
A Kottwitz parameter $(\gamma_0,a,[b_0])$ is \emph{$\bfb$-admissible} if the map $B(I_0)\to B(G)$ sends $[b_0]$ to $\bfb$. Note that we have fixed $\bfb$ in $B(G,\mu^{-1})$, so a $\bfb$-admissible Kottwitz parameter automatically satisfies item \ref{KP0} of Definition \ref{kottwitz parameter}.

To define acceptable, let $(\gamma_0,a,[b_0])$ be a Kottwitz parameter. Since $b_0$ (any representative of $[b_0]$) lies in $I_0=G_{\gamma_0}^\circ$, we see that $\gamma_0$ centralizes $b_0$; since $\gamma_0$ is rational, we see that $\gamma_0$ furthermore $\sigma$-centralizes $b_0$, so we can regard $\gamma_0$ as an element of $J_{b_0}(\QQ_p)$.

Say that $(\gamma_0,a,[b_0])$ is \emph{acceptable} if $\gamma_0$ is acceptable as an element of $J_{b_0}(\QQ_p)$. This does not depend on the choice of representative $b_0$ of $[b_0]$, as a different representative $b_0'$---being $\sigma$-conjugate to $b_0$---will admit an isomorphism $J_{b_0}\overset{\sim}{\to}J_{b_0'}$ sending $\gamma_0$ to a conjugate, and acceptability is insensitive to conjugation.

\paragraph{}
\label{par: Kottwitz invariant}
Let $\fc=(\gamma_0,a,[b_0])$ be a Kottwitz parameter, $I_0=G_{\gamma_0}^\circ$, and consider the group
\[
\fE(I_0,G;\AA/\QQ)=\coker\Big(H^0_\ab(\AA,G)\to H^0_\ab(\AA/\QQ,I_0\to G)\Big)
\]
where $H^0_\ab$ is the abelianized Galois cohomology of \cite{AST_1999__257__R1_0}. We write $\alpha(\fc)\in\fE(I_0,G;\AA/\QQ)$ for the \emph{Kottwitz invariant} of $\fc$ as defined in \cite[1.7]{KSZ}.

\paragraph{}
\label{par: map LRP->KP}
We now define a map from LR pairs to Kottwitz parameters, following \cite[3.5.1]{KSZ}.

Let $(\phi,\varepsilon)$ be a semi-admissible LR pair, and $\tau(\phi)\in I_\phi^\ad(\AA_f^p)$. After conjugation by an element of $G(\overline{\AA}_f^p)$ we may assume that $(\phi,\varepsilon)$ is gg---we verify in Lemma \ref{lemma: strong conjugacy} that the result is insensitive to this conjugation. We will define a Kottwitz parameter $\bft(\phi,\varepsilon,\tau(\phi))=(\gamma_0,a,[b_0])$ associated to $(\phi,\varepsilon)$ and $\tau(\phi)$.

Define $\gamma_0=\varepsilon$. By the gg condition, $\varepsilon$ is contained in $G(\QQ)$ and is semi-simple and elliptic in $G(\RR)$, verifying the requirements for $\gamma_0$ in Definition \ref{kottwitz parameter}. We write $I_0=G_{\gamma_0}^\circ=G_\varepsilon^\circ$.

Next we consider $a$. Recall the cocycles $\zeta_\phi^{p,\infty}$ and $\zeta_{\phi,\ell}$ of \ref{par: X^p(phi) and cocycles}. The gg condition
\begin{equation}
\label{eq: gg condition}
\phi(q_\rho)=g_\rho\rtimes\rho\text{ has }g_\rho\in G_\varepsilon^\circ=I_0\text{ for }q_\rho\text{ any lift of }\rho\in \Gal(\overline{\QQ}/\QQ)
\end{equation}
implies that $\zeta_\phi^{p,\infty}$ is valued in $I_0(\overline{\AA}_f^p)$.

Choose a lift $\tilde{\tau}\in I_\phi(\overline{\AA}_f^p)$ of $\tau(\phi)$, and define a cocycle $A:\Gal(\overline{\QQ}/\QQ)\to I_0(\overline{\AA}_f^p)$ by
\[
A(\rho)=t_\rho\zeta_\phi^{p,\infty}(\rho)
\]
where $t_\rho=\tilde{\tau}^{-1}\rho(\tilde{\tau})\in Z_{I_\phi}(\overline{\AA}_f^p)$, acting by $\rho$ via the $\QQ$-structure of $I_\phi$---we can regard $t_\rho$ as an element of $I_0(\overline{\AA}_f^p)$ because the natural embedding $Z_{I_\phi}\to G$ factors through $I_0$.

The cocycle $A$ splits over $G(\overline{\AA}_f^p)$. To see this, we write $\tilde{\tau}_G$ for the image of $\tilde{\tau}$ in $G$---we distinguish these because $\tilde{\tau}$ is subject to the Galois action given by the $\QQ$-structure on $I_0$, while $\tilde{\tau}_G$ is subject to that given by the $\QQ$-structure on $G$. With respect to the Galois action on $G$, the element $\rho(\tilde{\tau})$ becomes $\zeta_\phi^{p,\infty}(\rho)\rho(\tilde{\tau}_G)\zeta_\phi^{p,\infty}(\rho)^{-1}$, and so
\begin{align*}
A(\rho)
&=\tilde{\tau}^{-1}\rho(\tilde{\tau})\zeta_\phi^{p,\infty}(\rho)
\\
&=\tilde{\tau}_G^{-1}\zeta_\phi^{p,\infty}(\rho)\rho(\tilde{\tau}_G)\zeta_\phi^{p,\infty}(\rho)^{-1}\zeta_\phi^{p,\infty}(\rho)
\\
&=\tilde{\tau}_G^{-1}\zeta_\phi^{p,\infty}(\rho)\rho(\tilde{\tau}_G).
\end{align*}
Combined with the fact that $\zeta_\phi^{p,\infty}$ splits in $G(\overline{\AA}_f^p)$ (realized as $\rho\mapsto x\rho(x)^{-1}$ for $x\in X^p(\phi)$), this shows that $A$ splits in $G(\overline{\AA}_f^p)$ as well.

We define $a\in \fD(I_0,G;\AA_f^p)$ in our Kottwitz parameter to be the class defined by the image of $A$. This does not depend on the choice of lift $\tilde{\tau}$, because two choices differ by an element of $Z_{I_\phi}(\overline{\AA}_f^p)$ which commutes with $\zeta_\phi^{p,\infty}$.

Finally we construct $[b_0]$. The same gg condition \eqref{eq: gg condition} above implies that $\phi$ factors
\[
\phi:\fQ\overset{\phi_0}{\longrightarrow}\fG_{I_0}\to\fG_G.
\]
Then $(\phi_0,\varepsilon)$ is again an LR pair, and taking a $p$-adic realization $(b_\theta,\varepsilon')$ produces a class $[b_0]=[b_\theta]\in B(I_0)$ (independent of the choice of $p$-adic realization, as different choices of $b_\theta$ will still be $\sigma$-conjugate in $I_0$). This finishes the construction of $\bft(\phi,\varepsilon,\tau(\phi))=(\gamma_0,a,[b_0])$.

Note that we have taken $\tau(\phi)\in I_\phi^\ad(\AA_f^p)$, but by \cite[Prop 3.5.2]{KSZ} this construction only depends on its image in $\fE^p(\phi)=I_\phi(\AA_f^p)\backslash I_\phi^\ad(\AA_f^p)$, so we have a well-defined Kottwitz parameter $\bft(\phi,\varepsilon,\tau(\phi)$ associated to a semi-admissible pair $(\phi,\varepsilon)$ and an element $\tau\in \Gamma(\fE^p)$.

We want to show that this construction only depends on the conjugacy class of the LR pair. In particular, after conjugating we have worked in the case that our LR pair is gg, so we want to show that if two gg pairs are conjugate then the resulting Kottwitz parameters are isomorphic. For this we need to assume that the $\tau$-twists are well-behaved under conjugation as well.

\begin{lemma}
\label{lemma: strong conjugacy}
Let $(\phi,\varepsilon)$ and $(\phi',\varepsilon')$ be gg LR pairs, and $\tau\in\Gamma(\fE^p)_1$. Write $\bft(\phi,\varepsilon,\tau(\phi))=(\gamma_0,a,[b_0])$ and $\bft(\phi',\varepsilon',\tau(\phi'))=(\gamma_0',a',[b_0'])$ for the associated Kottwitz parameters. If $u\in G(\overline{\QQ})$ conjugates $(\phi,\varepsilon)$ to $(\phi',\varepsilon')$, then
\begin{enumerate}
\item $u\rho(u)^{-1}\in G_{\varepsilon'}^\circ$ for all $\rho\in\Gal(\overline{\QQ}/\QQ)$, and
\item $u$ gives an isomorphism $(\gamma_0,a,[b_0])\cong(\gamma_0',a',[b_0'])$.
\end{enumerate}
\end{lemma}
The proof is essentially the same as in \cite[Prop 3.5.3]{KSZ}.

Thus for any $\tau\in\Gamma(\fE^p)_1$ we have a well-defined map
\begin{align*}
\bft:\LRP_\sa/\text{conj.} & \to \KP/\text{isom.} \\
(\phi,\varepsilon) & \mapsto \bft(\phi,\varepsilon,\tau(\phi))=(\gamma_0,a,[b_0]).
\end{align*}

\begin{lemma}
\label{lemma: LRP iff KP}
Let $(\phi,\varepsilon)$ a semi-admissible pair, $\tau(\phi)\in I_\phi^\ad(\AA_f^p)$, and let $\bft(\phi,\varepsilon,\tau(\phi))=(\gamma_0,a,[b_0])$ be the associated Kottwitz parameter.
\begin{enumerate}
\item $(\gamma_0,a,[b_0])$ is $\bfb$-admissible if and only if $(\phi,\varepsilon)$ is $\bfb$-admissible.
\item $(\gamma_0,a,[b_0])$ is acceptable if and only if $(\phi,\varepsilon)$ is acceptable.
\end{enumerate}
\end{lemma}
\begin{proof}
For the first claim, recall that a semi-admissible pair $(\phi,\varepsilon)$ is $\bfb$-admissible if $\phi(p)\circ\zeta_p:\fG_p\to\fG_G(p)$ is conjugate to an unramified morphism $\theta$ with $b_\theta\in\bfb$. On the other hand, $[b_0]$ is the class (in $B(I_0)$) defined by precisely such a $b_\theta$, and $(\gamma_0,a,[b_0])$ is defined to be $\bfb$-admissible if $[b_0]$ maps to the class $\bfb$ under $B(I_0)\to B(G)$. These conditions are manifestly equivalent.

For the second claim, by conjugation we may assume that $(\phi,\varepsilon)$ is gg (as acceptability only depends on the conjugacy or isomorphism class). Recall that $(\phi,\varepsilon)$ is acceptable if for one (hence any) $p$-adic realization $(b_\theta,\varepsilon')$ we have $\varepsilon'$ acceptable in $J_{b_\theta}(\QQ_p)$; and that $(\gamma_0,a,[b_0])$ is acceptable if $\gamma_0$ is acceptable in $J_{b_0}(\QQ_p)$. Thus it suffices to show that $(b_0,\gamma_0)$ is a $p$-adic realization of $(\phi,\varepsilon)$. This follows from the fact that $b_0$ is defined by an $I_0(\overline{\QQ}_p)$-conjugate of $\phi(p)\circ\zeta_p$, which commutes with $\gamma_0=\varepsilon$.
\end{proof}

\paragraph{}
\label{par: LR -> classical KP}

We can also associate a classical Kottwitz parameter $(\gamma_0,\gamma,\delta)$ to an LR pair. Let $(\phi,\varepsilon)$ be a gg acceptable $\bfb$-admissible LR pair and $\tau(\phi)\in I_\phi(\overline{\AA}_f^p)$. By the gg condition, $\varepsilon$ lies in $G(\QQ)$, so we can let $\gamma_0=\varepsilon$. Then we define $\gamma$ and $\delta$ as the image of $\Int(\tau(\phi))\varepsilon$ in $G(\AA_f^p)$ and $J_b(\QQ_p)$ respectively.

\paragraph{}
\label{par: KP -> classical KP}
Next we associate a classical Kottwitz parameter to a Kottwitz parameter. Let $(\gamma_0,a,[b_0])$ be an acceptable $\bfb$-admissible Kottwitz parameter. We let $\gamma_0=\gamma_0$, i.e. the element $\gamma_0$ of our classical Kottwitz parameter we choose to be the element $\gamma_0$ of our Kottwitz parameter. The class $a$ determines a conjugacy class in $G(\AA_f^p)$ stably conjugate to $\gamma_0$, and we choose $\gamma$ to be an arbitrary element of this class. Finally, as our Kottwitz parameter is acceptable and $\bfb$-admissible, $\gamma_0$ can be considered as an element of $J_b(\QQ_p)$, and we set $\delta=\gamma_0$.

\begin{lemma}
\label{lem: same classical KP}
Let $(\phi,\varepsilon)$ a gg acceptable $\bfb$-admissible pair, $\tau(\phi)\in I_\phi(\overline{\AA}_f^p)$, and $\bft(\phi,\varepsilon,\tau(\phi))=(\gamma_0,a,[b_0])$ the associated Kottwitz parameter. Let $(\gamma_0,\gamma,\delta)$ be the classical Kottwitz parameter associated to $(\gamma_0,a,[b_0])$, and $(\gamma_0',\gamma',\delta')$ the classical Kottwitz parameter associated to $(\phi,\varepsilon)$ and $\tau(\phi)$. Then $(\gamma_0,\gamma,\delta)$ and $(\gamma_0',\gamma',\delta')$ are equivalent.
\end{lemma}
\begin{proof}
By construction we have $\gamma_0'=\varepsilon=\gamma_0$---in particular $\gamma_0$ is stably conjugate to $\gamma_0'$, which verifies the first condition of equivalence. Examining the cocycle
\[
A(\rho)=\tau(\phi)^{-1}\rho(\tau(\phi))\zeta_\phi^{p,\infty}(\rho)
\]
defining $a$, and noting that $\zeta_\phi^{p,\infty}$ commutes with $\varepsilon=\gamma_0$, we see that the conjugacy class of $\gamma$ is the conjugacy class of $\Int(\tau(\phi))$, which shows that $\gamma$ and $\gamma'$ are conjugate in $G(\AA_f^p)$. By construction, $\delta$ and $\delta'$ are conjugate as well.
\end{proof}

\subsection{The Image of the Map $\LRP\to\KP$}
\label{sec: image of LRP -> KP}

Recall from \S\ref{LR/KP section} we have defined a map $\bft:\LRP_\sa/\text{conj.} \to \KP/\text{isom.}$ for any $\tau\in\Gamma(\fE^p)_1$. The main result of this section is that the image of the set of acceptable $\bfb$-admissible LR pairs is the set of acceptable $\bfb$-admissible Kottwitz parameters with trivial Kottwitz invariant.

We begin by observing that the image of acceptable $\bfb$-admissible LR pairs lies in the set of acceptable $\bfb$-admissible Kottwitz parameters with trivial Kottwitz invariant. Let $(\phi,\varepsilon)$ an acceptable $\bfb$-admissible LR pair and $\tau\in\Gamma(\fE^p)$ a tori-rational element. By Lemma \ref{lemma: LRP iff KP}, its image $\bft(\phi,\varepsilon,\tau(\phi))$ is also acceptable and $\bfb$-admissible. The following lemma tells us that its image also has trivial Kottwitz invariant.

\begin{lemma}[{\cite[Prop. 3.6.3]{KSZ}}]
\label{lem: trivial invariant}
Let $(\phi,\varepsilon)$ be a semi-admissible LR pair, $\tau\in\Gamma(\fE^p)_1$ a tori-rational element, and $\fc=\bft(\phi,\varepsilon,\tau(\phi))$ the associated Kottwitz parameter. Then the Kottwitz invariant $\alpha(\fc)$ is zero.
\end{lemma}

\paragraph{}
Now we want to show that given an acceptable $\bfb$-admissible Kottwitz parameter $(\gamma_0,a,[b_0])$ and $\tau\in\Gamma(\fE^p)_0$ tori-rational, we can find an acceptable $\bfb$-admissible LR pair $(\phi,\varepsilon)$ such that $\bft(\phi,\varepsilon,\tau(\phi))=(\gamma_0,a,[b_0])$. 

If $\mu$ is a cocharacter of an algebraic group $H$, over $\QQ_p$, we define $[b_\bas(\mu)]\in B(H)$ to be the unique basic class in $B(H,\mu)$.

\begin{lemma}
\label{lem: KMPS analogue}
Let $(\gamma_0,a,[b_0])$ be an acceptable $\bfb$-admissible Kottwitz parameter, and write $I_0=G_{\gamma_0}^\circ$. Then there exists a maximal torus $T\subset I_0$ over $\QQ$ and $x\in X$ such that $h_x$ factors through $T_\RR$ (so $\mu_x\in X_*(T)$) and $[b_\bas(\mu_x^{-1})]\in B(T)$ maps to $[b_0]$ in $B(I_0)$.
\end{lemma}
This lemma and its proof are directly inspired from \cite[Prop 1.2.5]{KMPS}.
\begin{proof}
Since $\gamma_0$ is acceptable with respect to $[b_0]$, we have $I_0\subset M_{b_0}$ by Lemma \ref{not root argument}. This implies that $I_0$ centralizes the slope homomorphism of $[b_0]$, i.e. $[b_0]$ is basic in $B(I_0)$. Let $J_{b_0}^{I_0}$ be the inner form of $I_0$ defined by $[b_0]$.

Recall we can consider $\gamma_0$ as an element of $J_{b_0}(\QQ_p)$. The group $J_{b_0}^{I_0}$ is the connected centralizer of $\gamma_0$ in $J_{b_0}$. Choose an elliptic maximal torus $T'\subset J_{b_0}$ containing $\gamma_0$. Then $T'$ is contained in $J_{b_0}^{I_0}$, and being elliptic, transfers to a torus $T'\overset{\sim}{\to} T_p\subset I_0$ over $\QQ_p$.

By \cite[Cor 1.1.17]{KMPS}, there is a representative $\mu_p\in X_*(T_p)$ of $\{\mu_X\}$ (the conjugacy class of cocharacters induced by the Shimura datum $(G,X)$) such that $[b_\bas(\mu_p^{-1})]\in B(T_p)$ maps to $[b_0]$ in $B(I_0)$ (note that hypothesis (1.1.3.1) of that Corollary is satisfied because $[b_0]$ is basic in $B(I_0)$, and $[b_0]$ is $\mu^{-1}$-admissible because it maps to $\bfb\in B(G,\mu^{-1})$).

Choose a maximal torus $T_\infty\subset G_\RR$ containing $\gamma_0$, so in particular $T_\infty\subset I_{0,\RR}$. Since every element of $X$ factors through an elliptic maximal torus of $G_\RR$ and all such tori are $G(\RR)$-conjugate, we can choose $x\in X$ such that $h_x$ factors through $T_\infty$.

Having produced $T_p\subset I_{0,\QQ_p}$ and $T_\infty\subset I_{0,\RR}$ with the necessary properties, the rest of the proof can be read verbatim from the proof of \cite[Prop 1.2.5]{KMPS}, replacing their $G$ by our $I_0$.
\end{proof}

\begin{proposition}
\label{prop: hit gamma_0}
Let $(\gamma_0,a,[b_0])$ be an acceptable $\bfb$-admissible Kottwitz parameter. Then there is an admissible morphism $\phi_0$ which forms a gg acceptable $\bfb$-admissible LR pair $(\phi_0,\gamma_0)$.
\end{proposition}
\begin{proof}
By Lemma \ref{lem: KMPS analogue}, there exists a maximal torus $T\subset I_0$ over $\QQ$ and $x\in X$ such that $h_x$ factors through $T_\RR$ and $[b_\bas(\mu_x^{-1})]\in B(T)$ maps to $[b_0]\in B(I_0)$. In particular, $(T,h=h_x)$ forms a special point datum.

Let $\phi_0=i\circ\psi_{\mu_h}$ be the admissible morphism induced from the special point datum $(T,h)$ as in \ref{par: special point datum}. Since $T\subset I_0$ is maximal and $\gamma_0$ is central in $I_0$, we have $\gamma_0\in T$, so we can form the LR pair $(\phi_0,\gamma_0)$. By construction it is very special and in particular gg. It remains to show that this LR pair is acceptable and $\bfb$-admissible.

We check $\bfb$-admissibility first. Consider $\psi_{\mu_h}:\fQ\to \fG_T$ and its $p$-part $\psi_{\mu_h}(p)\circ\zeta_p:\fG_p\to\fG_T(p)$. As in \ref{unramified paragraph}, the latter is conjugate by some element $y\in T(\overline{\QQ}_p)$ to an unramified morphism $\theta:\fG_p\to\fG_T(p)$ which gives rise to an element $[b_\theta]\in B(T)$.

By \cite[Lemma 2.2.10]{KSZ}, the image $\kappa_T(b_\theta)\in X_*(T)_{\Gamma_p}$ of $[b_\theta]$ under the Kottwitz map is equal to the image of $\mu_h^{-1}\in X_*(T)$. Since the Kottwitz map is a bijection in the case of a torus, we conclude that $[b_\theta]=[b_\bas(\mu_h^{-1})]$ in $B(T)$, as the latter also maps to $\mu_h^{-1}$. By our conclusions from Lemma \ref{lem: KMPS analogue} above we find that $[b_\theta]\in B(T)$ maps to $[b_0]\in B(I_0)$ and therefore to $\bfb\in B(G)$ (by $\bfb$-admissibility of our Kottwitz parameter). 

Now conjugating $\phi_0(p)\circ\zeta_p=i\circ\psi_{\mu_h}(p)\circ\zeta_p$ by the same element $y\in T(\overline{\QQ}_p)$ produces the unramified morphism $i\circ\theta$, and the associated class $[b_{i\circ\theta}]\in B(G)$ is the image of the class $[b_\theta]\in B(T)$. Thus we have $[b_{i\circ\theta}]=\bfb$, and this verifies that $(\phi_0,\gamma_0)$ is $\bfb$-admissible.

Now we check acceptability. Note that we conjugated by an element $y\in T(\overline{\QQ}_p)$ to produce the unramified morphism $i\circ\theta$. In particular $y$ commutes with $\gamma_0$, so $(b_{i\circ\theta},\gamma_0)$ is a $p$-adic realization of $(\phi_0,\gamma_0)$. We have assumed that $\gamma_0$ is acceptable with respect to $[b_0]$, and shown that $[b_0]=[b_\theta]$ in $B(I_0)$. This shows that $\gamma_0$ is acceptable with respect to $[b_{i\circ\theta}]$, and $(\phi_0,\gamma_0)$ is acceptable.
\end{proof}

\begin{lemma}
\label{centralizers coincide}
Let $(\phi,\varepsilon)$ be an acceptable LR pair, and suppose that $\phi^\Delta$ is defined over $\QQ$ and $\varepsilon\in G(\QQ)$ (in particular, this applies to any acceptable gg pair). Then the inclusion $I_{\phi,\varepsilon}\hookrightarrow G_\varepsilon$ over $\overline{\QQ}$ is an isomorphism.
\end{lemma}
\begin{proof}
Recall from \ref{lem: gerb kernel facts} that $I_{\phi,\overline{\QQ}}$ is the centralizer of (the image of) $\phi^\Delta$ in $G_{\overline{\QQ}}$, so our goal is to show that any element commuting with $\varepsilon$ must commute with $\phi^\Delta$.

We begin by showing that any element commuting with $\varepsilon$ must commute with $\phi^\Delta\circ\nu(p)$ and $\phi^\Delta\circ\nu(\infty)$.

At $p$: let $g\in G(\overline{\QQ}_p)$ such that $\theta=\Int(g)\circ\phi(p)\circ\zeta_p$ is unramified, and let $b_\theta\in G(\QQ_p^\ur)$ defined by $\theta(d_\sigma)=b_\theta\rtimes\sigma$. Write $\varepsilon'=\Int(g)\varepsilon$. By Lemma \ref{unramified morphism facts} we have
\[
\Int(g)\circ\phi^\Delta\circ\nu(p)=(\Int(g)\circ\phi(p)\circ\zeta_p)^\Delta=-\nu_{b_\theta},
\]
so the centralizer of $\Int(g)\circ\phi^\Delta\circ\nu(p)$ is $M_{b_\theta}$. On the other hand, $(b_\theta,\varepsilon')$ is a $p$-adic realization of our acceptable pair $(\phi,\varepsilon)$---in particular, we can apply Lemma \ref{not root argument} to conclude that $G_{\varepsilon'}\subset M_{b_\theta}$.

Since the centralizer of $\Int(g)\circ\phi^\Delta\circ\nu(p)$ is $M_{b_\theta}$ we see the centralizer of $\phi^\Delta\circ\nu(p)$ is $\Int(g^{-1})M_{b_\theta}$, and likewise we have $G_\varepsilon=\Int(g^{-1})G_{\varepsilon'}$. The above analysis then shows that $\Int(g^{-1})G_{\varepsilon'}=\Int(g^{-1})M_{b_\theta}$, which is to say that the centralizer of $\varepsilon$ is contained in the centralizer of $\phi^\Delta\circ\nu(p)$, as desired.

At $\infty$: as in the proof of \cite[Lemma 3.1.9]{KSZ}, the fact is that $\phi^\Delta\circ\nu(\infty)$ is central in $G$, and therefore any element commuting with $\varepsilon$ trivially commutes with $\phi^\Delta\circ\nu(\infty)$.

Now, suppose that $g\in G(\overline{\QQ})$ commutes with $\varepsilon$, and we want to see that $g$ commutes with $\phi^\Delta$.

Recall that $\phi^\Delta$ is a morphism $Q\to G$ where $Q=\ilim_L Q^L$ is the kernel of $\fQ$. For each finite Galois $L/\QQ$, the torus $Q^L$ is generated by the $\Gal(L/\QQ)$-conjugates of the images of $\nu(p)^L$ and $\nu(\infty)^L$. Thus $Q$ is generated by the $\Gal(\overline{\QQ}/\QQ)$-conjugates of the images of $\nu(p)$ and $\nu(\infty)$.

For any $\rho\in\Gal(\overline{\QQ}/\QQ)$, the conjugate $\rho(g)$ again commutes with $\varepsilon$ by our hypothesis that $\varepsilon$ is rational, and therefore by the above arguments $\rho(g)$ commutes with $\phi^\Delta\circ\nu(v)$ for $v=p,\infty$. Applying $\rho^{-1}$ and using our hypothesis that $\phi^\Delta$ is defined over $\QQ$, we see that $g$ commutes with $\phi^\Delta\circ\rho^{-1}(\nu(v))$ for $v=p,\infty$. Since $\rho$ was arbitrary and the $\Gal(\overline{\QQ}/\QQ)$-conjugates of the images of $\nu(v)$ generate $Q$, this implies that $g$ commutes with $\phi^\Delta$, as desired.
\end{proof}

\paragraph{}
\label{compatible inner twists}
Let $(\phi,\varepsilon)$ be a gg acceptable pair, so that letting $q_\rho\in\fQ$ a lift of $\rho\in\Gal(\overline{\QQ}/\QQ)$, we have $\phi(q_\rho)=g_\rho\rtimes\rho$ with $g_\rho\in G_\varepsilon^\circ$.

The group $I_\phi$ is defined as an inner form of $Z_G(\phi^\Delta)$ by the cocycle
\[
\rho\mapsto \Int(g_\rho)\in\Aut((G_\varepsilon^\circ)_ {\overline{\QQ}})\qquad\rho\in\Gal(\overline{\QQ}/\QQ).
\]
In view of the isomorphism of Lemma \ref{centralizers coincide}, this produces---for any gg acceptable LR pair $(\phi,\varepsilon)$---compatible inner twists
\begin{align*}
(I_{\phi,\varepsilon}^\circ)_{\overline{\QQ}} & \overset{\sim}{\longrightarrow} (G_{\varepsilon}^\circ)_{\overline{\QQ}}, \\
(I_{\phi,\varepsilon})_{\overline{\QQ}} & \overset{\sim}{\longrightarrow} (G_{\varepsilon})_{\overline{\QQ}}
\end{align*}
defined by the same cocycle.

\begin{lemma}
\label{same nu---KP version}
Let $(\gamma_0,a,[b_0])$ and $(\gamma_0,a',[b_1])$ be acceptable $\bfb$-admissible Kottwitz parameters with the same element $\gamma_0$. Then $\nu_{b_0}=\nu_{b_1}$.
\end{lemma}
\begin{proof}
Both Kottwitz parameters are assumed to be $\bfb$-admissible, so $[b_0]$ and $[b_1]$ both map to $\bfb$ in $B(G)$. In particular, $b_0$ and $b_1$ are $\sigma$-conjugate in $G(L)$. Furthermore, the semi-simple element $\gamma_0$ as an element of $G(L)$ lies in both $J_{b_0}(\QQ_p)$ and $J_{b_1}(\QQ_p)$, and is acceptable with respect to both. Thus we are in the situation of Lemma \ref{same nu}, and we conclude that $\nu_{b_0}=\nu_{b_1}$.
\end{proof}

\begin{proposition}
\label{prop: hit whole KP}
Let $(\gamma_0,a,[b_0])$ be an acceptable $\bfb$-admissible Kottwitz parameter with trivial Kottwitz invariant, and suppose there is a gg acceptable $\bfb$-admissible pair $(\phi_0,\gamma_0)$. Then there is a gg acceptable $\bfb$-admissible pair $(\phi_1,\varepsilon_1)$ with
\[
\bft(\phi_1,\varepsilon_1,1)=(\gamma_0,a,[b_0]).
\]
\end{proposition}
This is Theorem 3.5.9 of \cite{KSZ}, except that their ``$p^n$-admissible'' hypothesis has been replaced by our ``acceptable $\bfb$-admissible'' hypothesis. We briefly sketch how their proof carries over to our case.
\begin{proof}
Write $\bft(\phi_0,\gamma_0,1)=(\gamma_0,a',[b_0'])$, and $I_0=G_{\gamma_0}^\circ$. By Lemma \ref{lemma: LRP iff KP}, our hypothesis that $(\phi_0,\gamma_0)$ is acceptable and $\bfb$-admissible implies that $(\gamma_0,a',[b_0'])$ is acceptable and $\bfb$-admissible.

Since $\gamma_0$ is acceptable with respect to $[b_0']$, by Lemma \ref{not root argument} we have $I_0\subset M_{b_0'}$, and therefore $\nu_{b_0'}$ is central in $I_0$. This is the first ingredient; we also need the fact from Lemma \ref{same nu---KP version} that $\nu_{b_0}=\nu_{b_0'}$; and the fact from \ref{compatible inner twists} of a canonical inner twisting $(I_{\phi_0,\gamma_0}^\circ)_{\overline{\QQ}}\overset{\sim}{\to} I_{0,\overline{\QQ}}$.

In the presence of these three ingredients, the proof of \cite[Thm 3.5.9]{KSZ} carries over without modification to show the existence of a gg semi-admissible LR pair $(\phi_1,\varepsilon_1)$ with $\bft(\phi_1,\varepsilon_1,1)=(\gamma_0,a,[b_0])$. By Lemma \ref{lemma: LRP iff KP}, $(\phi_1,\varepsilon_1)$ is acceptable and $\bfb$-admissible.
\end{proof}

The last step is to incorporate $\tau$-twists. We collect the full result in the following proposition.

\begin{proposition}
\label{prop: surjectivity}
Let $(\gamma_0,a,[b_0])$ be an acceptable $\bfb$-admissible Kottwitz parameter with trivial Kottwitz invariant, and $\tau\in\Gamma(\fE^p)_0$ a tori-rational element. Then there is a gg acceptable $\bfb$-admissible LR pair $(\phi,\varepsilon)$ such that
\[
\bft(\phi,\varepsilon,\tau(\phi))=(\gamma_0,a,[b_0]).
\]
\end{proposition}
\begin{proof}
This is our analogue of \cite[Prop 3.6.5]{KSZ}. The proof in their case applies equally well here, with our Propositions \ref{prop: hit gamma_0} and \ref{prop: hit whole KP} replacing their and Propositions 3.4.8 and 3.5.9 respectively.

The resulting LR pair $(\phi,\varepsilon)$ is furthermore acceptable and $\bfb$-admissible by Lemma \ref{lemma: LRP iff KP}.
\end{proof}

\paragraph{}
\label{par: image of LRP -> KP}
Having completed this result, and combining it with the discussion at the beginning of this section, we conclude that under the map (for $\tau\in\Gamma(\fE^p)_0$ tori-rational)
\[
\bft:\LRP_\sa/\text{conj.}\to\KP/\text{isom.}
\]
the image of the set of acceptable $\bfb$-admissible LR pairs is the set of acceptable $\bfb$-admissible Kottwitz parameters with trivial Kottwitz invariant.

\subsection{Point-Counting Formula}
\label{sec: final point counting}

Recall the result of our preliminary point counting as summarized in Proposition \ref{preliminary formula}, updated to the language of LR pairs. Let $\tau\in\Gamma(\cH)_0$ a tori rational element satisfying Theorem \ref{igusa LR}, which we may lift to a tori-rational element of $\Gamma(\fE^p)_0$ (still called $\tau$ by abuse of notation). For any acceptable function $f\in C^\infty_c(G(\AA_f^p)\times J_b(\QQ_p))$ we have
\begin{equation}
\label{eq: final counting start}
\tr(f\mid H_c(\Ig_\Sigma,\sL_\xi))=\sum_{(\phi,\varepsilon)} \frac{\vol\left(I^\circ_{\phi,\varepsilon}(\QQ)\backslash G_{b,\gamma\times\delta}^\circ\right)}{[I_{\phi,\varepsilon}(\QQ):I^\circ_{\phi,\varepsilon}(\QQ)]} O_{\gamma\times\delta}^{G(\AA_f^p)\times J_b(\QQ_p)}(f)\tr(\xi(\varepsilon)),
\end{equation}
where $(\phi,\varepsilon)$ ranges over conjugacy classes of $\bfb$-admissible LR pairs, and $\gamma,\delta$ are the elements appearing in the classical Kottwitz parameter associated to $(\phi,\varepsilon)$ and $\tau(\phi)$. In fact, as a first modification, we can restrict the sum to pairs $(\phi,\varepsilon)$ which are furthermore acceptable: since an acceptable function $f$ is supported on acceptable elements of $J_b(\QQ_p)$, the orbital integral will be zero unless $(\phi,\varepsilon)$ is acceptable.

\paragraph{}
\label{par: I_c}
Let $\fc=(\gamma_0,a,[b_0])$ be a Kottwitz parameter with $[b_0]\in B(I_0)$ basic, and with trivial Kottwitz parameter (note these hypotheses are satisfied by any Kottwitz parameter arising from an LR pair). We can define an inner form $I_\fc$ of $I_0=G_{\gamma_0}^\circ$ as follows. Writing $a=(a_\ell)_{\ell\neq p,\infty}$, let $I_\ell$ be the inner form of $I_0$ over $\QQ_\ell$ defined by $a_\ell$ (or to be precise, the image of $a_\ell$ in $H^1(\QQ_\ell,I_0^\ad)$). At $p$, let $I_p=J_{b_0}^{I_0}$ be the inner form of $I_0$ over $\QQ_p$ defined by the basic class $[b_0]\in B(I_0)$. At $\infty$, let $I_\infty$ be the inner form of $I_0$ over $\RR$ which is compact modulo $Z_G$. By \cite[Prop 1.7.12]{KSZ}, these local components determine a unique inner form $I_c$ of $I_0$ over $\QQ$ such that $I_\fc\otimes\QQ_v\cong I_v$ for all places $v$ of $\QQ$.

This is the group that will replace $I_{\phi,\varepsilon}$ in the volume term.

\paragraph{}
\label{par: orbital integrals agree}
Next we examine the orbital integral terms. Recall that in the setting of our preliminary formula \ref{eq: final counting start} we have defined $\gamma\times\delta\in G_b=G(\AA_f^p)\times J_b(\QQ_p)$ to be the image of $\Int(\tau(\phi))\varepsilon$; or in the language of \ref{par: LR -> classical KP}, $\gamma,\delta$ come from the classical Kottwitz parameter $(\gamma_0,\gamma,\delta)$ associated to $(\phi,\varepsilon)$ and $\tau(\phi)$. If $\bft(\phi,\varepsilon,\tau(\phi))=(\gamma_0,a,[b_0])$ is the Kottwitz parameter associated to our LR pair, then as in \ref{par: KP -> classical KP} we can construct a classical Kottwitz parameter from $(\gamma_0,a,[b_0])$, and by Lemma \ref{lem: same classical KP} this classical Kottwitz parameter produces the same elements $\gamma,\delta$ (up conjugacy in $G(\AA_f^p)$ or $J_b(\QQ_p)$ respectively). Thus, orbital integrals defined in terms of classical Kottwitz parameters on both sides will agree.

Define
\[
T_\xi^f(\phi,\varepsilon,\tau(\phi))
=
\vol\left(I^\circ_{\phi,\varepsilon}(\QQ)\backslash G_{b,\gamma\times\delta}^\circ\right) O_{\gamma\times\delta}^{G(\AA_f^p)\times J_b(\QQ_p)}(f)\tr(\xi(\varepsilon))
\]
where $\gamma,\delta$ are the elements in the classical Kottwitz parameter associated to $(\phi,\varepsilon)$ and $\tau(\phi)$; and
\[
T_\xi^f(\gamma_0,a,[b_0])
=
\vol\left(I^\circ_{\fc}(\QQ)\backslash I^\circ_{\fc}(\AA_f)\right) O_{\gamma\times\delta}^{G(\AA_f^p)\times J_b(\QQ_p)}(f)\tr(\xi(\gamma_0))
\]
where $\gamma,\delta$ are the elements in the classical Kottwitz parameter associated to $(\gamma_0,a,[b_0])$, and $I_\fc$ is the group associated to $\fc=(\gamma_0,a,[b_0])$ as in \ref{par: I_c}.

Our discussions thus far have demonstrated the following Lemma.

\begin{lemma}
\label{lem: same T_xi^f}
Let $(\phi,\varepsilon)$ be a gg acceptable $\bfb$-admissible LR pair, $\tau(\phi)\in I_\phi^\ad(\AA_f^p)$, and $\bft(\phi,\varepsilon,\tau(\phi))=(\gamma_0,a,[b_0])$ the associated Kottwitz parameter. Then
\[
T_\xi^f(\phi,\varepsilon,\tau(\phi))
=
T_\xi^f(\gamma_0,a,[b_0]).
\]
\end{lemma}

\paragraph{}
\label{par: twisting LR pairs}
We now return to our study of the map $\bft:\LRP/\text{conj.}\to\KP/\text{isom}$. We showed in \S\ref{sec: image of LRP -> KP} that the image of the set of acceptable $\bfb$-admissible LR pairs is the set of acceptable $\bfb$-admissible Kottwitz parameters with trivial Kottwitz invariant. Our next goal is to examine the fibers of this map; we will understand them in terms of cohomological twists.

Recall from \ref{lem: gerb kernel facts} that we can twist a morphism $\phi:\fQ\to\fG_G$ by a cocycle $e\in Z^1(\QQ,I_\phi)$ to get another morphism $e\phi$. We saw in \ref{prop: admissible twist} that if $\phi$ is admissible, then $e\phi$ is again admissible exactly when $e$ lies in $\Sha_G^{p,\infty}(\QQ,I_\phi)$.

We now consider twisting LR pairs. For an LR pair $(\phi,\varepsilon)$ and cocycle $e\in Z^1(\QQ,I_{\phi,\varepsilon})\subset Z^1(\QQ,I_\phi)$, we can define a twist by simply twisting the morphism $(e\phi,\varepsilon)$. By \cite[Lemma 3.2.5]{KSZ} this is again an LR pair---essentially we need to take a cocycle in $I_{\phi,\varepsilon}$ rather than $I_\phi$ to ensure that $\varepsilon$ still lies in $I_\phi$. As in the case of twisting morphisms, two twists $(e\phi,\varepsilon)$ and $(e'\phi,\varepsilon)$ are conjugate by $G(\overline{\QQ})$ exactly when $e,e'$ define the same class in $H^1(\QQ,I_{\phi,\varepsilon})$ (\cite[Lemma 3.2.6]{KSZ}).

Now, suppose that $(\phi,\varepsilon)$ is gg. Then \cite[Lemma 3.2.5]{KSZ} also tells us that if $e\in Z^1(\QQ,I_{\phi,\varepsilon}^\circ)\subset Z^1(\QQ,I_{\phi,\varepsilon})$ then the twist $(e\phi,\varepsilon)$ is again gg. This gives us, for any gg pair $(\phi,\varepsilon)$, a map
\[
Z^1(\QQ,I_{\phi,\varepsilon}^\circ)\to\LRP^\gg.
\]
If we only consider LR pairs up to conjugacy then this map factors through $H^1(\QQ,I_{\phi,\varepsilon}^\circ)$, giving
\[
H^1(\QQ,I_{\phi,\varepsilon}^\circ)\to\LRP^\gg/\text{conj}.
\]

The point of all this is that each fiber of the map $\bft$ is contained in the image of one such map---that is, two elements in the same fiber are always related by an $H^1$-twist. In order to show this, we need some preparatory lemmas.

\begin{lemma}
\label{lem: stably conjugate epsilon}
Let $(\phi,\varepsilon)$ a gg acceptable LR pair, and $\Int(g)\varepsilon\in G(\QQ)$ a rational element stably conjugate to $\varepsilon$, i.e. $g\in G(\overline{\QQ})$ and $g^{-1}\rho(g)\in G_\varepsilon^\circ$ for $\rho\in\Gal(\overline{\QQ}/\QQ)$. Then the conjugate $(\Int(g)\circ\phi,\Int(g)\varepsilon)$ is again gg and acceptable.
\end{lemma}
\begin{proof}
Acceptability is insensitive to conjugacy, so $(\Int(g)\circ\phi,\Int(g)\varepsilon)$ is acceptable. That it is also gg is proven as in \cite[Lemma 3.6.4]{KSZ}, with our \ref{compatible inner twists} replacing their 3.2.14.
\end{proof}

\begin{lemma}
\label{lem: same epsilon => same delta}
Let $(\phi,\varepsilon)$ and $(\phi',\varepsilon)$ be gg acceptable $\bfb$-admissible LR pairs with the same element $\varepsilon\in G(\QQ)$. Then $\phi^\Delta=\phi'^\Delta$.
\end{lemma}
\begin{proof}
By the same reasoning as the last paragraphs of the proof of Lemma \ref{centralizers coincide}, it suffices to show that $\phi^\Delta\circ\nu(v)=\phi'^\Delta\circ\nu(v)$ for $v=p,\infty$ (this implies that $\phi^\Delta$ and $\phi'^\Delta$ agree on a generating set for $\fQ^\Delta=Q$, and being morphisms they must then agree everywhere).

At $\infty$: by \cite[Lemma 3.1.9]{KSZ}, simply the fact that both $\phi$ and $\phi'$ are admissible implies that $\phi^\Delta\circ\nu(\infty)=\phi'^\Delta\circ\nu(\infty)$.

At $p$: since $(\phi,\varepsilon)$ is gg, we can factor
\[
\phi:\fQ\to\fG_{G_\varepsilon^\circ}\to\fG_G,
\]
and therefore we can conjugate $\phi(p)\circ\zeta_p$ to an unramified morphism $\theta$ by an element $u\in G_\varepsilon^\circ(\overline{\QQ})$. Let $b_\theta\in G_\varepsilon^\circ(\QQ_p^\ur)$ be the corresponding element as in \ref{unramified morphism}. In the same way we can conjugate $\phi'(p)\circ\zeta_p$ by an element $u'\in G_\varepsilon^\circ(\overline{\QQ})$ and produce an element $b_\theta'\in G_\varepsilon^\circ(\QQ_p^\ur)$. Then $(b_\theta,\varepsilon)$ is a $p$-adic realization of $(\phi,\varepsilon)$, and $(b_\theta',\varepsilon)$ is a $p$-adic realization of $(\phi',\varepsilon)$.

We have assumed that our LR pairs are acceptable and $\bfb$-admissible, which implies that $b_\theta,b_\theta'$ and $\varepsilon$ satisfy the hypotheses of Lemma \ref{same nu}, and we conclude that $\nu_{b_\theta}=\nu_{b_\theta'}$.

On the other hand, we have
\begin{align*}
\Int(u)\circ\phi^\Delta\circ\nu(p)
& = (\Int(u)\circ\phi(p)\circ\zeta_p)^\Delta
\overset{*}{=} -\nu_{b_\theta} \\
& = -\nu_{b_\theta'}
\overset{*}{=} (\Int(u')\circ\phi'(p)\circ\zeta_p)^\Delta
= \Int(u')\circ\phi'^\Delta\circ\nu(p)
\end{align*}
where the starred equalities are given by Lemma \ref{unramified morphism facts}. By our acceptable hypothesis, we can apply Lemma \ref{not root argument} to see that $G_\varepsilon$ commutes with $\nu_{b_\theta}=\nu_{b_\theta'}$. But the above equation demonstrates that $\phi^\Delta\circ\nu(p)$ and $-\nu_{b_\theta}$ and $-\nu_{b_\theta'}$ and $\phi'^\Delta\circ\nu(p)$ are all conjugate by $G_\varepsilon^\circ(\overline{\QQ})$, so they must all be equal, and in particular $\phi^\Delta\circ\nu(p)=\phi'^\Delta\circ\nu(p)$.
\end{proof}

Now we are prepared to show that points in the same fiber of $\bft$ are related by an $H^1$-twist.

\begin{lemma}
\label{lem: H^1 twist}
Suppose that $(\phi,\varepsilon)$ and $(\phi',\varepsilon')$ are gg acceptable $\bfb$-admissible LR pairs which give rise to isomorphic Kottwitz parameters. Then the conjugacy classes of $(\phi,\varepsilon)$ and $(\phi',\varepsilon')$ are related by twisting by an element of $H^1(\QQ,I_{\phi,\varepsilon}^\circ)$.
\end{lemma}
\begin{proof}
In fact we only need to assume that the rational elements, say $\gamma_0$ and $\gamma_0'$, appearing in the two Kottwitz parameters are stably conjugate (which follows from the isomorphism of Kottwitz parameters)---this justifies our neglect of $\tau$-twists in the statement, as $\tau$-twists do not affect the rational element $\gamma_0$ of the Kottwitz parameter.

So our hypothesis implies that $\varepsilon$ and $\varepsilon'$ are stably conjugate, and by Lemma \ref{lem: stably conjugate epsilon} we can conjugate $(\phi',\varepsilon')$ to a gg pair $(\phi_0,\varepsilon)$ which is again acceptable and $\bfb$-admissible.

By Lemma \ref{lem: same epsilon => same delta} we have $\phi^\Delta=\phi_0^\Delta$. As in Lemma \ref{lem: gerb kernel facts} we can choose $e\in Z^1(\QQ,I_\phi)$ so that $\phi_0=e\phi$. Now write
\[
\phi(q_\rho)=g_\rho\rtimes\rho,
\qquad\qquad
\phi_0(q_\rho)=e\phi(q_\rho)=e_\rho g_\rho\rtimes\rho
\]
where as usual $\rho\in\Gal(\overline{\QQ}/\QQ)$ and $q_\rho\in\fQ$ is a lift of $\rho$. Since our LR pairs are gg, we have $g_\rho\in G_\varepsilon^\circ$ and $e_\rho g_\rho\in G_\varepsilon^\circ$, so we conclude that $e_\rho\in G_\varepsilon^\circ$. By \ref{centralizers coincide} this shows that in fact $e\in Z^1(\QQ,I_{\phi,\varepsilon}^\circ)$, demonstrating that the conjugacy classes of $(\phi,\varepsilon)$ and $(\phi',\varepsilon')$ are related by twisting by $H^1(\QQ,I_{\phi,\varepsilon}^\circ)$.
\end{proof}

We also have an analogue of Proposition \ref{prop: admissible twist} characterizing which twists of a semi-admissible LR pair are semi-admissible.

\begin{proposition}[{\cite[Prop 3.2.19]{KSZ}}]
\label{prop: semi-admissible twist}
If $(\phi,\varepsilon)$ is gg and semi-admissible and $e\in Z^1(\QQ,I_{\phi,\varepsilon}^\circ)$, then the twist $(e\phi,\varepsilon)$ is gg and semi-admissible exactly when $e$ lies in $\Sha_G^{\infty}(\QQ,I_{\phi,\varepsilon}^\circ)$.
\end{proposition}

Combining Proposition \ref{prop: semi-admissible twist} with the discussion of \ref{par: twisting LR pairs}, we have for each gg acceptable $\bfb$-admissible pair $(\phi,\varepsilon)$ a map
\[
\eta_{\phi,\varepsilon}:\Sha_G^\infty(\QQ,I_{\phi,\varepsilon}^\circ)\to\LRP^\gg_\sa/\text{conj.}
\]
and by Lemma \ref{lem: H^1 twist} any two such LR pairs giving rise to isomorphic Kottwitz parameters must both lie in the image of one such map $\eta_{\phi,\varepsilon}$.

\paragraph{}
\label{par: fibers of LRP -> KP}
Let $\tau\in\Gamma(\fE^p)_0$ tori-rational, and $(\gamma_0,a,[b_0])$ an acceptable $\bfb$-admissible Kottwitz parameter with trivial Kottwitz invariant. By \ref{prop: surjectivity}, there is a gg acceptable $\bfb$-admissible LR pair $(\phi_0,\gamma_0)$ such that $\bft(\phi_0,\gamma_0,\tau(\phi_0))=(\gamma_0,a,[b_0])$. Furthermore combining Proposition \ref{prop: semi-admissible twist} and Lemma \ref{lem: H^1 twist}, we see that every gg acceptable $\bfb$-admissible LR pair $(\phi,\varepsilon)$ satisfying $\bft(\phi,\varepsilon,\tau(\phi))\cong(\gamma_0,a,[b_0])$ is contained in the image of the map $\eta_{\phi_0,\gamma_0}$, i.e. is related to $(\phi_0,\gamma_0)$ by twisting by $\Sha_G^\infty(\QQ,I_{\phi_0,\gamma_0}^\circ)$. Let $D(\phi_0,\gamma_0)\subset\Sha_G^\infty(\QQ,I_{\phi_0,\gamma_0}^\circ)$ be the subset of classes $e$ such that
\[
\bft(e\phi_0,\gamma_0,\tau(e\phi_0))\cong\bft(\phi_0,\gamma_0,\tau(\phi_0))=(\gamma_0,a,[b_0]),
\]
i.e. twists that do not change the Kottwitz parameter (up to isomorphism). Then the map
\[
D(\phi_0,\gamma_0)\hookrightarrow\Sha_G^\infty(\QQ,I_{\phi_0,\gamma_0}^\circ)\overset{\eta_{\phi_0,\gamma_0}}{\longrightarrow}\LRP^\gg_\sa/\text{conj.}
\]
is a surjection onto the set of conjugacy classes of LR pairs whose associated Kottwitz parameter is isomorphic to $(\gamma_0,a,[b_0])$. Note that it may not be a bijection, i.e. for each $e\in D(\phi_0,\gamma_0)$ we must account for the number of other elements of $D(\phi_0,\gamma_0)$ which give rise to the same conjugacy class of LR pairs, which number is $\#\{\text{fiber of $\eta_{\phi_0,\gamma_0}$ containing $e$}\}$. We proceed to count these sets, which is the final step en route to our final formula.

\begin{lemma}
\label{lem: count fibers of eta_phi0,gamma0}
In the setting of \ref{par: fibers of LRP -> KP}, we have
\[
\#\{\text{fiber of $\eta_{\phi_0,\gamma_0}$ containing $e$}\}=\frac{\lvert(I_{e\phi_0,\gamma_0}/I_{e\phi_0,\gamma_0}^\circ)(\QQ)\rvert}{[I_{e\phi_0,\gamma_0}(\QQ):I_{e\phi_0,\gamma_0}^\circ(\QQ)]}.
\]
\end{lemma}
Note that $(e\phi_0,\gamma_0)$ is a gg acceptable pair (since its corresponding Kottwitz parameter is acceptable), so we are in the situation of \ref{compatible inner twists} and we can equally well replace $I_{e\phi_0,\gamma_0}$ and $I_{e\phi_0,\gamma_0}^\circ$ by $G_{\gamma_0}$ and $G_{\gamma_0}^\circ$, respectively.
\begin{proof}
This is proven in the last paragraph of the proof of Lemma 3.7.6 in \cite{KSZ} (in their notation, $\iota_H(\varepsilon)=[H_\varepsilon(\QQ):H_\varepsilon^\circ(\QQ)]$ and $\overline{\iota}_H(\varepsilon)=\lvert(H_\varepsilon/H_\varepsilon^\circ)(\QQ)\rvert$).
\end{proof}

\begin{lemma}
\label{lem: count D(phi0,gamma0)}
In the setting of \ref{par: fibers of LRP -> KP}, we have
\[
\lvert D(\phi_0,\gamma_0)\rvert=\sum_{(a',[b_0'])}\lvert\Sha_G(\QQ,I_{\phi_0,\gamma_0}^\circ)\rvert
\]
where the sum runs over pairs $(a',[b_0'])$ such that $(\gamma_0,a',[b_0'])$ is an acceptable $\bfb$-admissible Kottwitz parameter with trivial Kottwitz invariant, and $\Sha_G$ is defined as in \ref{par: twisting admissible morphisms}.
\end{lemma}
As in Lemma \ref{lem: count fibers of eta_phi0,gamma0}, we can equally well replace $I_{\phi_0,\gamma_0}^\circ$ by $G_{\gamma_0}^\circ$.

\begin{proof}
Having fixed an LR pair $(\phi_0,\gamma_0)$ (and \emph{not} simply a conjugacy class), \cite[Prop 3.6.2(iii)]{KSZ} tells us that twisting by an element $e\in\Sha_G^\infty(\QQ,I_{\phi_0,\gamma_0}^\circ)$ results in a well-defined Kottwitz parameter $\bft(e\phi_0,\gamma_0,\tau(e\phi_0))$ (\emph{not} simply an isomorphism class). Thus we can write
\[
D(\phi_0,\gamma_0)=\coprod_{(a',[b_0'])}D_{(a',[b_0'])}
\]
where $(a',[b_0'])$ runs over all pairs for which $(\gamma_0,a',[b_0'])$ forms a Kottwitz parameter isomorphic to $(\gamma_0,a,[b_0])$, and $D_{(a',[b_0'])}\subset D(\phi_0,\gamma_0)$ is the subset of twists giving rise to the Kottwitz parameter $(\gamma_0,a',[b_0'])$.

If $D_{(a',[b_0'])}$ is non-empty, then by \cite[Prop 3.6.2]{KSZ} it must be a coset of $\Sha_G(\QQ,I_{\phi_0,\gamma_0}^\circ)$ inside $\Sha_G^\infty(\QQ,I_{\phi_0,\gamma_0}^\circ)$ (note that that proposition only uses their ``$p^n$-admissible'' hypothesis to show that $[b_0]$ is basic, which we know by our ``acceptable'' hypothesis). The proof that $D_{(a',[b_0'])}$ is indeed non-empty proceeds precisely as in the proof of \cite[Lemma 3.7.6]{KSZ} (where they call this set $D_i$), replacing their Lemma 3.6.2 with our Lemma \ref{lem: stably conjugate epsilon} and their Proposition 3.6.1 with our Lemma \ref{lemma: strong conjugacy}.
\end{proof}

Combining Lemmas \ref{lem: same T_xi^f}, \ref{lem: count fibers of eta_phi0,gamma0}, and \ref{lem: count D(phi0,gamma0)}, our point-counting formula is transformed into its final form.

\begin{theorem}
\label{point counting formula}
For any acceptable function $f\in C^\infty_c(G(\AA_f^p)\times J_b(\QQ_p))$, we have
\begin{align*}
\tr(f\mid H_c(\Ig_\Sigma,&\sL_\xi))
=
\sum_{\gamma_0\in\Sigma_{\Rell}(G)\phantom{W}}\sum_{(a,[b_0])\in\KP(\gamma_0)}
\\
&\frac{\lvert\Sha_G(\QQ,G_{\gamma_0}^\circ)\rvert}{\lvert(G_{\gamma_0}/G_{\gamma_0}^\circ)(\QQ)\rvert}
\vol\left(I^\circ_{\fc}(\QQ)\backslash I^\circ_{\fc}(\AA_f)\right) O_{\gamma\times\delta}^{G(\AA_f^p)\times J_b(\QQ_p)}(f)\tr(\xi(\gamma_0))
\end{align*}
where $I_\fc$ is the inner form of $G_{\gamma_0}^\circ$ associated to the Kottwitz parameter $\fc=(\gamma_0,a,[b_0])$ as in \ref{par: I_c}, and $\gamma,\delta$ are the elements belonging to the classical Kottwitz parameter $(\gamma_0,\gamma,\delta)$ associated to $\fc$ as in \ref{par: KP -> classical KP}.
\end{theorem}

\bibliographystyle{alpha}
\bibliography{igusa_v_arxiv}

\end{document}